\documentclass[11pt]{article}
\usepackage{epsfig}
\usepackage{amssymb}
\usepackage{amscd}
\usepackage{epsf}
\usepackage{amsfonts,amssymb,amsthm,amsmath}

\long\def\comment#1\endcomment{\relax}

\newcounter{subsubsubsection}
\newcounter{subsubsubsubsection}

\makeatletter
\@addtoreset{subsubsubsection}{subsubsection}
\@addtoreset{subsubsubsubsection}{subsubsubsection}

\makeatother

\newcommand{\sevafigc}[4]{\begin{figure}[h]\centerline{
 \epsfig{file=#1,width=#2,angle=#3}}
\bigskip\caption{#4}\end{figure}}

\DeclareMathSizes{11.1}{10}{8}{6}

\DeclareMathOperator{\Hom}{Hom}

\newtheorem*{theorem*}{Theorem}

\newtheorem{conjecture}{Conjecture}

\newtheorem*{lemma}{Lemma}

\newtheorem*{proposition}{Proposition}

\theoremstyle{remark}

\newtheorem*{remark}{Remark}

\newtheorem*{example}{Example}

\theoremstyle{definition}

\newtheorem*{definition}{Definition}

\DeclareMathOperator{\Id}{Id}

\DeclareMathOperator{\Alt}{Alt}

\newcommand{\mb}{{\bullet}}

\newcommand\End{\mathrm{End}}
\newcommand\Bi{\mathrm{Bialg}}

\newcommand\GS{\mathrm{GS}}

\newcommand\U{\mathcal U}

\newcommand{\codim}{\mathrm{codim}}

\newcommand{\Conf}{\mathrm{Conf}}

\newcommand{\g}{\mathfrak{g}}

\newcommand{\K}{\mathrm{K}}

\newcommand{\Hoch}{\mathrm{Hoch}}

\newcommand{\St}{\mathrm{St}}

\def\wtilde#1{\widetilde{#1}\vphantom{#1}}

\title{{\tt {\huge An explicit formula for the deformation quantization of Lie bialgebras}}}
\author{{\tt {\LARGE Boris Shoikhet}}}
\date{}

\begin{document}\maketitle

\begin{abstract}
{\tt A model of 3-dimensional topological quantum field theory is rigorously constructed.
The results are applied to an explicit formula for deformation quantization
of any finite-dimensional Lie bialgebra over the field of complex numbers.
This
gives an explicit construction of "quantum groups" from any Lie bialgebra,
which was proven without explicit formulas in [EK].}
\end{abstract}
\section*{\tt{Introduction}}
The most ideal goal of this paper would be a construction of an
$L_\infty$-structure on the Gerstenhaber-Schack complex $K_{GS}^\mb(A)$
of a (co)associative bialgebra $A$ and a proof of its formality for
$A=S(V^*)$, the free commutative cocommutative bialgebra of polynomial
functions on a finite-dimensional vector space $V$. Below in the
Introduction it is explained what part of this project is realized here.

First of all, let us recall the definitions. A (co)associative bialgebra $A$
is a vector space endowed with two operations, $*\colon A^{\otimes 2}\to A$
and $\Delta\colon A\to A^{\otimes 2}$, called the product and the coproduct,
correspondingly. These operations should obey the following 3 axioms:
\begin{itemize}
\item[(i)] $(a*b)*c=a*(b*c)$ for any $a,b,c\in A$ (the associativity),
\item[(ii)] $(\Delta\otimes 1)\circ\Delta (a)=(1\otimes\Delta)\circ\Delta
(a)$ for any $a\in A$ (the coassociativity),
\item[(iii)] $\Delta (a*b)=\Delta(a)*\Delta(b)$ for any $a,b\in A$ (the
compatibility)
\end{itemize}
(Here in the r.h.s. of (iii) the product $*$ on $A^{\otimes 2}$ is defined
as $(a_1\otimes a_2)*(b_1\otimes b_2)=(a_1*b_1)\otimes (a_2*b_2)$).

Notice that a (co)associative bialgebra could not have the (co)unit and the
antipode.

In the case of an associative algebra $A$ there is the well-known
construction of the Hochschild cohomological complex $\Hoch^\mb(A)$ and the
Gerstenhaber bracket on it which makes $\Hoch^\mb(A)$ a dg Lie algebra. This
dg Lie algebra plays a fundamental role in the deformation theory of
associative algebras. Namely, the deformation functor associated with this
dg Lie algebra describes the deformations of the algebra $A$ in the class of
associative algebras (more precisely, it describes the deformations of the
category of $A$-modules). Roughly, it means that for a cochain $\Psi\in
\Hom_{\mathbb{C}}(A^{\otimes 2},A)=\Hoch^1(A)$ the Maurer-Cartan equation
\begin{equation}\label{eq0001}
d\Psi+\frac12[\Psi,\Psi]=0
\end{equation}
in $\Hoch^\mb(A)$ is equivalent that the product $a\tilde{*}b=a*b+\Psi(a,b)$
is again associative (see, e.g. [K1] for details).

There is a complex which could be considered as analog of $\Hoch^\mb(A)$ in
the case of (co)associative bialgebras, namely, the Gerstenhaber-Schack
complex $K^\mb_{GS}(A)$ [GS]. Recall that
\begin{equation}\label{eq0002}
K^\mb_{GS}(A)=\bigoplus_{m,n\ge 1}\Hom_{\mathbb{C}}(A^{\otimes m},A^{\otimes
n})[-m-n+2]
\end{equation}
(in the agreement with the usual notations, $(L^\mb[s])^k=L^{s+k}$; in
particular, if $A$ has only degree 0, $\Hom_{\mathbb{C}}(A^{\otimes
m},A^{\otimes n})[-m-n+2]$ has degree $m+n-2$). In particular, if $A$ has
only degree 0, $K^1_{GS}(A)=\Hom(A^{\otimes 2},A)\oplus\Hom(A,A^{\otimes
2})$. In [GS], Gerstenhaber and Schack constructed a differential on
$K^\mb_{GS}(A)$ for any bialgebra $A$ such that the first cohomology
$\mathrm{H}^1(K^\mb_{GS}(A))$ is isomorphic to the infinitesimal deformations of the
(co)associative bialgebra structure on $A$.

Even the first attempt to construct a Lie algebra structure on
$K^\mb_{GS}(A)$ which would describe the global deformations of
(co)associative bialgebras via the Maurer-Cartan equation (\ref{eq0001})
fails. Indeed, consider $\Psi_1+\Delta_1\in \Hom(A^{\otimes
2},A)\oplus\Hom(A,A^{\otimes 2})=K^1_{GS}(A)$. The Maurer-Cartan equation
for any possible bracket is {\it quadratic} in $\Psi_1,\Delta_1$, while the
r.h.s. of the compatibility equation (iii) in the definition of
(co)associative bialgebra above is of the 4th degree in $\Psi_1,\Delta_1$.
Therefore, the best we can expect is the existence of an $L_\infty$-algebra
structure on $K^\mb_{GS}(A)$ with the first component equal to the
Gerstenhaber-Schack differential. Recall that an $L_\infty$-algebra on a
$\mathbb{Z}$-graded vector space $L^\mb$ is an odd vector field $Q$ of
degree +1 on the space $L^\mb[1]$ such that $Q^2=0$. There exists a
deformation theory associated with an $L_\infty$-algebra: the Maurer-Cartan
equation (\ref{eq0001}) is replaced by the equation
\begin{equation}\label{eq0003}
\{\gamma\in L^1\ {such\ that}\ Q|_{\gamma}=0\}
\end{equation}
It would be very nice to construct such an odd vector field $Q$ on
$K^\mb_{GS}(A)[1]$ such that the generalized Maurer-Cartan equation
(\ref{eq0003}) describes exactly the (co)associative bialgebras. This
problem still remains to be open.

If we would find an $L_\infty$-structure on $K^\mb_{GS}(A)$ and prove
explicitly the formality of it in the case when $A=S(V^*)$ (like it is done
in [K1] in the case of associative algebras), we immediately would get an
explicit construction of "quantum groups" from the infinitesimal datum--a
Lie bialgebra structure on $V$. The theorem that any Lie bialgebra can be
quantized was proven without explicit formulas in [EK].

In the present paper we find an explicit formula for the deformation
quantization of any Lie bialgebra $V$. For this, we find an analog of the
formality equation for $K^\mb_{GS}(S(V^*))$, but with unknown
$L_\infty$-structure. We derive this "formality" from the Stokes formula in
some rigorously defined in the paper "3-dimensional topological quantum
field theory". This "formality" equation is not the formality in the proper
sense, because we do not know the $L_\infty$ structure which is supposed
to be formal, but it is enough to quantize any Lie bialgebra $V$.

To construct the "3-dimensional topological quantum field theory", we
construct a compactification of the Kontsevich spaces $K(m,n)$ and of the
extended Kontsevich spaces $K(m,n;s)$. In the case of associative algebras,
the analog of $K(m,n)$ is the Stasheff space
\begin{equation}\label{eq0004}
\St_n=\{(p_1,\dots,p_n)\in\mathbb{R},p_i<p_j\ {for}\ i<j\}/G^{(2)}
\end{equation}
where the 2-dimensional group $G^{(2)}$ is the group of transformations
\begin{equation}\label{eq0005}
G^{(2)}=\{x\mapsto ax+b,\ a\in \mathbb{R}_+,b\in\mathbb{R}\}
\end{equation}

Actually, the whole deformation theory of associative algebras is contained
in the geometry of the Stasheff compactification $\overline{\St}_n$ of the
spaces $\St_n$. Namely, the chain operad $M^\mb=\bigoplus_{n\ge
2}C_\mb(\overline{\St}_n)$ in the realization by the Stasheff cells is free
and is a minimal model of the operad $\mathrm{Assoc}$ of associative
algebras. Then, the application of the Markl's construction from [M1] to
this minimal model gives exactly the Hochschild complex with the
Gerstenhaber bracket.

Unfortunately, we have not so nice description for our compactification
$\overline{K(m,n)}$. We hope to understand the geometry of this
compactification better in next papers. Now we have the description of all
strata of codimension 1. Just notice here that it is {\it not} the CROC
compactification from [Sh1] we suppose to be related with the unexisted
theory of non-commutative deformations, but rather a Stasheff-type
compactification. In particular, it is exactly the Stasheff compactification
when $m=1$ or $n=1$.

In the case of associative algebras, an approach alternative to the explicit
description of the chain operad as a minimal model (and which gives much more
strong results) is the Kontsevich approach [K1]. Kontsevich "extends" the
Stasheff space to a 2-dimensional configuration space $C_{m,n}$ and
constructs its compactification $\overline{C_{m,n}}$ extending the Stasheff
compactification. The space $C_{m,n}$ is the configuration space of $m$
non-coinciding points at the upper half-plane and of $n$ points at its
boundary $\mathbb{R}$. Kontsevich gives in [K1] from this compactification
a rigorous description of some particular case of the AKSZ model of
topological quantum field theory on an open disk. For this, he constructs an
appropriate "propagator" as a closed 1-form on $\overline{C_{2,0}}$. Within
this approach, the Stokes formula gives "some relation" even if we would not
know about the Gerstenhaber bracket. More precisely, the Gerstenhaber
bracket is the only one which makes this equation to be a formality
equation.

In this paper we follow the Kontsevich approach. We extend our
compactification of the 1-dimensional configuration space $K(m,n)$ to a
compactification of a 3-dimensional configuration space $K(m,n;s)$ (which is
the analog of the upper half-plane in the case of associative algebras). We
construct all ingredients of a model of the 3-dimensional topological quantum filed
theory from this compactification. We construct a propagator as a closed
2-form {\it with singularities} on $\overline{K(0,0;2)}$ which degenerates
to a closed 1-form in some limit. Then we associate some closed forms (with
singularities) to the admissible graphs and apply the Stokes formula.

In this way, we replace the Gerstenhaber-Schack complex $K^\mb_{GS}(A)$ to a
homotopically equivalent complex $\wtilde{K}^\mb_{GS}(A)$. The definition of
$\wtilde{K}^\mb_{GS}(A)$ itself depends on our compactification
$\overline{K(m,n)}$.

The Stokes formula gives an equation which we would like to interpret as the
"formality of $L_\infty$-morphism". The further research should shed some light to this
claim.

The good problems (if you like, it is a definition of a good problem) in
mathematics are valuable not only by themselves, but mostly by new ideas which
appear when one tries to solve them. The author is sure that deformation
theory of (co)associative bialgebras is a good problem. As far this problem
is still not solved, one can expect that it will grow many further ideas.
From this point of view the main idea of this paper is the introduction of
the complex $B^{\mb\mb}(m,n)$. Probably, the problem studied here is the
first problem where the introduction of it is really necessarily. We
formulated Conjectures 1,2 in Section 1 which formalize which properties we need from
this structure in a possible greater generality. We think that
technically the introduction of this bicomplex is the main new thing invented in
the paper.

\comment The goal of this paper is to construct
a deformation theory for (co)associative bialgebras. According to general
principles, it means that we are looking for a dg Lie algebra (or, more
generally, for an $L_\infty$ algebra) controlling the deformation theory of
a (co)associative bialgebra. In the case of the deformation theory of
associative algebras, such a dg Lie algebra controlling the deformations of
an associative algebra $A$, is the cohomological Hochschild complex of $A$
with the Gerstenhaber bracket. (More precisely, this Hochschild complex
controls the deformations of the category of $A$-modules).

First of all, recall that a (co)associative bialgebra is a vector space $A$
equipped with the maps $\star\colon A^{\otimes 2}\to A$ (the product) and
$\Delta\colon A\to A^{\otimes 2}$ (the coproduct). The product is supposed
to be associative and the coproduct is supposed to be coassociative.
Moreover, we suppose the following compatibility of them:
\begin{equation}\label{eq230.1}
\Delta(a\star b)=\Delta(a)\star\Delta(b)
\end{equation}
for any $a,b\in A$. (Here in the r.h.s. the product is the component product
in $A^{\otimes 2}$ defined as $(a\otimes b)\star(a_1\otimes b_1)=
(a\star a_1)\otimes (b\star b_1)$). Notice that we do not suppose the
existence of unit and counit in $A$.

Here we meet our first difficulty: the r.h.s. of (\ref{eq230.1}) is of the
4th degree and not quadratic. Recall that (little bit roughly) we associate
the deformation theory with a dg Lie algebra $\g^\mb$ as follows: we
consider the solutions of the Maurer-Cartan equation
\begin{equation}\label{eq230.2}
d\alpha+\frac12[\alpha,\alpha]=0
\end{equation}
for $\alpha\in \g^1$ modulo the action of the gauge group associated with
$\g^0$ on these solutions. (Because of possible divergences in the action of
the gauge group, we say instead of this direct construction that the
deformation functor is a functor from the category of the Artinian algebras
to the category of sets).

It is known that the Gerstenhaber-Schack complex [GS] associated with a
bialgebra $A$ is a deformation complex of the bialgebra structure on $A$.
It means that the first cohomology of this complex are isomorphic to the
infinitesimal deformations on $A$. To pass from the infinitesimal
deformations to the global ones, one needs to have an appropriate dg Lie algebra
structure on the Gerstenhaber-Schack complex (or, more generally, an
$L_\infty$-structure). Recall here that as a vector space, the
Gerstenhaber-Schack complex of $A$ is
\begin{equation}\label{eq230.3}
K^\mb_{GS}=\bigoplus_{m,n\ge 1}\Hom(A^{\otimes m},A^{\otimes n})[-m-n+2]
\end{equation}

In particular, in degree 1 we have: $K^1_{GS}=\Hom(A^{\otimes
2},A)\oplus\Hom(A,A^{\otimes 2})$. We could expect that for some dg Lie
algebra structure on $K^\mb_{GS}$ the Maurer-Cartan equation (\ref{eq230.2})
for the element $*_1\oplus \Delta_1\in\Hom(A^{\otimes
2},A)\oplus\Hom(A,A^{\otimes 2})$ means exactly that
$(*+*_1,\Delta+\Delta_1)$ defines a new (co)associative bialgebra structure
on $A$.

But it is impossible: because the r.h.s. of the equation (\ref{eq230.1}) is
of the 4th degree in $*_1$ and $\Delta_1$, while the Maurer-Cartan equation
(\ref{eq230.2}) is quadratic. It means that the best we could expect is to
have an $L_\infty$ algebra structure on $K^\mb_{GS}$ (which looks quite
complicated). This crucial observation was explained to the author by Boris
Tsygan about 3 years ago.

This idea is one source of the theory developed here in this paper. Another
source is the Kontsevich spaces $K(m,n)$. The reader can find the definition
of them in Section 2 of the paper. The original Kontsevich motivation when
he invented these spaces was the following:
the space $K(2,2)$ is the configuration space on two independent lines, we
have 2 points in each line modulo independent common shift on each line, and
modulo the following action of $\mathbb{R}^*_+$ on this space: for
$\lambda\in\mathbb{R}^*_+$, we dilatate the first line with the scale
$\lambda$ and the second with the coefficient $\lambda^{-1}$. Then $K(2,2)$
is a 1-dimensional space: we have an interval on the first line, an interval
on the second (we identify the intervals of the same length), and we
identify such two configurations with the same product of the lengths of
the two intervals. Therefore, the configuration has the only one module--the
product of the lengths of the intervals. Before compactification, it is
isomorphic to $\mathbb{R}_+$.

Now we compactify the space $K(2,2)$ to the closed interval. The two limit
configurations are shown in the Figure 1:
\sevafigc{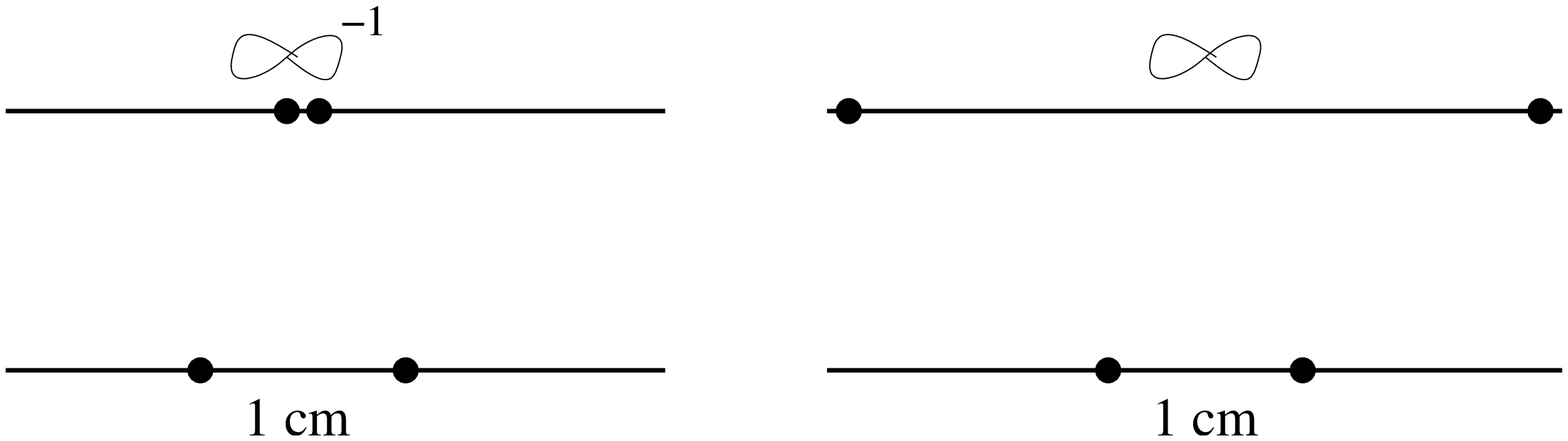}{100mm}{0}{The two limit points in $\overline{K(2,2)}$}
The Kontsevich's insight was that the left picture should give the left-hand
side of the compatibility equation (\ref{eq230.1}), while the right picture
should give the right-hand side of (\ref{eq230.1}). We say "should give"
having in mind  a
construction of the type of Kontsevich formality (see Sections 3-4).

Now let us outline the line of this paper. We think that the
spaces $K(m,n)$ and their 3-dimensional extensions $K(m,n;s)$ introduced by Maxim Kontsevich
are geniously
remarkable and should play the crucial role in the subject. One of the main
technical achievements of this paper is the construction of the
Stasheff-type compactification of these spaces. In the simplest case of
$K(0,0;2)$ we obtain a 3-dimensional analog of "The Eye" of Kontsevich [K1]. We
construct a propagator as a 2-form {\it with singularities} on the
3-dimensional Eye. In some limit we obtain after the integration a 1-form
which is also a part of the propagator. Then we construct from this
propagator differential forms associated with admissible graphs. Using all
these data, we construct the components $\\U_k\colon \wedge^k\g_1\to \g_2[1-k]$
where
$\g_1=\oplus_{m,n\ge 1}\Hom(\wedge^m(V), \wedge^n(V))[-m-n+2]$ and
$\g_2=\oplus_{m,n\ge 1}\Hom(A^{\otimes m},A^{\otimes n})[-m-n+2]$ where $V$
is a finite-dimensional vector space and $A=S(V^*)$ is the symmetric algebra
(considered as free commutative cocommutative bialgebra). We consider
$\g_1$ as the Poisson Lie algebra (the Poisson Lie algebra of functions on
$V[1]\oplus V^*[1]$ where the even symplectic structure is given by the canonical pairing of $V$ and $V^*$).
Then consider an element $\alpha\in
\Hom(\wedge^2(V),V)\oplus\Hom(V,\wedge^2(V))=\g_1^1$. The Maurer-Cartan
equation $[\alpha,\alpha]=0$ in the graded Lie algebra $\g_1$ is equivalent
to the condition that $\alpha$ defines a Lie bialgebra structure on $V$.

We found using the Stokes theorem the
family of relations on all $\U_k$ which we call the formality relations. We
do not know it is the formality of which structure on $\g_2$, it is a
subject for the further research. But to know these formality relations is
enough to quantize any Lie bialgebra structure on $V$ in the sense of the
Etingof-Kazhdan theorem. It means that we write down explicitly a
(co)associative bialgebra structure on $S(V^*)$ depending formally on two
parameters $\hbar_1$ and $\hbar_2$.

In Section 1 we introduce the concept of $\frac23$PROP.
It would be better to call this concept
$\frac12$PROP-$\frac12$CROC,
by for simplicity we use the name $\frac23$PROP. This concept is not
necessarily for the sequel but it is natural to describe the Stasheff-type
compactification of the spaces $K(m,n)$ using the concept of $\frac23$PROP.
We do it in Section 2. The first two Sections are inserted from [Sh2] (we
just added the last example in Section 2) for the convenience of the reader.
In Section 3 we introduce the extended Kontsevich space $K(m,n;s)$ and
construct its compactification $\overline{K(m,n;s)}$. In particular, we
construct the 3-dimensional Eye $\overline{K(0,0,2)}$ and construct the
non-homogeneous Propagator differential form on it. In Section 4 we
construct the components $\U_k$ and prove the formality relation. Also we
consider the particular case which gives an explicit formula for the
Etingof-Kazhdan quantization.
\endcomment
\comment
\section{\tt{The concept of $\frac23$PROP}}
\subsection{{\tt The definition}}
Here we define our main technical tool--$\frac23$PROPs. The name
$\frac23$PROPs indicates that this concept is a further generalization (or
simplification) of the concept of $\frac12$PROP due to Maxim Kontsevich
(see [K3], [MV]).

We define a $\frac23$PROP of vector spaces, a $\frac23$PROP of dg vector
spaces, of topological spaces,... can be defined analogously.
\begin{definition}
A pre-$\frac23$PROP of vector spaces consists of the following data:
\begin{itemize}
\item[(i)] a collection of vector spaces $F(m,n)$ defined for
$m\ge1,n\ge1,m+n\ge3$, with an action of symmetric groups
$\Sigma_m^\vee\times\Sigma_n$ on $F(m,n)$,
\item[(ii)] a collection of vector spaces $F^{1,1,\dots,1\ (n\ times)}_m$
defined for $n\ge 2,m\ge 2$, with an action of the symmetric group
$\Sigma_m^\vee$ on $F^{1,1,\dots,1\ (n\ times)}_m$,
\item[(iii)] a collection of vector spaces $F_{1,1,\dots,1\ (m\ times)}^n$
defined for $m\ge 2,n\ge 2$, with an action of the symmetric group
$\Sigma_n$ on $F_{1,1,\dots,1\ (m\ times)}^n$,
\item[(iv)] compositions $\circ_i\colon F(m,n)\otimes F(1,n_1)\to
F(m,n+n_1-1)$, $1\le i\le n$,
\item[(v)] compositions $_j\circ\colon F(m_1,1)\otimes F(m,n)\to
F(m+m_1-1,n)$, $1\le j\le m$,
\item[(vi)] compositions $\circledcirc_i\colon F_{1,1,\dots,1\ (m\
times)}^n\otimes F_{1,1,\dots,1\ (m\
times)}^{n_1}\to F_{1,1,\dots,1\ (m\
times)}^{n+n_1-1}$, $1\le i\le n$,
\item[(vii)] compositions $_j\circledcirc\colon F^{1,1,\dots,1\ (n\
times)}_{m_1}\otimes F^{1,1,\dots,1\ (n\
times)}_{m}\to F^{1,1,\dots,1\ (n\
times)}_{m+m_1-1}$, $1\le j\le m$,
\item[(viii)] compositions $\circledcirc\colon F^{1,1,\dots,1\ (n\
times)}_{m}\otimes F_{1,1,\dots,1\ (m\
times)}^n\to F(m,n)$,
\item[(ix)] all the compositions are equivariant with respect to the actions
of the symmetric groups.
\end{itemize}
This data should obey the following properties:
\begin{itemize}
\item[(1)] the composition $\circ_j\star\circ_i\colon F(m,n)\otimes
F(1,n_1)\otimes F(1,n_2)\to F(n,n+n_1+n_2-2)$ is associative for
$i\le j\le i+n_1$: in this case $\circ_{j-i}\star\circ_i\colon
F(m,n)\otimes
(F(1,n_1)\otimes F(1,n_2))\to F(n,n+n_1+n_2-2)$ coincides with
$\circ_j\star\circ_i \colon  (F(m,n)\otimes
F(1,n_1))\otimes F(1,n_2)\to F(n,n+n_1+n_2-2)$,
\item[(2)] in the notations of (1), if $j<i$ or $j>i+n_1$, we have
$\circ_j\star\circ_i\colon (F(m,n)\otimes
F(1,n_1))\otimes F(1,n_2)\to F(n,n+n_1+n_2-2)$ is equal to
$\circ_{i_1}\star\circ_{j_1}\colon (F(m,n)\otimes
F(1,n_2))\otimes F(1,n_1)\to F(n,n+n_1+n_2-2)$ where
$j_1=j, i_1=i+n_2$ if $j<i$, and $j_1=j-n_1, i_1=i$ if $j>i+n_1$ {\it (the commutativity)},
\item[(3)] the property analogous to (1) for $_j\circ$,
\item[(4)] the property analogous to (2) for $_j\circ$,
\item[(5)-(8)] the analogous properties for $\circledcirc_i$ and for
$_j\circledcirc$.
\end{itemize}
\end{definition}

Notice that this structure without $F^{1,1,\dots,1}_m$ and
$F_{1,1,\dots,1}^n$ is exactly the Kontsevich's $\frac12$PROP structure.
\begin{definition}
A pre-$\frac23$PROP is a $\frac23$PROP if the operations of $\circ$-type are
compatible with the operations of $\circledcirc$-type, as follows:

There are extra maps $\curlywedge_{m_1\to m_2}^n\colon F_{1,1,\dots,1\ (m_1\
times)}^n\to F_{1,1,\dots,1\ (m_2\ times)}^n$, $(m_1\le m_2)$ and maps
$\curlyvee_m^{n_1\to n_2}\colon F_m^{1,1,\dots,1\ (n_1\ times)}\to F_m^{1,1,\dots,1\ (n_2\
times)}$, $n_1\le n_2$, which are supposed to be equivariant with respect to
the actions of symmetric groups. We also suppose that these maps are isomorphisms. Then we have:
\begin{itemize}
\item[(A)]$\circ_i \star\circledcirc\colon(F^{1,1,\dots,1\ (n\ times)}_m\otimes F_{1,1,\dots,1\
(m\ times)}^n)\otimes F(1,n_1)\to F(m,n+n_1-1)$ is equal to
$\circledcirc\star\circledcirc_i\colon F_m^{1,1,\dots,1\ (n+n_1-1 \ times)}\otimes (F_{1,1,\dots,1\ (m\ times)}^n\otimes
F_{1,1,\dots,1\ (m\ times)}^{n_1})\to F(m,n+n_1-1)$. More precisely, let $\alpha\in F^{1,1,\dots,1\ (n\
times)}_m$, $\beta\in F_{1,1,\dots,1\
(m\ times)}^n)$, and $\gamma\in F(1,n_1)$. Then
\begin{equation}\label{eq23.01}
(\alpha\circledcirc\beta)\circ_i\gamma=  \curlyvee_m^{n\to
n+n_1-1}(\alpha)\circledcirc_i(\beta\circledcirc \curlywedge^{n_1}_{1\to
m}(\gamma))
\end{equation}
\item[(B)] the analogous compatibility with $_j\circ$.
\end{itemize}
\end{definition}
We can imagine what is a free $\frac23$PROP. It consists from all "free"
words of the following two forms:
\begin{equation}\label{eq23.1}
\dots\circ F(m_k,1)\circ\dots\circ F(m_1,1)\circ F(m,n)\circ
F(1,n_1)\circ\dots\circ F(1,n_\ell)\circ\dots
\end{equation}
and
\begin{multline}\label{eq23.2}
 F^{1,1,\dots,1\ (n_1+n_2+\dots-k+1 \ times)}_{m_k}\circledcirc
\dots\circledcirc F^{1,1,\dots,1\ (n_1+n_2+\dots-k+1 \
times)}_{m_1}\circledcirc\\
\circledcirc
F_{1,1,\dots,1 \ (m_1+m_2+\dots-l+1 \ times)}^{n_1}\circledcirc
 F_{1,1,\dots,1 \ (m_1+m_2+\dots -l+1\
times)}^{n_2}\circledcirc\dots\circledcirc F_{1,1,\dots,1 \ (m_1+m_2+\dots \ times)}^{n_l}
\end{multline}

We can draw these free elements as "two-sided trees", see Figure 2:
\sevafigc{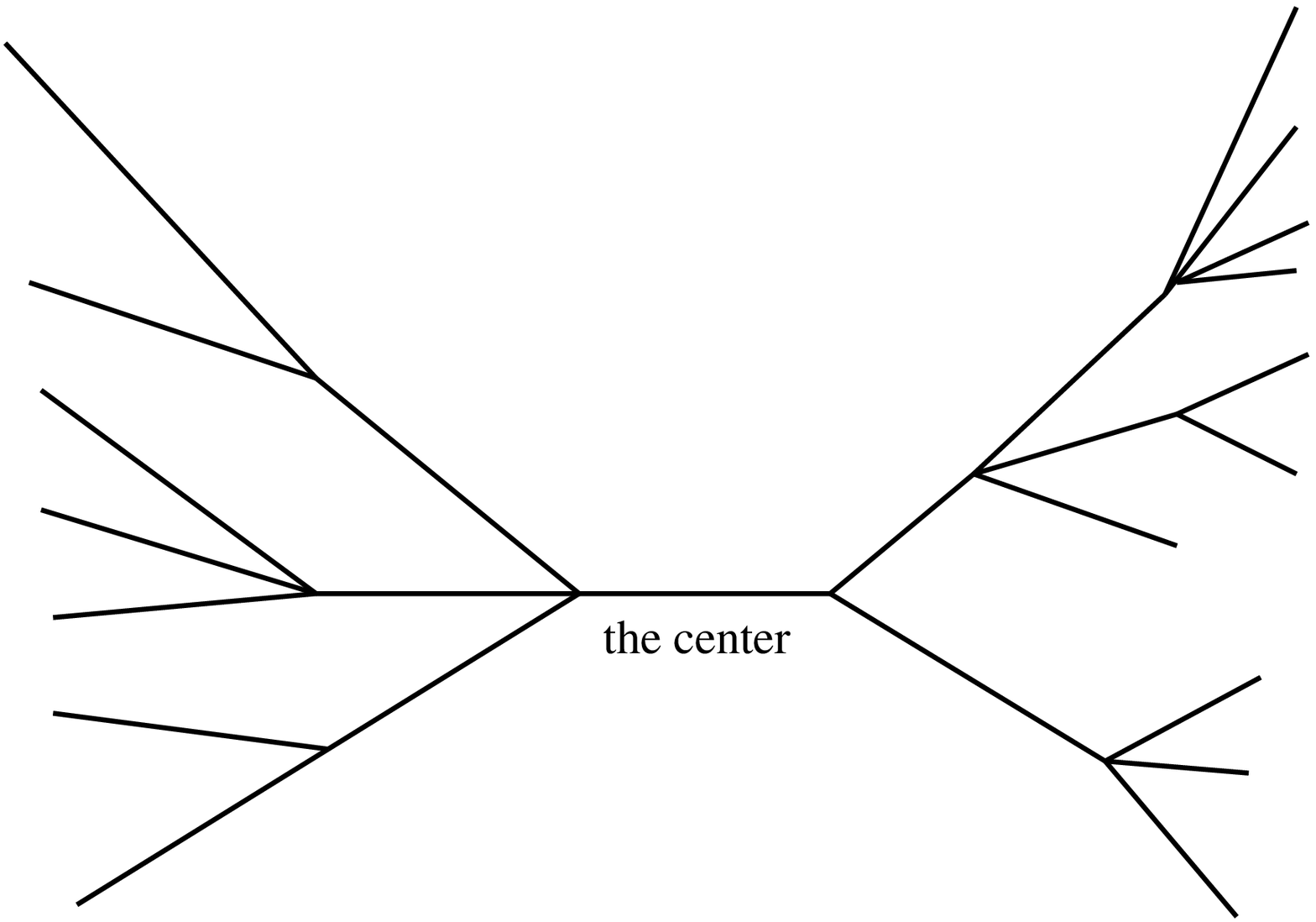}{90mm}{0}{A typical element of a free $\frac23$PROP}
The valences of the vertices are greater or equal than 3. We have trees of
"two different colors" according to the equations
(\ref{eq23.1})-(\ref{eq23.2}) above.

The definition of (pre-)$\frac23$PROP is motivated by the geometry of the
Kontsevich spaces $K(m,n)$ (see the next Section). The reader who is
interested in the origin of this definition can pass directly to Section 2.
We tried to construct a compactification of these spaces and to formalize
the operations among the strata. The advantage of $\frac23$PROPs is that the
chain complex of the compactification $\overline{K(m,n)}$ is a free dg
$\frac23$PROP. Its homology $\frac23$PROP is exactly the $\frac23$PROP $\Bi$
controlling the (co)associative bialgebras. Finally, any (co)associative
bialgebra structure on a vector space $V$ gives a map of the
pre-$\frac23$PROPs $\Bi\to\End(V)$.

\subsection{{\tt The pre-$\frac23$PROP $\End(V)$}}
Here we define the pre-$\frac23$PROP $\End(V)$ for a vector spaces $V$. We use
here the notations $\End(m,n)$, $\End_{1,1,\dots,1}^n$, and
$\End^{1,1,\dots,1}_m$.

We set:
\begin{equation}\label{eq23.3}
\begin{aligned}
\ & \End(m,n)=\Hom(V^{\otimes m},V^{\otimes n}),\\
& \End_{1,1,\dots,1\ (m\ times)}^n=(\Hom(V, V^{\otimes n}))^{\otimes m},\\
& \End^{1,1,\dots,1\ (n\ times)}_m=(\Hom(V^{\otimes m},V))^{\otimes n}
\end{aligned}
\end{equation}
We should now define the compositions $\circ_i,\ _j\circ,\
\circledcirc_i,\ _j\circledcirc$, and $\circledcirc$.

The case of the
composition $\circ_i\colon \End(m,n)\otimes\End(1,n_1)\to\End(m,n+n_1-1)$
can be schematically shown as follows (see Figure 3):
\sevafigc{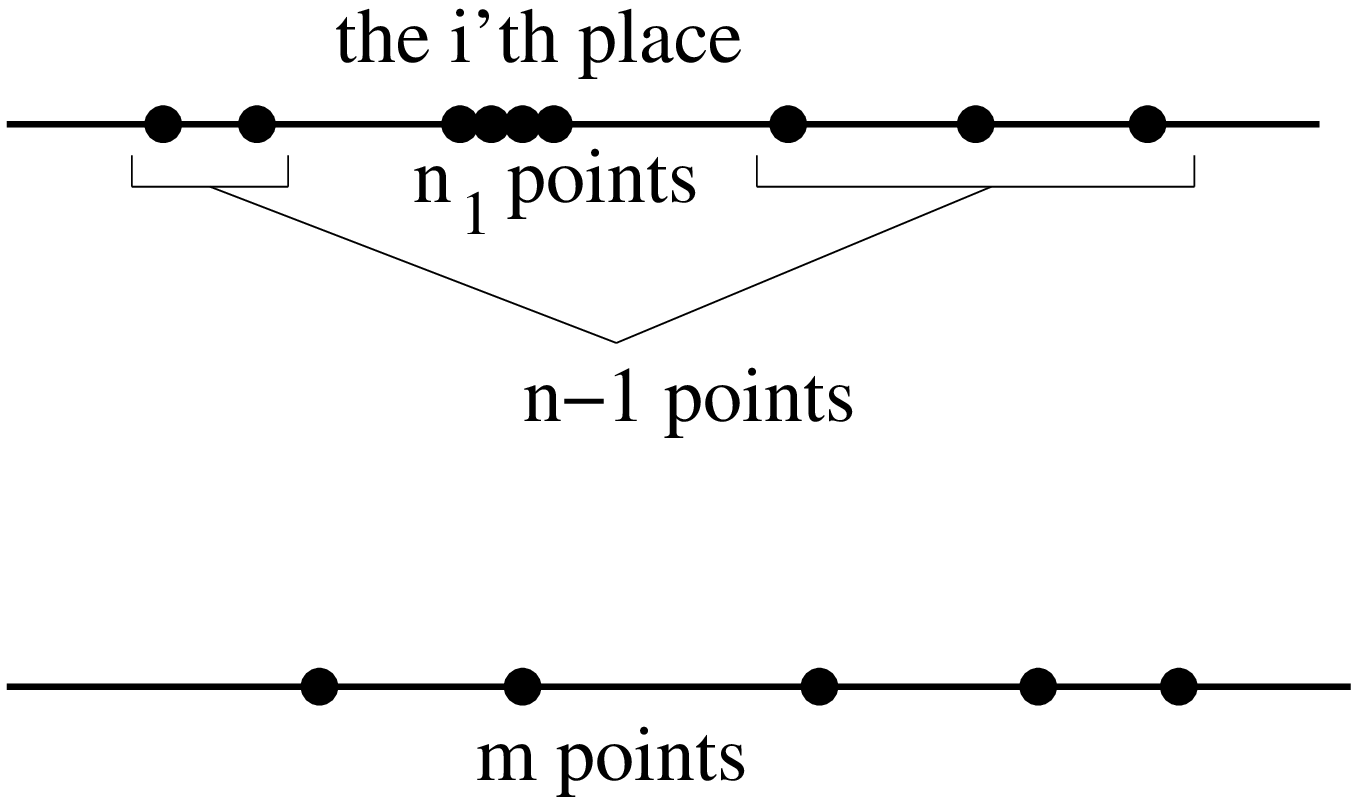}{80mm}{0}{The composition $\circ_i$}
Let $\Psi\in\Hom(V^{\otimes m},V^{\otimes n})$, $\Theta\in\Hom(V,V^{\otimes
n_1})$. Their composition $\Psi\circ_i\Theta\in\Hom(V^{\otimes m},
V^{\otimes n+n_1-1})$ is defined as
\begin{equation}\label{eq232.1}
\Psi\circ_i\Theta (v_1\otimes\dots\otimes v_m)=
(\Id\otimes\dots\otimes\Id\otimes\Theta\otimes\dots\otimes\Id)\circ\Psi(v_1\otimes\dots\otimes
v_m)
\end{equation}
(Here $\Theta$ stands at the $i$th place).

The picture for the composition
$_j\circ\End(m_1,1)\otimes\End(m,n)\to\End(m+m_1-1,n)$ is the following (see
Figure 4):
\sevafigc{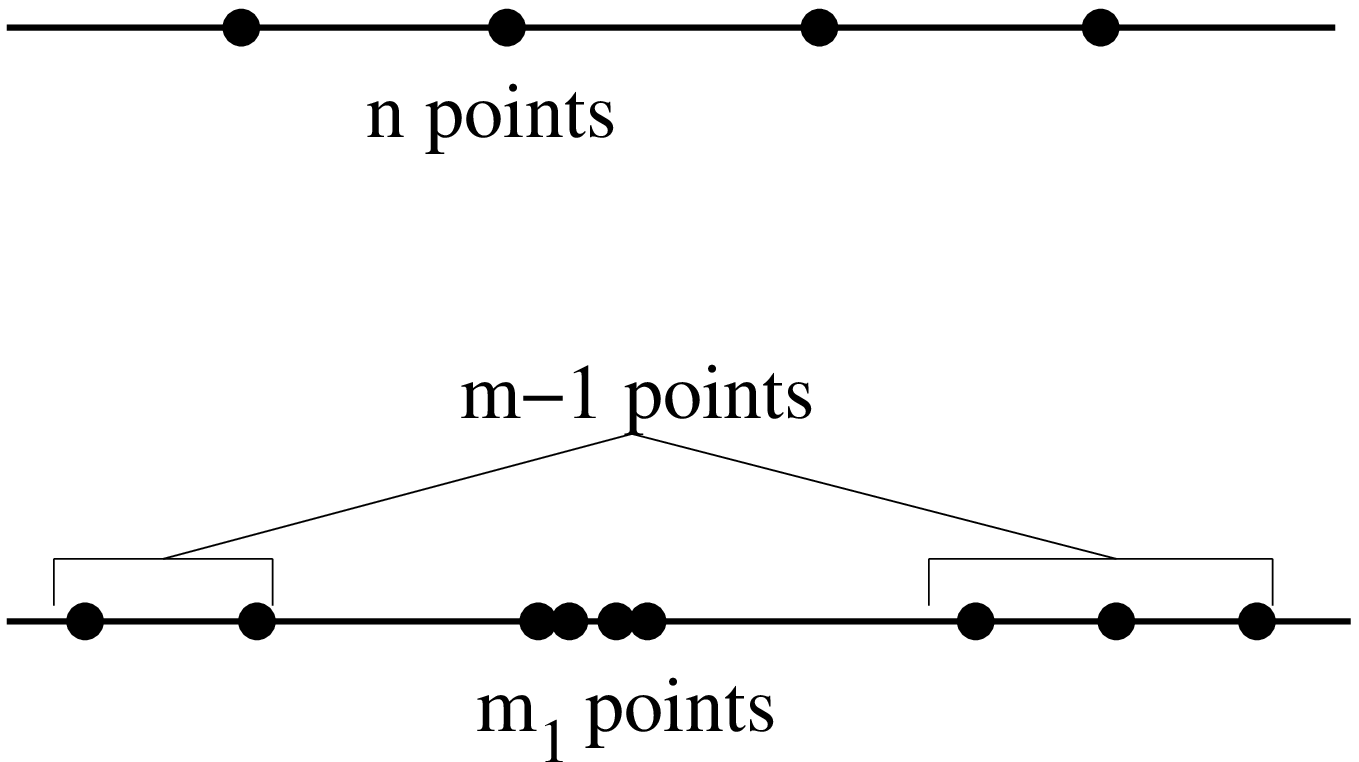}{80mm}{0}{The composition $_j\circ$}
For $\Theta\in \Hom(V^{\otimes m_1}, V)$ and $\Psi\in\Hom(V^{\otimes
m},V^{\otimes n})$, their composition
$\Theta(_j\circ)\Psi\in\Hom(V^{\otimes m+m_1-1},V^{\otimes n})$ is
\begin{multline}\label{eq232.2}
\Theta(_j\circ)\Psi(v_1\otimes\dots\otimes
v_{m+m_1-1})=\\
\Psi(v_1\otimes\dots\otimes
v_{j-1}\otimes\Theta(v_j\otimes\dots\otimes v_{j+m_1-1})\otimes
v_{j+m_1}\otimes\dots\otimes v_{m+m_1-1})
\end{multline}

We have the particular case of the product $\circ_i$ when $m=1$. Denote it
by $\circ_i^1$. By definition, the composition
$\circledcirc_i\colon\End_{1,1,\dots,1\ (m\ times}^n\otimes \End_{1,1,\dots,1\ (m\
times}^{n_1}\to \End_{1,1,\dots,1\ (m\ times}^{n+n_1-1}$ is the $m$th tensor
power of the composition $\circ_i^1$. Analogously we define $_j\circ^1$ and
the composition $_j\circledcirc\colon \End_{m_1}^{1,1,\dots,1\ (n\
times)}\otimes \End_{m}^{1,1,\dots,1\ (n\
times)}\to \End_{m_1+m-1}^{1,1,\dots,1\ (n\
times)}$ as the $n$th tensor power of the composition $_j\circ^1$.

It remains to define the composition
$\circledcirc$. We define the composition $\circledcirc\colon
(\Hom(V^{\otimes m},V))^{\otimes n} \otimes (\Hom(V,V^{\otimes n}))^{\otimes
m}\to\Hom(V^{\otimes m},V^{\otimes n})$. Suppose
$\Psi_1\otimes\dots\otimes\Psi_n\in(\Hom(V^{\otimes m},V))^{\otimes n}$, and
$\Theta_1\otimes\dots\otimes\Theta_m\in(\Hom(V,V^{\otimes n}))^{\otimes
m}$ We are going to define their composition
$(\Psi_1\otimes\dots\otimes\Psi_n)\circledcirc(\Theta_1\otimes\dots\otimes\Theta_m)\in\Hom(V^{\otimes
m},V^{\otimes n}$.
Denote
\begin{equation}\label{eq232.3}
F(v_1\otimes\dots\otimes
v_m)=\Theta_1(v_1)\otimes\dots\otimes\Theta_m(v_m)\in V^{\otimes mn}
\end{equation}
Next, we define a map $G\in\Hom(V^{\otimes mn}, V^{\otimes n})$ as follows:
\begin{multline}\label{eq232.4}
G(w^1_1\otimes\dots\otimes w_1^n\otimes w_2^1\otimes\dots\otimes
w_2^n\otimes\dots\otimes w_m^1\otimes\dots\otimes w_m^n)=\\
\Psi_1(w_1^1\otimes w_2^1\otimes
w_m^1)\otimes\Psi_2(w^2_1\otimes\dots\otimes w^2_m)\otimes\dots\otimes
\Psi_n(w^n_1\otimes\dots\otimes w^n_m)\in V^{\otimes n}
\end{multline}
Now we set
\begin{equation}\label{eq232.5}
(\Psi_1\otimes\dots\otimes\Psi_n)\circledcirc(\Theta_1\otimes\dots\otimes\Theta_m)=
(G\circ F)(v_1\otimes\dots\otimes v_m)
\end{equation}
It is clear that these compositions define a pre-$\frac23$PROP structure on
$\End(V)$.
\begin{remark}
M.~Markl communicated to the author that our composition $\circledcirc$ is a
particular case of his "fractions" composition [M2].
\end{remark}
\subsection{\tt{The $\frac23$PROP of (co)associative bialgebras $\Bi$}}
For a (pre-)$\frac23$PROP $F$ we define an $F$-algebra structure on a vector
space $V$ as a map of pre-$\frac23$PROPs $F\to\End(V)$. We are going to
construct now a $\frac23$PROP $\Bi$ such that a $\Bi$-algebra structure on
$V$ is exactly a (co)associative bialgebra structure on $V$.

Let $\Sigma_n$ be the symmetric group on $n$ points, and for a group $G$
denote by $G^\vee$ the dual group.

We can consider $\Bi$ as $\frac23$PROP of sets, or, if we like, as the
corresponding $\frac23$PROP of vector spaces (generated by these sets).
We here consider $\Bi$ as a $\frac23$PROP of sets. Later it will
appear also as the homology $\frac23$PROP of the topological $\frac23$PROP
$\overline{K(m,n)}$, then we consider it as the corresponding $\frac23$PROP
of vector spaces. This functor replaces the direct product $\times$ to the
tensor product $\otimes$.

We set:
\begin{equation}\label{eq232.6}
\begin{aligned}
\ &\Bi(m,n)=\Sigma_m^\vee\times\Sigma_n,\\
& \Bi_{1,1,\dots,1\ (m\ times)}^n=\Sigma_n,\\
&\Bi^{1,1,\dots,1\ (n\ times)}_m=\Sigma_m^\vee
\end{aligned}
\end{equation}
The $\frac23$PROP maps $\curlywedge^n_{m_1\to m}$ and $\curlyvee_m^{n_1\to n}$
are the identity maps. We define the compositions $\circ_i$, $_j\circ$, $\circledcirc_i$, $_j\circledcirc$ and $\circledcirc$ as follows:

Consider any of these compositions for the pre-$\frac23$PROP $\End(V)$, a
composition $\bigstar$.
Suppose that $\Psi\in\End_\alpha$ and $\Theta\in\End_\beta$ are its arguments. These compositions were
defined in the previous Subsection. There is the
$\Sigma_{i_1}^\vee\times\Sigma_{j_1}$-action on $\End_\alpha$ and the
$\Sigma_{i_2}^\vee\times\Sigma_{j_2}$-action on $\End_\beta$. Suppose
$\sigma_1^\vee\times\sigma_1\in\Sigma_{i_1}^\vee\times\Sigma_{j_1}$,
and $\sigma_2^\vee\times\sigma_2\in\Sigma_{i_2}^\vee\times\Sigma_{j_2}$. We
are going to define the composition $(\sigma_1^\vee\times
\sigma_1)\bigstar(\sigma_2^\vee\times\sigma_2)$ in $\Bi$. For this, consider
the composition $((\sigma_1^\vee\times
\sigma_1)\Psi)\bigstar((\sigma_2^\vee\times\sigma_2)\Theta$. It is clear that
it is equal to the action of some $\sigma\in\Bi$ on the product of $\Psi$
and $\Theta$ in $\End$:
\begin{equation}\label{eq232.7}
((\sigma_1^\vee\times
\sigma_1)\Psi)\bigstar((\sigma_2^\vee\times\sigma_2)\Theta=\sigma(\Psi\bigstar\Theta)
\end{equation}
The last equation holds for any $\Psi$ and $\Theta$ in the corresponding
components of $\End$, that is, $\sigma$ does not depend on the choice of
$|psi$ and $\Theta$. We define the composition $(\sigma_1^\vee\times
\sigma_1)\bigstar(\sigma_2^\vee\times\sigma_2)$ as $\sigma$.

It clear that this definition is correct, and in this way we define a
$\frac23$PROP $\Bi$.
\begin{lemma}
A map $\phi\colon \Bi\to\End(V)$ of pre-$\frac23$PROPs is the same that a
(co)associative bialgebra structure on $V$.
\begin{proof}
First, let $V$ be a (co)associative bialgebra with the product $\star$ and
the coproduct $\Delta$. We define a map of pre-$\frac23$PROPs
$\phi_{\star,\Delta}\colon\Bi\to\End(V)$. Let
$\sigma^\vee\times\sigma\in\Bi(m,n)$. We put
$\phi_{\star,\Delta}(\sigma^\vee\times\sigma)(v_1\otimes\dots\otimes v_m)=
\Delta^{n-1}\circ \star^{m-1}(v_1\otimes\dots\otimes v_m)$. Here in the
formula $\Delta^n$ and $\star^m$ are the composition powers of the coproduct
and of the product, correspondingly. Because of the (co)associativity, these
powers are well-defined. Next, for $\sigma\in\Bi_{1,1,\dots,1\
(m\times}^n$ we set $\phi_{\star,\Delta}(\sigma)\in(\Hom(V,V^{\otimes
n}))^{\otimes m}$ is the $m$th tensor power of the map
$v\mapsto \sigma(\Delta^{n-1}(v))$. Analogously, using the product, we
define $\phi_{\star,\Delta}$ on $\Bi^{1,1,\dots,1}_m$. Now we explain why
without the compatibility (\ref{eq230.1}) this definition would be
incorrect.

Consider many identity permutations: $\Id_2\in \Sigma_2=\Bi(1,2)$,
$\Id^2\in \Sigma_2^\vee=\Bi(2,1)$, $\Id_2^{1,1}\in \Bi_{1,1}^2$, and
$\Id^2_{1,1}\in\Bi^{1,1}_2$. Than we have the following identity in $\Bi$:
\begin{equation}\label{eq232.8}
\Id^2\circ\Id_2=\Id^2_{1,1}\circledcirc\Id_2^{1,1}=\Id^\vee\times\Id\in\Bi(2,2)
\end{equation}
It is clear that this identity in $\Bi$ follow some identity in the images
of these elements $\Id_2,\ \Id^2,\ \Id_2^{1,1},\ \Id^2_{1,1}$ by the map
$\phi_{\star,\Delta}$ of pre-$\frac23$PROPs. The reader can easily verify
that this identity is exactly the compatibility (\ref{eq230.1}) in a
(co)associative bialgebra. One can prove also that if the compatibility
holds, the definition of $\phi_{\star,\Delta}$ is correct.

Vice versa, suppose we have a map $\phi\colon\Bi\to\End(V)$ of
pre-$\frac23$PROPs. Denote $a\star b:=\phi(\Id^2)(a\otimes b)$ and
$\Delta(a):=\phi(\Id_2)(a)$. The compatibility follows from (\ref{eq232.8}).
The reader can easily find analogous identities in $\Bi$ which imply the
associativity of $\star$ and the coassociativity of $\Delta$.
\end{proof}
\end{lemma}
\endcomment
\section{\tt{The Kontsevich spaces $K(m,n)$, their Stasheff-type
compactification, and the homotopical Gerstenhaber-Schack complex.}}
\subsection{\tt{The compactification}}
First of all, recall the definition of the spaces $\K(m,n)$ due to Maxim
Kontsevich (see also [Sh]). We show in the sequel that these spaces and its
compactification introduced below play a crucial role in the deformation
theory of (co)associative bialgebras.

First define the space $\Conf (m,n)$. By definition, $m,n\ge 1$, $m+n\ge 3$,
and
\begin{multline}\label{eq1.1}
\Conf (m,n)=\{p_1,\dots, p_m\in \mathbb{R}^{(1)}, p_i<p_j\ \  for\ \  i<j;\\
q_1,\dots,q_n\in\mathbb{R}^{(2)}, q_i<q_j\ \ for\ \ i<j\}
\end{multline}
Here we denote by $\mathbb{R}^{(1)}$ and by $\mathbb{R}^{(2)}$ two different
copies of a real line $\mathbb{R}$.

Next, define a 3-dimensional group $G^3$ acting on $\Conf(m,n)$. This group
is a semidirect product $G^3=\mathbb{R}^2\ltimes\mathbb{R}_+$ (here
$\mathbb{R}_+=\{x\in\mathbb{R}, x>0\}$) with the following group law:
\begin{equation}\label{eq1.2}
(a^{\prime},b^{\prime},\lambda ^{\prime})\circ(a,b,\lambda)=
(\lambda ^{\prime} a+a^{\prime},(\lambda ^{\prime})^{-1} b+b^{\prime},\lambda\lambda ^{\prime} )
\end{equation}
where $a,b,a^\prime ,b^\prime\in\mathbb{R}, \lambda,\lambda^\prime\in\mathbb{R}_+$.
This group acts on the space $\Conf(m,n)$ as
\begin{multline}\label{eq1.3}
(a,b,\lambda)\cdot (p_1,\dots,p_m;q_1,\dots,q_n)=
(\lambda p_1+a,\dots,\lambda
p_m+a;\lambda^{-1}q_1+b,\dots,\lambda^{-1}q_n+b)
\end{multline}
In other words, we have two independent shifts on $\mathbb{R}^{(1)}$ and
$\mathbb{R}^{(2)}$ (by $a$ and $b$), and $\mathbb{R}_+$ dilatates
$\mathbb{R}^{(1)}$ by $\lambda$ and dilatates $\mathbb{R}^{(2)}$ by $\lambda
^{-1}$.

In our conditions $m,n\ge 1, m+n\ge 3$, the group $G^3$ acts on $\Conf(m,n)$
freely. Denote by $\K(m,n)$ the quotient-space. It is a smooth manifold of
dimension $m+n-3$.

We will need also the spaces  $\K_{m_1,\dots,m_{\ell_1}}^{n_1,\dots,n_{\ell_2}}$  introduced below.
Recall here our definition of the space
$\K_{m_1,\dots,m_{\ell_1}}^{n_1,\dots,n_{\ell_2}}$ (generalizing the
Kontsevich space $\K(m,n)$) from [Sh]:

Fist define the space
$\Conf^{m_1,\dots,m_{\ell_1}}_{n_1,\dots,n_{\ell_2}}$. By definition,
\begin{multline}\label{eq11.4}
\Conf^{m_1,\dots,m_{\ell_1}}_{n_1,\dots,n_{\ell_2}}=\\
\{p^1_1,\dots,p^1_{m_1}\in\mathbb{R}^{(1,1)},
p^2_1,\dots,p^2_{m_2}\in\mathbb{R}^{(1,2)},\dots,
p^{\ell_1}_1,\dots,p^{\ell_1}_{m_{\ell_1}}\in\mathbb{R}^{(1,\ell_1)};\\
q^1_1,\dots,q^1_{n_1}\in\mathbb{R}^{(2,1)},
q^2_1,\dots,q^2_{n_2}\in\mathbb{R}^{(2,2)}\dots,
q^{\ell_2}_1,\dots,q^{\ell_2}_{n_{\ell_2}}\in\mathbb{R}^{(2,\ell_2)}|\\
p^j_{i_1}<p^j_{i_2}\ \ for\ \ i_1<i_2;q^j_{i_1}<q^j_{i_2}\ \ for\ \
i_1<i_2\}
\end{multline}
Here $\mathbb{R}^{(i,j)}$ are copies of the real line $\mathbb{R}$.
Now we have an $\ell_1+\ell_2+1$-dimensional group $G^{\ell_1,\ell_2,1}$
acting on $\Conf^{m_1,\dots,m_{\ell_1}}_{n_1,\dots,n_{\ell_2}}$.
It contains $\ell_1+\ell_2$ independent shifts
$$
p_i^j\mapsto p_i^j+a_j, i=1,\dots, m_j, a_j\in\mathbb{R};
q_i^j\mapsto q_i^j+b_j, i=1,\dots, n_j, b_j\in\mathbb{R}
$$
and {\it one } dilatation
$$
p_i^j\mapsto \lambda\cdot p_i^j\ \ for\ \ all\ \ i,j;q_i^j\mapsto
\lambda^{-1}\cdot q_i^j\ \ for\ \ all\ \ i,j.
$$
This group is isomorphic to $\mathbb{R}^{\ell_1+\ell_2}\ltimes\mathbb{R}_+$.
We say that the lines
$\mathbb{R}^{(1,1)},\mathbb{R}^{(1,2)},\dots,\mathbb{R}^{(1,\ell_1)}$
(corresponding to the factor $\lambda$) are the lines of the first type, and
the lines $\mathbb{R}^{(2,1)},\mathbb{R}^{(2,2)},\dots,\mathbb{R}^{(2,\ell_2)}$
(corresponding to the factor $\lambda^{-1}$) are the lines of the second
type.

Denote
\begin{equation}\label{eq11.5}
\K^{m_1,\dots,m_{\ell_1}}_{n_1,\dots,n_{\ell_2}}=\Conf
^{m_1,\dots,m_{\ell_1}}_{n_1,\dots,n_{\ell_2}}/G^{\ell_1,\ell_2,1}
\end{equation}

We are going to construct a compactication $\overline{K(m,n)}$
of the space $K(m,n)$. First consider the simplest examples.
\subsubsection*{{\tt Example}}
Consider the case when $(m,n)$ is $(1,3)$ or $(3,1)$. Then the spaces
$K(1,3)$ or $K(3,1)$ are the Stasheff spaces with 3 points. We compactify
them as the corresponding Stasheff polyhedra to the unit closed interval.
More generally, the space $\overline{K(m,1)}$ (or $\overline{K(1,m)}$
in the compactification defined below is the Stasheff polyhedron $\St_m$
with $m$ points.
\subsubsection*{{\tt Example}}
Let $m=n=2$. Then the space $\K(2,2)$ is 1-dimensional. It is easy to see
that $(p_2-p_1)\cdot (q_2-q_1)$ is preserved by the action of $G^3$, and it
is the only invariant of the $G^3$-action on $\K(2,2)$. Therefore, $\K(2,2)\simeq
\mathbb{R}_+$. There are two "limit" configurations: $(p_2-p_1)\cdot
(q_2-q_1)\rightarrow 0$ and $(p_2-p_1)\cdot (q_2-q_1)\rightarrow \infty$.
Therefore, the compactification $\overline{\K(2,2)}\simeq [0,1]$. See Figure
1 below:
\sevafigc{2o3-0.eps}{100mm}{0}{The two limit points in $\overline{K(2,2)}$}
\bigskip

\begin{remark}
The Kontsevich's insight when he introduced the spaces $K(m,n)$ was that
the left picture in Figure 1 "should give" the l.h.s of the compatibility
equation in the definition of (co)associative bialgebras, $\Delta(a*b)$, and
the right picture in Figure 1 should give the r.h.s. $\Delta(a)*\Delta(b)$
of the compatibility equation. It will be more clear in the next Section
when we define the compactification of the extended Kontsevich spaces
$\overline{K(m,n;s)}$.
\end{remark}
Now having in mind the two previous examples, we define the compactification
in the general case.

Consider the set $\Sigma_{1,3}$ consisting from all possible samples of 1
point among the $m$ points at the first line in $K(m,n)$ and of 3 points
among the $n$ points at the second line in $K(m,n)$. For any
$\sigma\in\Sigma_{1,3}$ we have a map
\begin{equation}\label{eq0010}
r_{\sigma}\colon K(m,n)\to K(1,3)
\end{equation}
Analogously, we define the sets $\Sigma_{3,1}$ and $\Sigma_{2,2}$ as the
sets of all samples of 3 points at the first line and 1 point at the second
line, and of all samples of 2 points at each line, correspondingly.
For any $\sigma\in \Sigma_{3,1}$ we have the map $r_{\sigma}\colon K(m,n)\to
K(3,1)$, and for any $\sigma\in \Sigma_{2,2}$ we have the map
$r_{\sigma}\colon K(m,n)\to K(2,2)$.

Now consider the map which is the product of all $r_\sigma$ over all
possible $\sigma$. It is a map
\begin{multline}\label{eq0011}
r=\prod_{\sigma\in\Sigma_{1,3}\bigsqcup\Sigma_{3,1}\bigsqcup\Sigma_{2,2}}r_{\sigma}\colon K(m,n)\to\prod_{\sigma\in
\Sigma_{1,3}}K(1,3)_{\sigma}\times\prod_{\sigma\in\Sigma_{3,1}}K(3,1)_{\sigma}
\times\prod_{\sigma\in\Sigma_{2,2}}K(2,2)_{\sigma}
\end{multline}
(Here the lower index $\sigma$ in the r.h.s. just indicates the copy of the
space associated with $\sigma$).

It is clear that the map $r$ is {\it an imbedding}. Now we can compactify
the image, $\prod_{\sigma\in
\Sigma_{1,3}}K(1,3)_{\sigma}\times\prod_{\sigma\in\Sigma_{3,1}}K(3,1)_{\sigma}
\times\prod_{\sigma\in\Sigma_{2,2}}K(2,2)_{\sigma}$, compactifying each
factor as in the Examples above. We get an imbedding $\overline{r}$ of
$K(m,n)$ to a compact space (actually, a product of closed intervals). We
define the compactification $\overline{K(m,n)}$ as the closure of the image
of the map $\overline{r}$. The following lemma shows that the obtained
compactification is of the Stasheff-type:

\begin{lemma}
\begin{itemize}
\item[(i)] The space $\overline{K(m,n)}$ is a manifold with corners,
\item[(ii)] any stratum of codimension 1 has the form as is shown in the
Figure~2 below:
Here in the Figure~2 $m_0$ points at the first line move close to each
other in the scale $\frac1\infty$ and other points at this line are in
finite distance from each other, and $n_1$ "external" points at the second
line move infinitely far from each other with the scale $\infty$ (here this
$\infty$ and $\infty$ in the fraction $\frac1\infty$ above are "the same"),
and other "internal" $n_0=n-n_1$ points are in finite distance from each
other. Notice that after the application of the transformation
$(0,0,\infty)\in"\overline{G^(3)}"$ the picture on the first line after the
transformation becomes as the picture on the second line before the
transformation, and wise versa. This stratum of codimension 1 is canonically
isomorphic to $K_{m_1+1}^{1,1,\dots,1,n_0,1,\dots,1}\times
K^{n_1+1}_{1,1,\dots,1,m_0,1,\dots,1}$.
\end{itemize}
\sevafigc{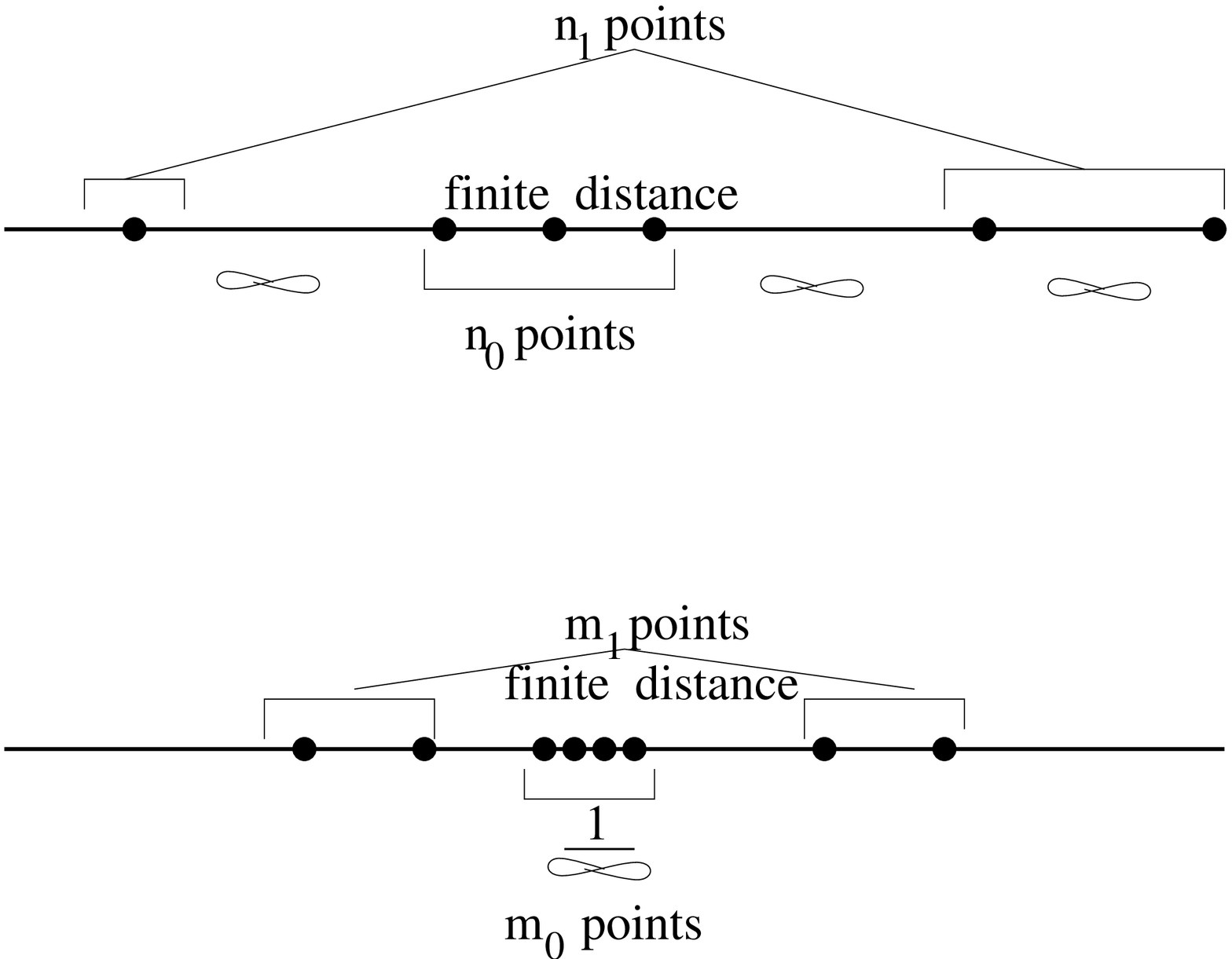}{100mm}{0}{A typical stratum of codimension 1}
\begin{proof}
We call {\it a 4-point ratio} the image of any map $r_\sigma$ defined above.
It corresponds a real number to a sample of 1 point on the first line and
3 points on the second line, or to a sample of 3 points on the first line and
1 point on the second line, or to 2 points on each line.

Consider the maximal open strata in $\overline{K(m,n)}$ which is, by
definition, the image of $K(m,n)$ under the imbedding $r$ (before taking
the closure). It is clear that the minimal number of 4-point ratios we
should know to reconstruct the configuration which belongs to the maximal
open stratum is exactly $m+n-3$, the dimension of this open stratum.
Moreover, we can introduce coordinates on it using the 4-point ratios.

Now define all $k$-dimensional strata in $\overline{K(m,n)}$ as the limit
configurations which we can reconstruct from $k$ of 4-point
ratios ($k\le m+n-3$) and can not reconstruct by any $l<k$ 4-point ratios.
It is in general nonconnected space, the connected components of which we
call the strata of dimension $k$. We can use the 4-point ratios to introduce
coordinates on it. It is clear that in this way we get a {\it manifold with
corners}.

One easily sees that the strata like the stratum drawn in Figure~2 are
uniquely defined by $m+n-4$ 4-point ratios and not less of them, and they
exhaust all such strata.
\end{proof}
\end{lemma}

Now we are going to attach to the strata of codimension 1 some operations on
the Gerstenhaber-Schack space $K^\mb_{GS}(A)$ for any $A$.
\subsection{\tt The strata of codimension 1 and operations on $K^\mb_{GS}(A)$}
We first define the operations in a bit bigger generality than we will
really need in the sequel. These operations is a particular case of Markl's
fractions in [M2], and we independently introduced in the context of CROCs
in [Sh1].

Let $V$ be a vector space.
Suppose we have
$$
\Psi_1\in\Hom(V^{\otimes \ell_1},V^{\otimes N_1}),\Psi_2\in\Hom(V^{\otimes \ell_1},V^{\otimes
N_2}),\dots,\Psi_{\ell_2}\in\Hom(V^{\otimes \ell_1},V^{\otimes
N_{\ell_2}})
$$
and
$$
\Theta_1\in\Hom(V^{\otimes
M_1},V^{\otimes \ell_2}),\Theta_2\in\Hom(V^{\otimes
M_2},V^{\otimes \ell_2}),\dots,\Theta_{\ell_1}\in\Hom(V^{\otimes
M_{\ell_1}},V^{\otimes \ell_2}),
$$
we are going to define their composition which belongs to
$\Hom(V^{\otimes {M_1+\dots+M_{\ell_1}}},V^{\otimes {N_1+\dots+N_{\ell_2}}})$. Denote
$m=M_1+\dots+M_{\ell_1}$, $n=N_1+\dots+M_{\ell_2}$. The construction is
as follows:

First define
\begin{multline}\label{eqc2.6}
F(v_1\otimes\dots\otimes v_m)\colon =
\Theta_1(v_1\otimes\dots\otimes
v_{M_1})\bigotimes\Theta_2(v_{M_1+1}\otimes\dots\otimes
v_{M_1+M_2})\bigotimes\dots\\ \bigotimes
\Theta_{\ell_2}(v_{M_1+\dots+M_{\ell_1-1}+1}\otimes\dots\otimes v_{{M_1}+\dots+{M_{\ell_1}}})
\in V^{\otimes \ell_1\ell_2}
\end{multline}
Now we apply $\{\Psi_i\}$'s to this element in $V^{\otimes \ell_1\ell_2}$:
we define an element $G\colon V^{\otimes \ell_1\ell_2}\to V^{\otimes n}$ as follows:
\begin{multline}\label{eqc2.7}
G(v_1\otimes v_2\otimes\dots\otimes v_{\ell_1\ell_2}):=
\Psi_1(v_1\otimes v_{\ell_2+1}\otimes\dots\otimes
v_{\ell_2(\ell_1-1)+1})\bigotimes\\
\Psi_2(v_2\otimes v_{\ell_2+2}\otimes\dots\otimes
v_{\ell_2(\ell_1-1)+2})\bigotimes\dots
\bigotimes\Psi_{\ell_2}(v_{\ell_2}\otimes v_{2\ell_2}\otimes\dots\otimes
v_{\ell_2\ell_1})\in V^{\otimes n}.
\end{multline}
Define now
\begin{equation}\label{eqc2.8}
Q(v_1\otimes \dots\otimes v_m):=G\circ
F(v_1\otimes\dots\otimes v_m)\in V^{\otimes n}
\end{equation}

By definition, the element $Q$ is the composition
$\frac{\Psi_1\Psi_2\dots\Psi_{\ell_2}}{\Theta_1\Theta_2\dots\Theta_{\ell_1}}\in
\Hom(V^{\otimes m},V^{\otimes n})$.

In [Sh1], we introduced these compositions to associate operations on
$K^\mb_{GS}$ with the strata of codimension 1 in $\overline{K(m,n)}_C$ in
the CROC compactification. In the case of the compactification introduced
here, we need only the following particular case of the construction above.
Consider $\ell_1=m_1+1$, $\ell_2=n_1+1$ in the notations of Lemma above, and
all from $M_i$'s are equal to 1 except a one which is equal to $m_0$, and
also all $N_j$'s are equal to 1 except one which is equal to $n_0$. The
reader could remind the Figure ... above.

Our (unrealized) goal is to construct an $L_\infty$ structure on
$K^\mb_{GS}(V)$ (we can consider $V$ as a (co)associative bialgebra with 0
product and 0 coproduct). According to this goal, we would like to consider
the operations
$\frac{\Psi_1\Psi_2\dots\Psi_{\ell_2}}{\Theta_1\Theta_2\dots\Theta_{\ell_1}}$
as candidates for components of an $L_\infty$ structure on $K_{GS}^\mb(V)$.
The first thing we should check is the grading condition, and even here we
runs to troubles. Indeed,
an $L_\infty$ operation is an operation of the form
$\wedge^{m_1+n_1+2}\g^\mb\to\g^\mb[2-(m_1+1)-(n_1+1)]$ where $\g^\mb=K_{GS}^\mb(V)$
with the
natural grading
\begin{equation}\label{eq0015}
\deg\Hom(V^{\otimes m},V^{\otimes n})=m+n-2
\end{equation}
Denote $A=(\sum_i\deg\Psi_i+\sum_j\deg\Theta_j)-m_1-n_1$ and $B=m_0+m_1+n_0+n_1-2$. We
would have the correct $L_\infty$ gradings iff $A=B$. Denote the defect
$A-B$ by $D$. We have: $A=2m_1n_1+m_0+n_0-2$, $D=2m_1n_1-m_1-n_1$. We see
that the defect $D=0$ iff $(m_1,n_1)=(0,0)$ or $(1,1)$.

In the next Subsection we introduce a complex $\wtilde{K}^\mb_{GS}(V)$
quasi-isomorphic to the graded space $K^\mb_{GS}(V)$ and redefine the
operations
$\frac{\Psi_1\Psi_2\dots\Psi_{\ell_2}}{\Theta_1\Theta_2\dots\Theta_{\ell_1}}$
on it such that they will be compatible with the $L_\infty$ rule of degrees.
\subsection{\tt The homotopical Gerstenhaber-Schack space}
For each $(m,n),m,n\in\mathbb{Z}_+,m+n\ge 3$ introduce the bicomplex
$B^{\mb\mb}(m,n)$ as follows:

Denote by $\Omega_+(\sigma)$ the de Rham complex of smooth differential
forms on a stratum $\sigma$ of the manifold with corners $\overline{K(m,n)}$
which can be continued to the {\it normalization} of the pair $(\sigma,\overline{\sigma})$.
Here by the normalization we mean a compact space $\overline{\sigma}_{norm}$
with a projection $\overline{\sigma}_{norm}\to\overline{\sigma}$ which
is a 1-1 map over the open stratum $\sigma$ and which separates the points
of the boundary $\overline{\sigma}\setminus\sigma$ which are limits of points which are far from each other on $\sigma$.
For example, consider the stratification of the circle drawn in Figure~3
below.
\sevafigc{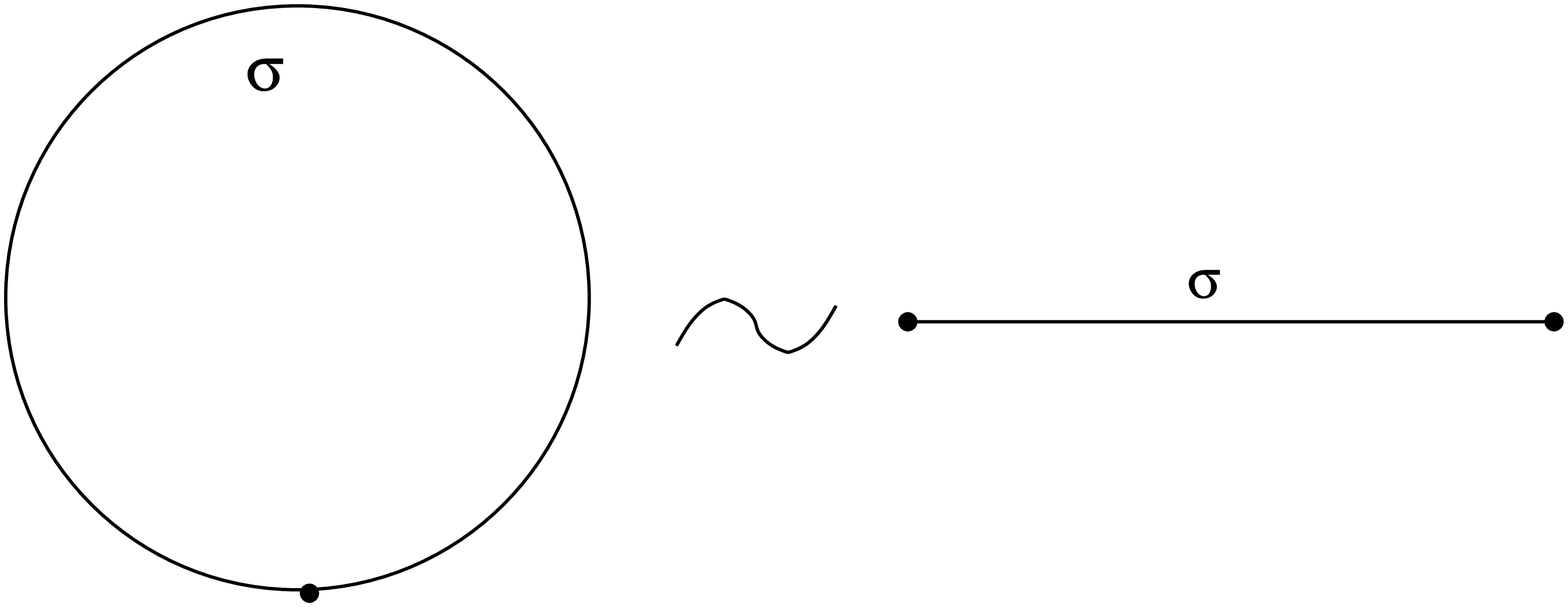}{100mm}{0}{A stratification of the circle and its normalization}
Here we have the two strata: the stratum $\sigma$ and the point. The
normalization of $\overline{\sigma}$ here is {\it not} the circle
$\overline{\sigma}$, but an interval. Namely, we separate the two boundary
points of $\overline{\sigma}$. In particular, $\Omega_+(\sigma)$ are {\it
not} the forms on the circle, but the forms on the closed interval, they
could have two different limits on the two boundary points. In the {\it
chain} differential below we take the difference of the restrictions of a
function (a 1-form) in the two boundary points in Figure~3.

We set:
\begin{equation}\label{eq0020}
B^{\mb\mb}(m,n)=\bigoplus_{\sigma\in\overline{K(m,n)}}\Omega_+(\sigma)[\dim\sigma]
\end{equation}
We consider the two differentials of degree +1 on $B^{\mb\mb}(m,n)$: the
first is the de Rham differential, and the second $-\partial$ is the chain differential
in $\overline{K(m,n)}$ acting on $\sigma$'s {\it with the opposite sign}. The total cohomology of these
two differentials is clearly 1-dimensional in the grading 0 and 0 otherwise.
Denote by $[\omega]_\sigma$ a differential form $\omega$ on a stratum
$\sigma$ considered as an element in $B^{\mb\mb}(m,n)$. It is clear that
$\deg [\omega]_\sigma=\deg\omega-\dim\sigma$. In particular, the bicomplex
$B^{\mb\mb}(m,n)$ is $\mathbb{Z}_{\le 0}$-graded.
We denote just by $[\omega]$ a differential form on the top degree open
stratum.

Introduce
\begin{equation}\label{eq0021}
\wtilde{K}^\mb_{GS}(V)=\Hom(V,V)[0]\oplus\bigoplus_{m+n\ge 3}\Hom(V^{\otimes
m},V^{\otimes n})\otimes B^{\mb\mb}(m,n)
\end{equation}
It is clear that $\wtilde{K}^\mb_{GS}(V)$ is quasi-isomorphic to the graded
space $K^\mb_{GS}(V)$.

Now we are going to define "the right" operations $\frac{\Psi_1\Psi_2\dots\Psi_{\ell_2}}{\Theta_1\Theta_2\dots\Theta_{\ell_1}}$
on $\wtilde{K}_{GS}^\mb(V)$. First we consider the forms on the open top
dimensional stratum.

Let $\wtilde{\Psi}_i=\Psi_i\otimes[\omega_i]$ and
$\wtilde{\Theta}_j=\Theta_j\otimes[\omega^\prime_j]$ where $\Psi_i$'s and
$\Theta_j$'s are as above. We consider the total degree:
$\deg\wtilde{\Psi}_i=\deg\Psi_i+\deg[\omega_i]$, and analogously for
$\wtilde{\Theta}_j$. Here all forms $\omega_i$ except the one are
differential forms on the top degree stratum in $\overline{K(m_1+1,1)}$ and
the remaining one is a form on the top dimensional stratum in
$\overline{K(m_1+1,n_0)}$. Analogously, the all forms $\omega_j^\prime$'s
except one are forms on the top dimensional stratum in
$\overline{K(1,n_1+1)}$, and the remaining one is a form on the top
dimensional stratum of $\overline{K(m_0,n_1+1)}$.

The operation $\frac{\Psi_1\Psi_2\dots\Psi_{\ell_2}}{\Theta_1\Theta_2\dots\Theta_{\ell_1}}$
itself is associated with a stratum $\sigma$ of codimension 1 in
$\overline{K(m,n)}$. The stratum $\sigma$ is {\it canonically} isomorphic to
$K_{m_1+1}^{1,1,\dots,1,n_0,1,\dots,1}\times
K^{n_1+1}_{1,1,\dots,1,m_0,1,\dots,1}$ (here in the upper index at the first
factor there are $n_1$ of 1's and $n_0$ is at the $s$'th place, and in the
second factor in the lower index there are $m_1$ of 1's and $m_0$ is at the
$s^\prime$'th place.
First of all, construct a differential form $\omega_{tot}$ on this stratum
$\sigma$ of degree $\sum_i\deg\omega_i+\sum_j\deg\omega_j$, starting from
the forms $\omega_i$'s and $\omega_j^\prime$'s. The construction is as
follows:

Consider the projections $p_1,\dots,p_{s-1},p_{s+1},\dots,p_{n_1+1}\colon
K_{m_1+1}^{1,1,\dots,1,n_0,1,\dots,1}\to K(m_1+1,1)$ which is just the
forgetting map which forgets the all upper points except the $i$'th for
$p_i$, and the forgetting map $p_s\colon
K_{m_1+1}^{1,1,\dots,1,n_0,1,\dots,1}\to K(m_1+1,n_0)$. Then take the wedge
product $\Omega_1=\wedge_{i=1}^{n_1+1}p_i^*(\omega_i)$ on $K_{m_1+1}^{1,1,\dots,1,n_0,1,\dots,1}$.
Analogously, define
the projections $p_j^\prime$ and define the form
$\Omega_2=\wedge_{j=1}^{m_1+1}p_j^{\prime *}\omega_j^\prime$ on
$K^{n_1+1}_{1,1,\dots,1,m_0,1,\dots,1}$. Then consider the form
$\omega:=\Omega_1\boxtimes\Omega_2$ on the product $K_{m_1+1}^{1,1,\dots,1,n_0,1,\dots,1}\times
K^{n_1+1}_{1,1,\dots,1,m_0,1,\dots,1}$. This is, by definition, the form
$\omega$ defined on the open stratum $\sigma=K_{m_1+1}^{1,1,\dots,1,n_0,1,\dots,1}\times
K^{n_1+1}_{1,1,\dots,1,m_0,1,\dots,1}$ of codimension 1 in
$\overline{K(m,n)}$. It is clear, that
$\deg\omega=\sum_i\deg\omega_i+\sum_j\deg\omega_j$.

Now define the composition
$\frac{\wtilde{\Psi}_1\wtilde{\Psi}_2\dots\wtilde{\Psi}_{n_1+1}}{\wtilde{\Theta}_1\wtilde{\Theta}_2\dots
\wtilde{\Theta}_{m_1+1}}$ as
\begin{equation}\label{eq0022}
\frac{\wtilde{\Psi}_1\wtilde{\Psi}_2\dots\wtilde{\Psi}_{n_1+1}}{\wtilde{\Theta}_1\wtilde{\Theta}_2\dots
\wtilde{\Theta}_{m_1+1}}:=\frac{\Psi_1\Psi_2\dots\Psi_{n_1+1}}{\Theta_1\Theta_2\dots\Theta_{m_1+1}}\otimes[\omega]_\sigma
\end{equation}

It is clear that in this definition the $L_\infty$ degree condition holds.
Indeed,
we need to prove that
\begin{multline}\label{eq0023}
\deg\wtilde{\Psi}_1+\dots+\deg\wtilde{\Psi}_{n_1+1}+\deg
\wtilde{\Theta}_1+\dots+\deg\wtilde{\Theta}_{m_1+1}+(2-(m_1+n_1+2))=\\
\deg\frac{\Psi_1\Psi_2\dots\Psi_{n_1+1}}{\Theta_1\Theta_2\dots\Theta_{m_1+1}}
+\deg[\omega]_\sigma
\end{multline}

We have:
$\deg\wtilde{\Psi}_i=\deg\Psi_i+\deg\omega_i-(m_1+1+1-3)=1+\deg\omega_i$ for $i\ne s$
and $\deg\wtilde{\Psi}_s=\deg\Psi_s+\deg\omega_s-(m_1+1+n_0-3)=1+\deg\omega_s$.
Analogously,
$\deg\wtilde{\Theta}_j=\deg\Theta_j+\deg\omega_j^\prime-(n_1+1+1-3)=1+\deg\omega^\prime_j$ for
$j\ne s^\prime$, and
$\deg\wtilde{\Theta}_{s^\prime}=\deg\Theta_{s^\prime}+\deg\omega_{s^\prime}^\prime-
(n_1+1+m_0-3)=1+\deg\omega^\prime_{s^{\prime}}$. Also,
$\deg[\omega]_\sigma=\sum_{i=1}^{n_1+1}\deg\omega_i+\sum_{j=1}^{m_1+1}\deg\omega_j^\prime-(m+n-4)$
where $m=m_0+m_1$ and $n=n_0+n_1$ (here $m+n-4$ is the dimension of the
stratum $\sigma$), and
$\deg\frac{\Psi_1\Psi_2\dots\Psi_{n_1+1}}{\Theta_1\Theta_2\dots\Theta_{m_1+1}}=m+n-2$.
Then (\ref{eq0023}) reduces to
\begin{multline}\label{eq0024}
m_1+n_1+2-(m_1+n_1)+\sum_i\deg\omega_i+\sum_j\deg\omega_j=\\
m+n-2+\sum_i\deg\omega_i+\sum_j\deg\omega_j-(m+n-4)
\end{multline}
which surely holds.

\begin{remark}
It is clear now the origin of our problem with the "naive" definition of the
composition on the level of $K^\mb_{GS}(V)$. Namely, our $\Psi$ and $\Theta$
there is identical to $\Psi\otimes[\omega_0]$ and
$\Theta\otimes[\omega_0^\prime]$ where $\omega_0$ and $\omega_0^\prime$
are {\it top degree differential forms}. Then the resulting form $\omega$ on
the stratum $\sigma$ is 0 by dimensional reasons except few simplest cases
(when $(m_1,n_1)=(0,0)$ or $(1,1)$). Then, the naively defined operation
should be 0 from our point of view except these 2 cases. On the other hand,
when the degrees of the forms $\omega_i$'s and $\omega_j^\prime$'s are
sufficiently small, the answer could be non-zero.
\end{remark}

It remains to consider the general case, when the forms $\omega_i$'s and
$\omega_j^\prime$'s are defined on strata of any codimension (the previous
case is the case of codimension 0).
\subsubsection{The case of arbitrary strata}
First of all, define a compactification of the space
$K_{m_1,\dots,m_{\ell_1}}^{n_1,\dots, n_{\ell_2}}$.
For each $1\le i\le \ell_1$ and $1\le j\le \ell_2$ we have a projection
$p_{ij}\colon K_{m_1,\dots,m_{\ell_1}}^{n_1,\dots, n_{\ell_2}}\to
K(m_i,n_j)$ (we suppose that $m_i+n_j\ge 3$.
Consider the map $\prod_{ij}p_{ij}\colon K_{m_1,\dots,m_{\ell_1}}^{n_1,\dots,
n_{\ell_2}}\to \prod_{ij}K(m_i,n_j)$. Clearly it is an embedding. Next,
consider another embedding $i\colon \prod_{ij}K(m_i,n_j)\to
\prod_{ij}\overline{K(m_i,n_j)}$. The composition $i\circ(\prod_{ij}p_{ij})$
is again an embedding. It embeds the open space $K_{m_1,\dots,m_{\ell_1}}^{n_1,\dots,
n_{\ell_2}}$ to a compact space. Define the compactification $\overline{K_{m_1,\dots,m_{\ell_1}}^{n_1,\dots,
n_{\ell_2}}}$ as the closure of the image of this embedding. It is a
manifold with corners with the natural stratification, defined analogously
with the stratification of $\overline{K(m,n)}$.

We give the following definition:
\begin{definition}
Let $M_1,M_2$ be two manifold with corners, and let $s\colon M_1\to M_2$ be
a continuous map. We say that $s$ is a map of manifold with corners iff
the preimage of each stratum in $M_2$ is a manifold with corners, and the
manifold with corners $M_1$ is glued from these preimages (as a manifold
with corners). Moreover, we demand that the restriction of $s$ to the
preimage of each stratum $\sigma$ in $M_2$ is a trivial bundle, the fiber of which is
a manifold with corners, and the total space of this bundle, as a manifold
with corners, is the product of the open space $\sigma$ with this fiber.
\end{definition}
\begin{conjecture}
\begin{itemize}
\item[(i)] the map $p_{ij}$ defines a map of manifold with corners
$\overline{p_{ij}}\colon \overline{K_{m_1,\dots,m_{\ell_1}}^{n_1,\dots,
n_{\ell_2}}}\to\overline{K(m_i,n_j)}$,
\item[(ii)] the imbedding $i_{\sigma}\colon
K_{m_1+1}^{1,\dots,1,n_0,1\dots,1}\times
K_{1,\dots,1,m_0,1,\dots,1}^{n_1+1}\to K(m,n)$ ($m=m_0+m_1, n=n_0+n_1$) is
continued to a map of manifold with corners $\overline{i_\sigma}\colon
\overline{K_{m_1+1}^{1,\dots,1,n_0,1\dots,1}}\times
\overline{K_{1,\dots,1,m_0,1,\dots,1}^{n_1+1}}\to\overline{K(m,n)}$.
\end{itemize}
\end{conjecture}

Define now the compositions $\frac{\wtilde{\Psi}_1\wtilde{\Psi}_2\dots\wtilde{\Psi}_{n_1+1}}{\wtilde{\Theta}_1\wtilde{\Theta}_2\dots
\wtilde{\Theta}_{m_1+1}}$ in the general case as follows:

Consider a stratum $\sigma$ of some codimension in
$\overline{K_{m_1+1}^{1,\dots,1,n_0,1\dots,1}}$. All images
$\overline{p_i}\colon \overline{K_{m_1+1}^{1,\dots,1,n_0,1\dots,1}}\to
\overline{K(m_1+1,N_i)}$ (here $N_i=1$ for all $i$ except $i=s$ for which
$N_s=n_0$) are single strata in $\overline{K(m_1+1,N_i)}$ according to
Conjecture 1(i). Define the image of $\sigma$ with respect to the map $\overline{p_i}$
by $\sigma_i$. As well, consider a stratum $\sigma^\prime$ of some
codimension in $\overline{K_{1,\dots,1,m_0,1,\dots,1}^{n_1+1}}$. Consider the
maps $\overline{p_j}^\prime\colon
\overline{K_{1,\dots,1,m_0,1,\dots,1}^{n_1+1}}\to \overline{K(M_j,n_1+1)}$.
Denote the stratum in $\overline{K(M_j, n_1+1)}$ which is the image of
$\sigma^\prime$ with respect to the map $\overline{p_j}^\prime$ by
$\sigma_j^\prime$. We start with forms  $\omega_i$ (of some degrees) on the
strata $\sigma_i$ in $\overline{K(m_1+1, N_i)}$, and with forms
$\omega_j^\prime$ (of some degrees) on the strata $\sigma^\prime_j$ in
$\overline{K(M_j, n_1+1)}$. Let
$\wtilde{\Psi}_i=\Psi_i\otimes[\omega_i]_{\sigma_i}$ and
$\wtilde{\Theta_j}=\Theta_j\otimes[\omega^\prime_j]_{\sigma^\prime_j}$ where
$\Psi_i$'s and $\Theta_j$'s are as above.

We define the composition $\frac{\wtilde{\Psi}_1\wtilde{\Psi}_2\dots\wtilde{\Psi}_{n_1+1}}{\wtilde{\Theta}_1\wtilde{\Theta}_2\dots
\wtilde{\Theta}_{m_1+1}}$  as
\begin{equation}\label{eq2004_1}
\frac{\wtilde{\Psi}_1\wtilde{\Psi}_2\dots\wtilde{\Psi}_{n_1+1}}{\wtilde{\Theta}_1\wtilde{\Theta}_2\dots
\wtilde{\Theta}_{m_1+1}}:=\frac{\Psi_1\Psi_2\dots\Psi_{n_1+1}}{\Theta_1\Theta_2\dots\Theta_{m_1+1}}\otimes[\Omega]_\Sigma
\end{equation}
where the stratum $\Sigma$ in $\overline{K(m,n)}$ and a form $\omega$ on it
are defined as follows:

Define the form $\Omega_1$ on stratum $\sigma$ in $\overline{K_{m_1+1}^{1,\dots,1,n_0,1\dots,1}}$
as $\Omega_1=\wedge_{i=1}^{n_1+1}\overline{p_i}^*\omega_i$, and define the
form $\Omega_2$ on stratum $\sigma^\prime$ in
$\overline{K_{1,\dots,1,m_0,1,\dots,1}^{n_1+1}}$ as
$\Omega_2=\wedge_{j=1}^{m_1+1}\overline{p_j}^{\prime *}\omega^\prime_j$.
Now consider the form $\Omega_1\boxtimes\Omega_2$ on the stratum
$\sigma\times\sigma^\prime$ in $\overline{K_{m_1+1}^{1,\dots,1,n_0,1\dots,1}}\times
\overline{K_{1,\dots,1,m_0,1,\dots,1}^{n_1+1}}$. According to Conjecture~1(ii), the image of the stratum $\sigma\times\sigma^\prime$ in $\overline{K_{m_1+1}^{1,\dots,1,n_0,1\dots,1}}\times
\overline{K_{1,\dots,1,m_0,1,\dots,1}^{n_1+1}}$
is a single stratum $\Sigma$ in $\overline{K(m,n)}$. Define now $\Sigma$ as
the image (it can be not the isomorphic image) of
$\sigma\times\sigma^\prime$ with respect to the map
$\overline{i}_\sigma\colon\overline{K_{m_1+1}^{1,\dots,1,n_0,1\dots,1}}\times
\overline{K_{1,\dots,1,m_0,1,\dots,1}^{n_1+1}}\to\overline{K(m,n)}$, and
define the form $\Omega$ on $\overline{K(m,n)}$ as the direct image (the
integration along the fiber) of the form $\Omega_1\boxtimes\Omega_2$ with
respect to the restriction of $\overline{i}_\sigma$ to the stratum
$\sigma\times\sigma^\prime$ (which is a trivial bundle).
{\it We consider only the case when
$\sum_i\codim_{\overline{K(m_1+1,N_i)}}\sigma_i+\sum_j\codim_{\overline{K(M_j,n_1+1)}}\sigma_j^\prime=
\codim_{\overline{K_{m_1+1}^{1,\dots,1,n_0,1\dots,1}}\times
\overline{K_{1,\dots,1,m_0,1,\dots,1}^{n_1+1}}}\sigma\times\sigma^\prime$.}
In this case, the formula above for  $\frac{\wtilde{\Psi}_1\wtilde{\Psi}_2\dots\wtilde{\Psi}_{n_1+1}}{\wtilde{\Theta}_1\wtilde{\Theta}_2\dots
\wtilde{\Theta}_{m_1+1}}$
obeys the $L_\infty$ grading condition. The proof is straightforward.

We define the composition $\frac{\wtilde{\Psi}_1\wtilde{\Psi}_2\dots\wtilde{\Psi}_{n_1+1}}{\wtilde{\Theta}_1\wtilde{\Theta}_2\dots
\wtilde{\Theta}_{m_1+1}}$ as 0 in all other cases.

The following conjecture has the crucial meaning, at the moment we have no
ways to prove it.
\begin{conjecture}
The compositions  $\frac{\wtilde{\Psi}_1\wtilde{\Psi}_2\dots\wtilde{\Psi}_{n_1+1}}{\wtilde{\Theta}_1\wtilde{\Theta}_2\dots
\wtilde{\Theta}_{m_1+1}}$ defined above are compatible with the (total) differential
in $B^{\mb\mb}(m,n)$'s.
\end{conjecture}

\comment
\begin{equation}\label{eq2323.1}
\begin{aligned}
\ &F(m,n)=\overline{K(m,n)}\\
&F_{1,1,\dots,1}^n=\overline{K_{1,1,\dots,1}^n}\\
&F_m^{1,1,\dots,1}=\overline{K^{1,1,\dots,1}_m}
\end{aligned}
\end{equation}
The compactifications will be defined below.

First of all, let us describe all strata of codimension 1 in
$\overline{K(m,n)}$. There are boundary strata of codimension 1 of two
different types. The first two strata are shown in Figures 3 and 4. In the
picture in Figure 3 $n_1$ points on the upper line move infinitely close to each other,
and the "scale" of this infinitely small number is irrelevant (we have in
mind here the CROC compactification from [Sh] where this scale is relevant),
$2\le n_1\le n$. In Figure 4 $m_1$ points on the lower line move close to
each other, $2\le m_1\le m$. The remaining stratum of codimension 1 in
$\overline{K(m,n)}$ (there is the only one such stratum) is shown in Figure
5:
\sevafigc{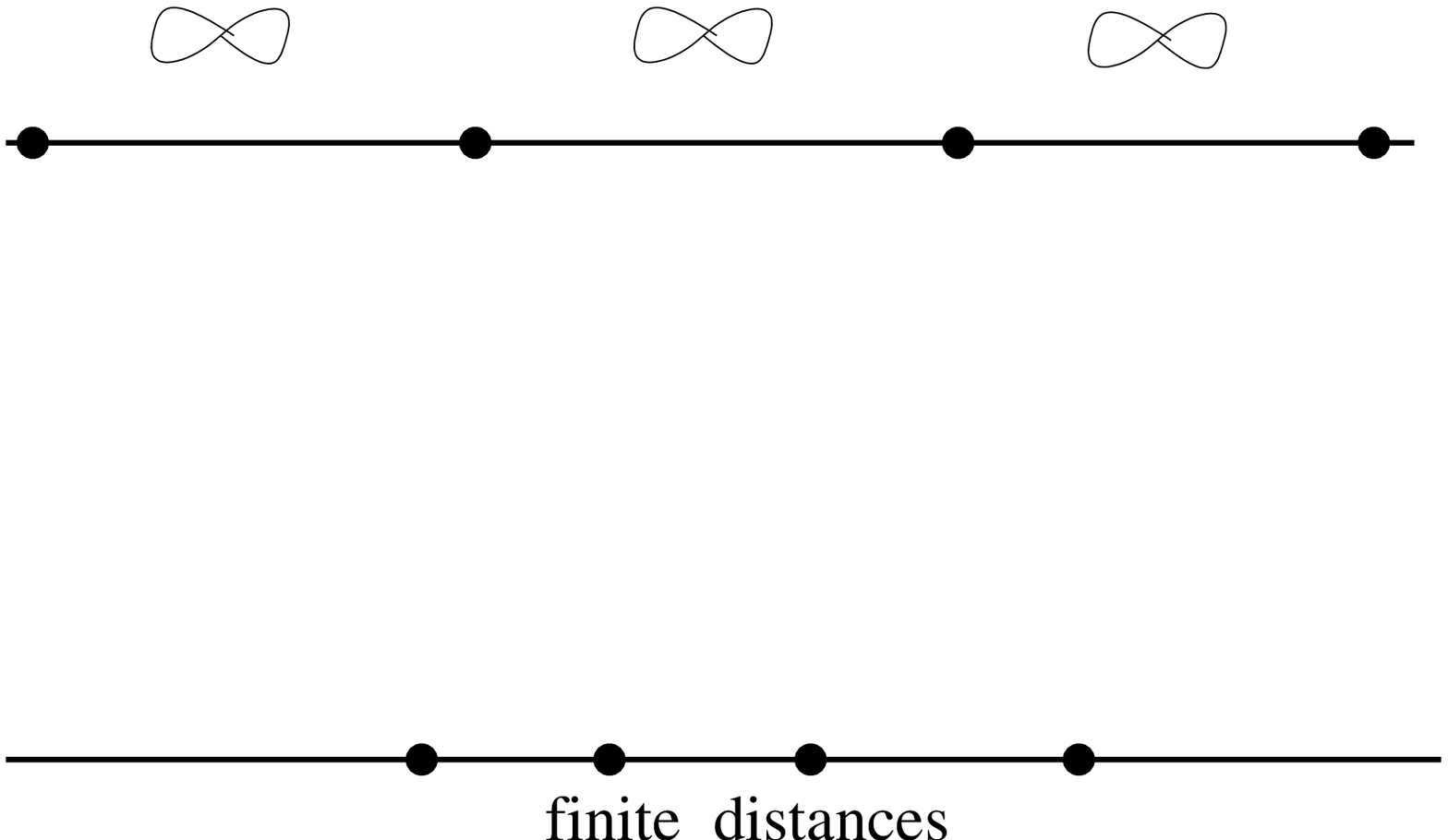}{80mm}{0}{The stratum of codimension 1 in $\overline{K(m,n)}$ of the third type}
Here in the Figure all points on the upper line move infinitely far from
each other {\it will a finite ratio of any two among these infinite
distances}. The all distances on the lower line are finite. Of course, the stratum when
the points on the lower line are infinite and the points on the upper line
are in finite distances from each other, is the same: one stratum can be
obtained from another by the application of the element
$(0,0,\infty)\in G^3$.

The strata in Figures 3 and 4 are isomorphic to $K(m,n)\times K(1,n_1)$ and
$K(m_1,1)\times K(m,n)$, correspondingly. The stratum in Figure 5 is
isomorphic to $K_m^{1,1,\dots,1\ (n\ times)}\times K^n_{1,1,\dots,1\ (m\
times)}$. Let us explain the last formula: the distances between the points
at the upper line are infinite, but their ratios are finite. Therefore, we
should count these ratios. For this, we apply to Figure 5 the transformation
$(0,0,\infty)\in G^3$. Then the infinite distances will become finite, and
then we can count the ratios.

Now we claim that using the 3 operations shown in Figures 3,4,5 we can
obtain {\it any} limit configuration (here by a limit configuration we mean
a configuration where some distances are infinitely large or/and infinitely
small). Moreover, we can apply the configuration in Figure 5 not more than 1
time. Let us explain it:

Apply a transformation from $G^3$ such that there are no infinite distances
at the lower line, and the diameter of the configuration of points at the
lower line is finite (not infinitely small). We distinguish the following
two cases: in the first case in the obtained configuration there are no infinitely large distances in the
upper line, and in the second some distances are infinitely large. It is
clear that we can reach any limit configuration of the first type by the
applying several times the degenerations shown in Figure 3 and Figure 4.
Analogously, in the second case, we first apply degenerations in Figure 3
(several times), then apply the transformation from Figure 5 (with the scale
of the infinity depending on the configuration), and then apply several
times the degenerations from Figure 4. It is clear that in this case we can
get {\it any} limit configuration.

Denote the operations shown in Figures 3,4,5 by $\circ_i$, $_j\circ$, and
$\circledcirc$, correspondingly. It remains to define our operations
$\circledcirc_i$ and $_j\circledcirc$ in the Definition in Section 1.1, that
is, to compactify the spaces $K_{1,1,\dots,1}^n$ and $K_m^{1,1,\dots,1}$.

Notice that the space $K^{1,1,\dots,1\ (n\ times)_m}$ is isomorphic to the
Stasheff polyhedron $\St_m$ for {\it any} $n$, as well the space
$K_{1,1,\dots,1\ (m\ times)}^n$ is isomorphic to the Stasheff polyhedron
$\St_n$ for any $m$. In particular, we define our maps $\curlywedge_{m_1\to
m_2}^n$ and $\curlyvee_m^{n_1\to n_2}$ as the identity maps. Furthermore, we
compactify these spaces as usual in the Stasheff compactification, and the
compositions $\circledcirc_i$ and $_j\circledcirc$ are defined in the
natural way. Moreover, it is clear that the formulas (A) and (B) in the
Definition of $\frac23$PROP hold.

Thus, we constructed a topological $\frac23$PROP $\overline{K(m,n)}$.
\begin{lemma}
\begin{itemize}
\item[(i)] The corresponding chain dg $\frac23$PROP (formed by the chain
complexes of the spaces in our stratification) is a free $\frac23$PROP of
graded vector spaces (when we forget about the differential),
\item[(ii)] the corresponding homology $\frac23$PROP is the $\frac23$PROP of
bialgebras $\Bi$
\end{itemize}
\begin{proof}
(i) follows from the fact that any limit configuration can be written as the
composition of the degenerations in Figures 3,4,5 in a unique way, (ii)
follows from the fact that all spaces $\overline{K(m,n)}$ are contractible
(and, therefore, have only 0-th nontrivial homology).
\end{proof}
\end{lemma}

\begin{example}
In our compactification, the left picture in Figure 1 is $K(2,1)\times
K(1,2)$, and the right picture is $K_{1,1}^2\times K_2^{1,1}$.
\end{example}

\begin{example}
Consider the configuration drawn in Figure 6 below.
\sevafigc{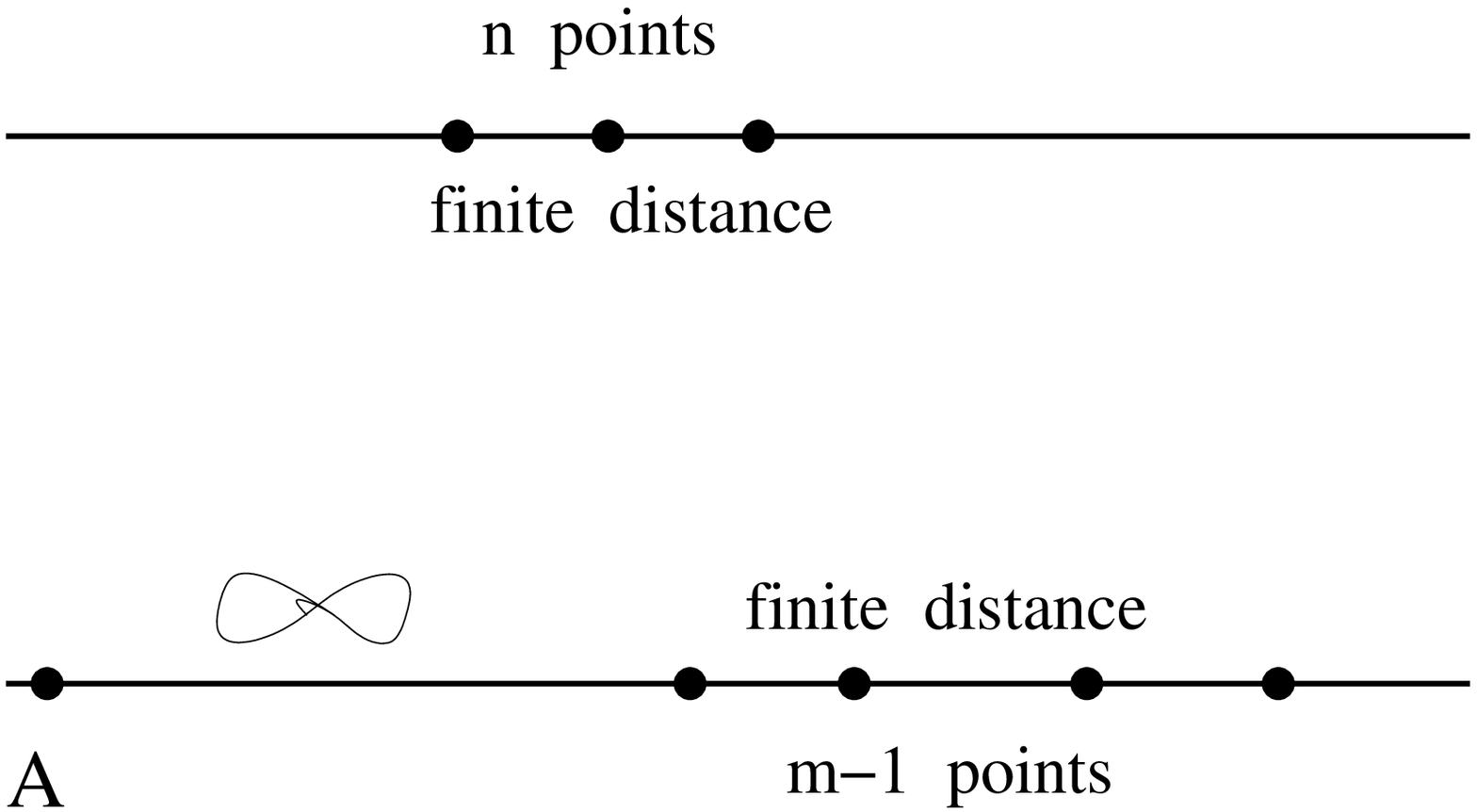}{80mm}{0}{A stratum of codimension 2 for $n>1$}
Here the most left point on the lower line is infinitely far from other
points which are in finite distance from each other, and $n$ points on the
upper line are in finite distance from each other.
This is a stratum of codimension 1 when $n=1$ and of codimension 2 when
$n>1$. Indeed, we can reach the configuration in two steps in our
compactification, see Figure 7.
\sevafigc{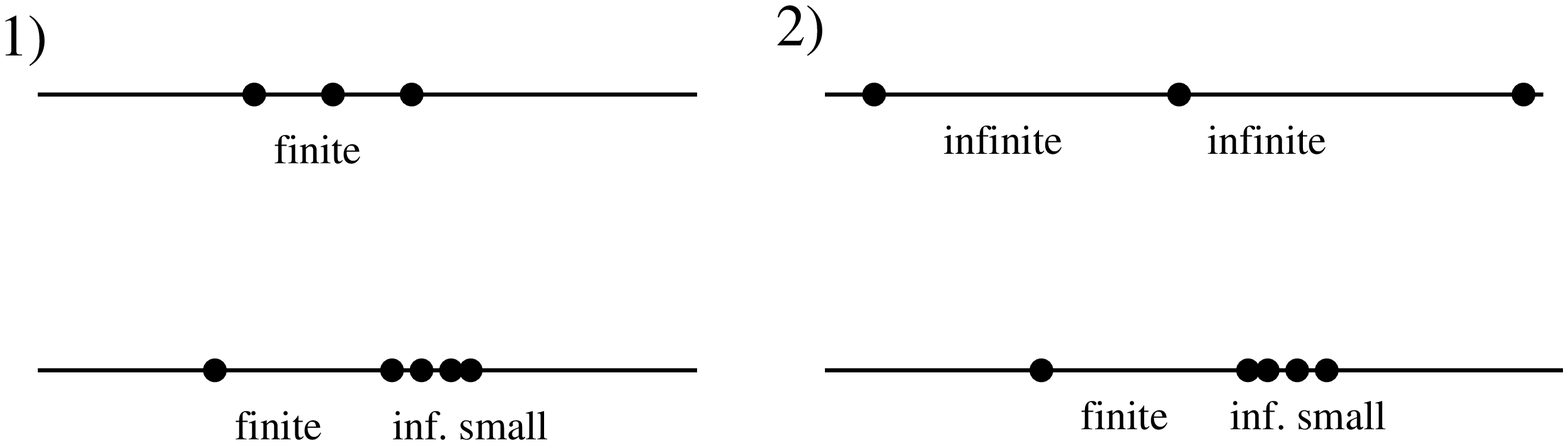}{120mm}{0}{How to find the stratum from Figure 6 in our
compactification}
At first, we use the stratum $K(m-1,1)\times K(2,n)\hookrightarrow
\overline{K(m,n)}$, and then the operation $K_{1,1,\dots,1\ (m\
times)}^n\times K^{1,1,\dots,1\ (n\ times)}_m\to\overline{K(m,n)}$. (See the
fragments 1) and 2) in Figure 7. The second operation is of codimension 0
when $n=1$.
\end{example}
To go further, we introduce the extended Kontsevich space $K(m,n;s)$ (also
due to Maxim Kontsevich) and construct its compactification.
\endcomment
\section{\tt The extended Kontsevich space $K(m,n;s)$ and the Propagator}
\subsection{\tt The definition of the space $K(m,n;s)$}
Consider the direct product $P=\mathbb{R}^2\times\mathbb{R}_{+}$ of the plane $\mathbb{R}^2$ with coordinates $(x,y)$ with the half-line
$\mathbb{R}_{+}$ with the coordinate $\lambda$, $\lambda>0$.
We denote by $(x,y,\lambda)$ the coordinates of a point in $\mathbb{R}^2\times\mathbb{R}_{+}$.
Consider the disjoint union of $P$ with the two lines: $P_1=
P\bigsqcup\{(\mathbb{R},0,0)\}\bigsqcup\{(0,\mathbb{R},\infty)\}$. Define
the following configuration space of points in $P_1$:
\begin{equation}\label{eqw1.1}
\begin{aligned}
\ &\Conf_{m,n;s}=\{p_1,\dots,p_m\in(\mathbb{R},0,0),\  p_i\ne p_j \ for\ i\ne
j,\\
&q_1,\dots,q_n\in(0,\mathbb{R},\infty),\ q_i\ne q_j\ for\ i\ne j,\\
&t_1,\dots,t_s\in\mathbb{R}^2\times\mathbb{R}_+,\ t_i\ne t_j\ for\ i\ne j\}
\end{aligned}
\end{equation}
The action of the 3-dimensional group $G^3$ on $\Conf_{m,n}$ (see Section 2)
can be continued to an action on $\Conf_{m,n;s}$. Indeed, define a product
on $\mathbb{R}^2\times\mathbb{R}_+$ as follows:
\begin{equation}\label{eqw1.2}
(x_1,y_1,\lambda_1)\cdot(x,y,\lambda)=(\lambda_1 x+x_1,\lambda_1^{-1}y+y_1,\lambda\lambda_1)
\end{equation}
This group $\mathbb{R}^2\ltimes\mathbb{R}$ is exactly $G^3$. Then it is
clear that it acts on $\Conf_{m,n;s}$ (see formula (\ref{eq1.3})). This action is free when $m+n+s\ge
3$. Denote in this case the quotient space $\Conf_{m,n;s}/G^3$ by
$K(m,n;s)$. It is a smooth manifold of dimension $m+n+3s-3$. We are going to
compactify this space in a way compatible with the compactification of
$K(m,n)$ introduced in Section 2.
\subsection{\tt The 3-dimensional Eye and the Propagator}
\subsubsection{\tt The 3-dimensional Eye}
Here we compactify the space $K(0,0;2)$ which is 3-dimensional. We call the
space $\overline{K(0,0;2)}$ the {\it 3-dimensional Eye} by the analogy with
the Kontsevich Eye ([K1], Sect. 5.2). The Propagator constructed below will
be a form (with singularities) on the 3-dimensional Eye. On the other hand,
this simplest example illustrates (the all) new ideas which appear in the
compactification of the extended Kontsevich spaces. We describe all strata
of codimension 1 in $\overline{K(m,n;s)}$ in the next Subsection.

We consider an oriented pair of points $(t_1,t_2)$ in
$\mathbb{R}^2\times\mathbb{R}_+$. Then using the action of the group $G^3$
we can fix one of them, say $t_1$. Then we suppose that $t_1=(0,0,1)$ and
there is no group action. Then we should compactify the space
\begin{equation}\label{eq11.1}
P=\{t=(x,y,\lambda)\in\mathbb{R}^2\times\mathbb{R}_+,\ (x,y,\lambda)\ne
(0,0,1)\}
\end{equation}
It is easy to compactify near the point $(0,0,1)$: we just cut-off a sphere
around this point. The other boundary components arise when the point $t$ tends
to $\infty$. Let us describe them:

First consider the case when the point $t$ tends to infinity when the
coordinate $\lambda$ is finite. Consider the projection $p\colon
(x,y,\lambda)\mapsto (x,y)$. The the two lines (for the values $\lambda=0$
and $\lambda=\infty$ are the coordinate axis (see Figure 8 below). The
degenerations of codimension 1 are then the configurations when {\it only
one} among the two coordinates after the projection tends to infinity. They
are "separated" by faces of codimension 2 when the both coordinates $x,y$
tend to $\infty$. Let us explain how we compute the dimension of the strata:
when one coordinate tends to $\infty$, it is defined only up to a finite
summand. Therefore the only coordinates on the moduli is another coordinate
on the plane (after the projection), and the coordinate $\lambda$.
Therefore, the strata of codimension 1 are "separated" by the strata of
codimension 2 corresponding to the cases when both $x,y$ tend to $\infty$.
More precisely, they are separated by the lines $x=y$ and $x=-y$ in Figure
8. Thus, we have 4 faces of codimension 1 which are corresponded to the
following 4 cases:
\begin{itemize}
\item[(i)] $x\gg 0$, $y$ is finite (positive or negative),
\item[(ii)] $y\gg 0$, $x$ is finite (positive or negative),
\item[(iii)] $x\ll 0$, $y$ is finite (positive or negative),
\item[(iv)] $y\ll 0$, $x$ is finite (positive or negative).
\end{itemize}
\sevafigc{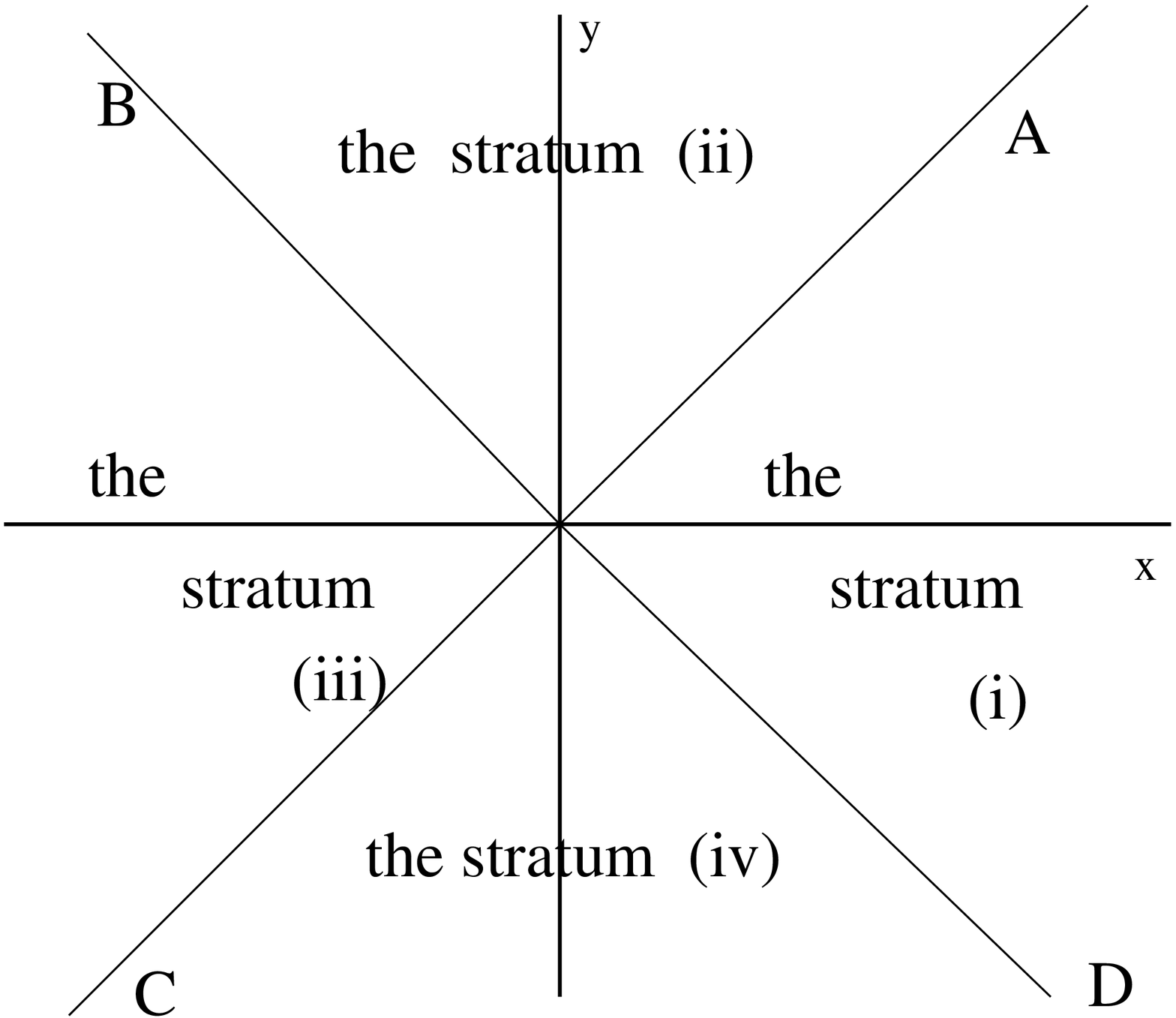}{100mm}{0}{The strata of codimension 1 in 2-dimensional
projection}
These are the all strata of codimension 1. There are 6 strata of codimension
2. First 4 among are:
\begin{itemize}
\item[(A)] $x\gg 0$, $y\gg 0$,
\item[(B)] $x\ll 0$, $y\gg 0$,
\item[(C)] $x\ll 0$, $y\ll 0$,
\item[(D)] $x\gg 0$, $y\ll 0$.
\end{itemize}
The last 2 faces of codimension 2 are obtained when $\lambda$ tends to 0 and
to $\infty$. We show these strata in Figure 8. Finally, The 3-dimensional
Eye is the tetrahedron without a small ball inside (corresponded to the case
when $t$ is close to $(0,0,1)$ (see Figure 9). We designated the strata of
codimension 1 in Figure 9.
\sevafigc{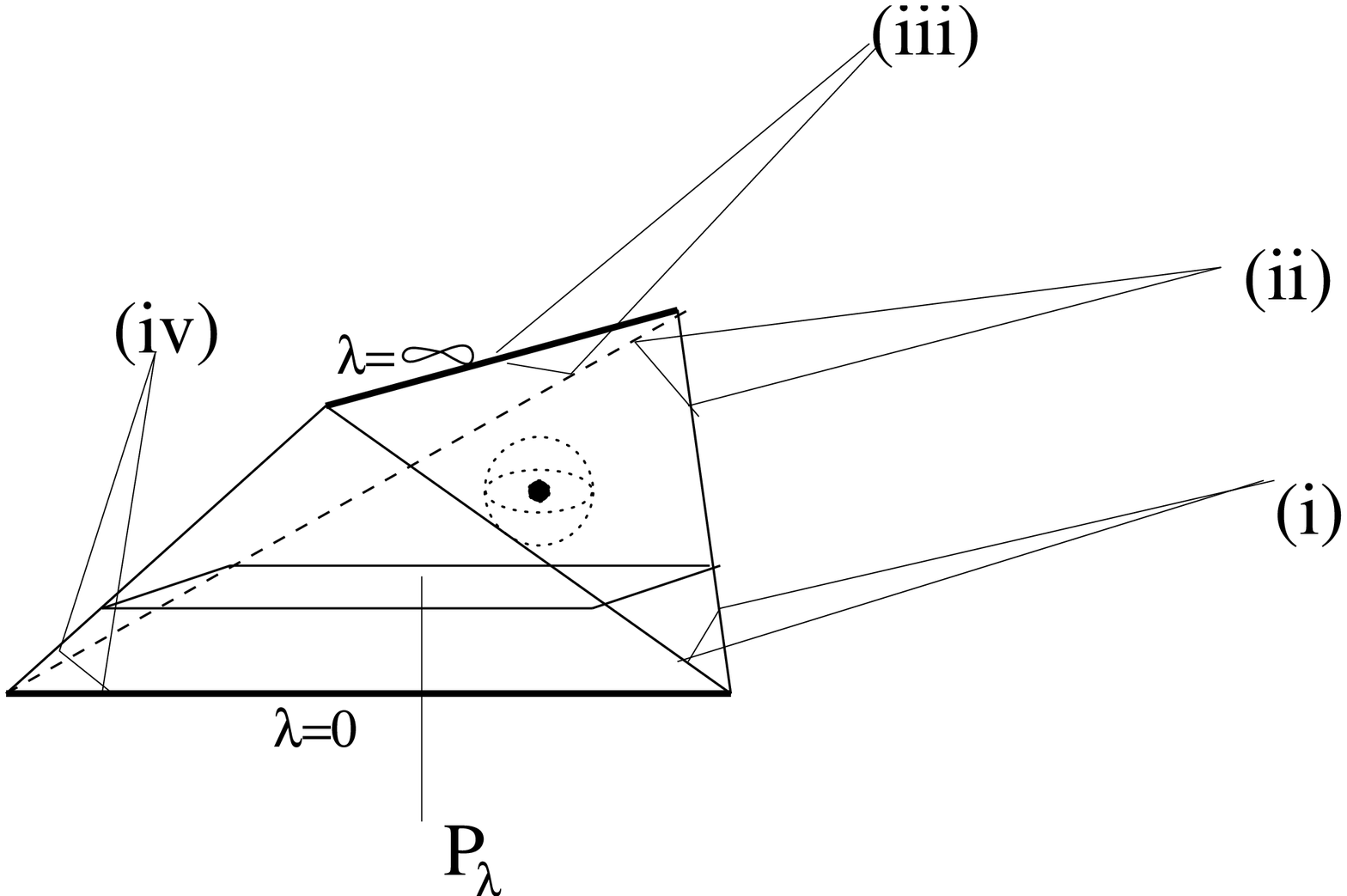}{120mm}{0}{The 3-dimensional Eye}
\subsubsection{\tt The Propagator}
We define a Propagator 2-form as a {\it closed} 2-form $\phi$ {\it with singularities} on the 3-dimensional Eye
such that:
\begin{itemize}
\item[1)] the form $\phi$ has singularities only at the edge
$\{\lambda=0\}$ of the tetrahedron,
\item[2)] the restriction of the form $\phi$ to the 2-dimensional sphere
(around the point $(0,0,1)$) is the volume form on the sphere normalized
such that the integral over the sphere is equal to 1,
\item[3)] the restriction of the form $\phi$ to the other boundary
components (besides the sphere and the interval $\{\lambda=0\}$) is 0,
\item[4)] the singularity at the edge $\{\lambda=0\}$ has the form described
below.
\end{itemize}
Consider the rectangle $P_\lambda$ which is the horizontal section of the
tetrahedron by the plane $z=\lambda$ where $\lambda$ is close to 0 (see
Figure 9). Let $x$ and $y$ be the coordinates on the rectangle. We say that
$x\in [-1,1]$ and $y\in [-\epsilon,\epsilon]$ for a "very small" $\epsilon$.
Roughly speaking, we want that in the limit $\lambda\to 0$ the restriction of
the Propagator 2-form to our rectangle would be $f(x)\delta(y)dx\wedge dy$
where $\delta(y)$ is the Dirac delta-function with the support at $y=0$, and
$f(x)$ is any {\it positive} function with the support in the open interval
$(-1,1)$ and such that $\int_{-1}^1f(x)dx=1$. Notice that the Propagator
2-form should be not defined when $t\in\{\lambda=0\}$, there we have the
Propagator 1-form which is, by definition, the 1-form $f(x)dx$ where $f(x)$
is as above. So, the above description is {\it not} a description at
$\lambda=0$, but when $\lambda\to 0$.

Let us prove that such a Propagator 2-form exists (it is more or less clear
that if it exists it is defined up to a homotopy in an appropriate sense):

Consider the part of the 3-dimensional Eye which is above the rectangle
$P_\lambda$ (see Figure 9), where $\lambda$ is a sufficiently small number
(it is important only that the cute-off sphere is above it). We call this
part the {\it truncated} 3-dimensional Eye.

\begin{definition}
We say that a map of the truncated 3-dimensional Eye to the 2-dimensional
sphere $S^2$ is {\it spherical} it maps the sphere inside the tetrahedron
homotetically to the sphere $S^2$, maps all the boundary of the truncated tetrahedron
{\it except} the rectangle $P_\lambda$ to a point $p\in S^2$, and it maps
the domain $D\subset P_\lambda$ to the sphere $S^2$ and $P_\lambda\setminus D$
to the point $p\in S^2$ such that the factor-space $D/\partial D$ maps
isomorphically to $S^2$ in the way compatible with the orientation. Here the
domain $D$ is any domain like it is shown in Figure 10. It should be any
simply-connected domain having a non-empty intersection with $\{y=0\}\subset
P_\lambda$.
\end{definition}
\sevafigc{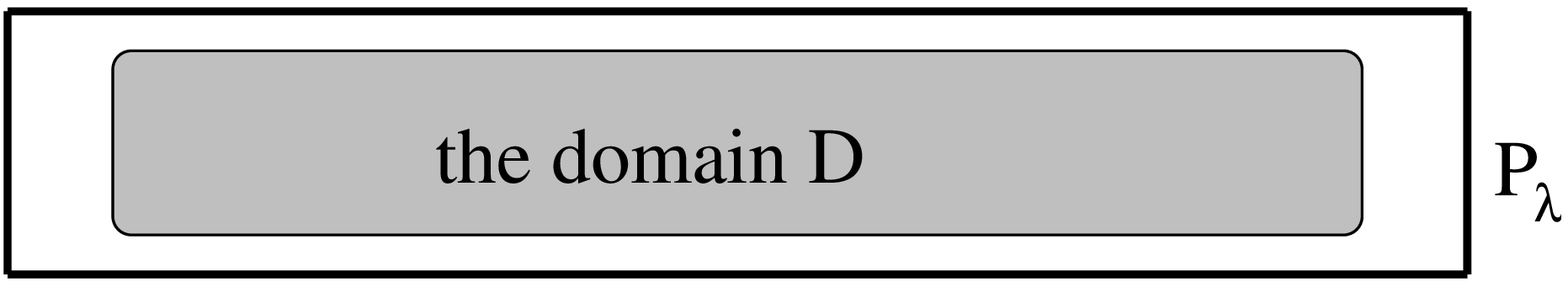}{80mm}{0}{The domain $D$ in the rectangle $P_\lambda$}

It is clear that such a map exists and is homotopically unique.
Then, starting with a spherical map, we can easily construct a Propagator
2-form. First, define it on the truncated rectangle as $\pi^*(\omega)$ where
$\pi$ is a spherical map, and $\omega$ is the volume form on $S^2$
normalized such that $\int_{S^2}\omega=1$. Then consider the cone which is
the complement to the truncated tetrahedron in the tetrahedron, and continue
the Propagator form to a form with singularities in a natural way. We should
be careful about the smoothness of the form, but it is possible to make it
smooth. Then we define the Propagator 1-form at the edge $\{\lambda=0\}$ as
follows: we consider the direct image of the restriction of the Propagator
2-form on $P_\lambda$ with respect to the projection of $P_\lambda$ to the
$P_\lambda \cap \{y=0\}$.
\subsection{\tt The compactification $\overline{K(m,n;s)}$}
Here we define a compactification $\overline{K(m,n;s)}$ of the space
$K(m,n;s)$.

Recall our notations $\Sigma_{1,3}$, $\Sigma_{2,2}$, and $\Sigma_{3,1}$ for
the possible samples of 4 points on the 2 boundary lines (see Section 1.1).
Now denote by $\Sigma^2$ the all possible pairs of points in $K(m,n;s)$
among the $s$ inner points, and denote by $\Sigma^1_1$ the all possible
pairs of points in $K(m,n;s)$ one of which is an inner point, and another is
a points on the upper or lower boundary line.
For $\sigma\in\Sigma^2$ we have a map $r_\sigma\colon K(m,n;s)\to K(0,0;2)$,
and as well for $\sigma\in\Sigma_1^1$ we have a map $r_\sigma\colon
K(m,n;s)\to K(1,0;1)$ or to $K(0,1;1)$. In the sequel we will denote the
last two spaces by one symbol $K(\P;1)$. We have the direct product of the
maps $r_\sigma$:
\begin{multline}\label{eq777.1}
r=\prod_{\sigma\in\Sigma_{1,3}\bigsqcup\Sigma_{3,1}\bigsqcup\Sigma_{2,2}\bigsqcup\Sigma^2\bigsqcup\Sigma^1_1}
r_\sigma\colon\\
K(m,n;s)\to\prod_{\sigma\in\Sigma_{1,3}}K(1,3)\times\prod_{\sigma\in\Sigma_{2,2}}K(2,2)\times
\prod_{\sigma\in\Sigma_{3,1}}\K(3,1)\times\\
\prod_{\sigma\in\Sigma^2}K(0,0;2)\times\prod_{\sigma\in \Sigma_1^1}K(\P;1)
\end{multline}
It is clear that the map $r$ is an imbedding. We know how to compactify each
space in the right-hand side. It is the 3-dimensional eye for $K(0,0;2)$,
and a part of the 3-dimensional eye (which we imbed into the 3-dimensional
eye) for $K(\P;1)$. Then we can imbed each space in the r.h.s to its
compactifacation, and then take the closure of the image in this imbedding.
This is, by definition, our compactification $\overline{K(m,n;s)}$.
\begin{remark}
In the Kontsevich paper on deformation quantization, one also can compactify
the space $C_{m,n}$ in this way. But for a pair of inner points we should
consider the corresponding point of the Kontsevich eye. Kontsevich
attached just the corresponding angle for a pair of inner points. In this
way, he does not get an imbedding, and therefore automatically he gets a
wrong compactifacation, which does not coincide with his right
compactification described in the terms of trees. When one considers a point
of the Kontsevich eye instead of the corresponding angle, we get the right
described in trees compactification.
\end{remark}

We claim the following lemma:
\begin{lemma}
\begin{itemize}
\item[(i)] The space $\overline{K(m,n;s)}$ is naturally a manifold with
corners,
\item[(ii)] The projection $p\colon K(m,n;s)\to K(m,n)$ can be uniquely
continued to a map
$\overline{p}\colon\overline{K(m,n;s)}\to\overline{K(m,n)}$ which is a map
of manifolds with corners.
\end{itemize}
\qed
\end{lemma}

Now we describe all strata of codimension 1 in $\overline{K(m,n;s)}$. We
define the dimension "of a point" as follows: it is equal to $k$ if we need
to know $k_1$ 1-dimensional 4-point ratios, $k_2$ 3-dimensional 4-point
ratios, $k_1+3k_2=k$, and this number can be not made less.

The images with respect to $\overline{p}$ of the strata of codimension 1 in
$\overline{K(m,n;s)}$ are either the stratum of codimension 0 in
$\overline{K(m,n)}$ (we call these strata of codimension 1 in
$\overline{K(m,n;s)}$ {\it the strata of codimension 1 of the first type}),
or are strata of codimension 1 in $\overline{K(m,n)}$ (we call them {\it
the strata of codimension 1 of the second type}).
First list all the strata of codimension 1 in $\overline{K(m,n;s)}$.

{\it Strata of codimension 1 of the first type}
\begin{itemize}
\item[S1.1] $s_1$ points from the $s$ "inner" points move close to each
other, to a finite point, and for a finite value of $\lambda$, $s_1\ge 2$,
\item[S1.2] $s_1$ "inner" points at a finite distance from each other move
to infinity for a finite $\lambda$, and such that only one coordinate among
$x,y$ tends to $\infty$.
\end{itemize}
{\it Strata of codimension 1 of the second type}
\begin{itemize}
\item[S2] A typical stratum is drawn in Figure~7. The  reader should
remember that the two lines are {\it crossing} (that means two lines in a
3-dimensional space which do not intersect). On these two boundary lines we
have exactly the picture from Figure 2. The points inside the circle {\bf A}
are in finite distances everywhere. The points inside the circles {\bf B}
and {\bf C} are infinitely far to the right and to the left from the circle
{\bf A}. When we apply an infinite shift to the circle {\bf B} the points
inside it will be everywhere in the finite distances, the same is true for
the circle {\bf C}. (Remember that the boundary lines are {\it crossing}!).
On the lower line we have several groups of points infinitely close to each
other, the distances between the groups are finite. The lower and the upper
infinities have the same order. When we apply the transform
$(0,0,\lambda)\in G^3$ for infinite $\lambda$ to Figure~7, the lower line
will look like the upper, and vise versa.
\end{itemize}
\begin{remark}
It was a remarkable insight of Maxim Kontsevich that the lines should be
crossing.
\end{remark}
\sevafigc{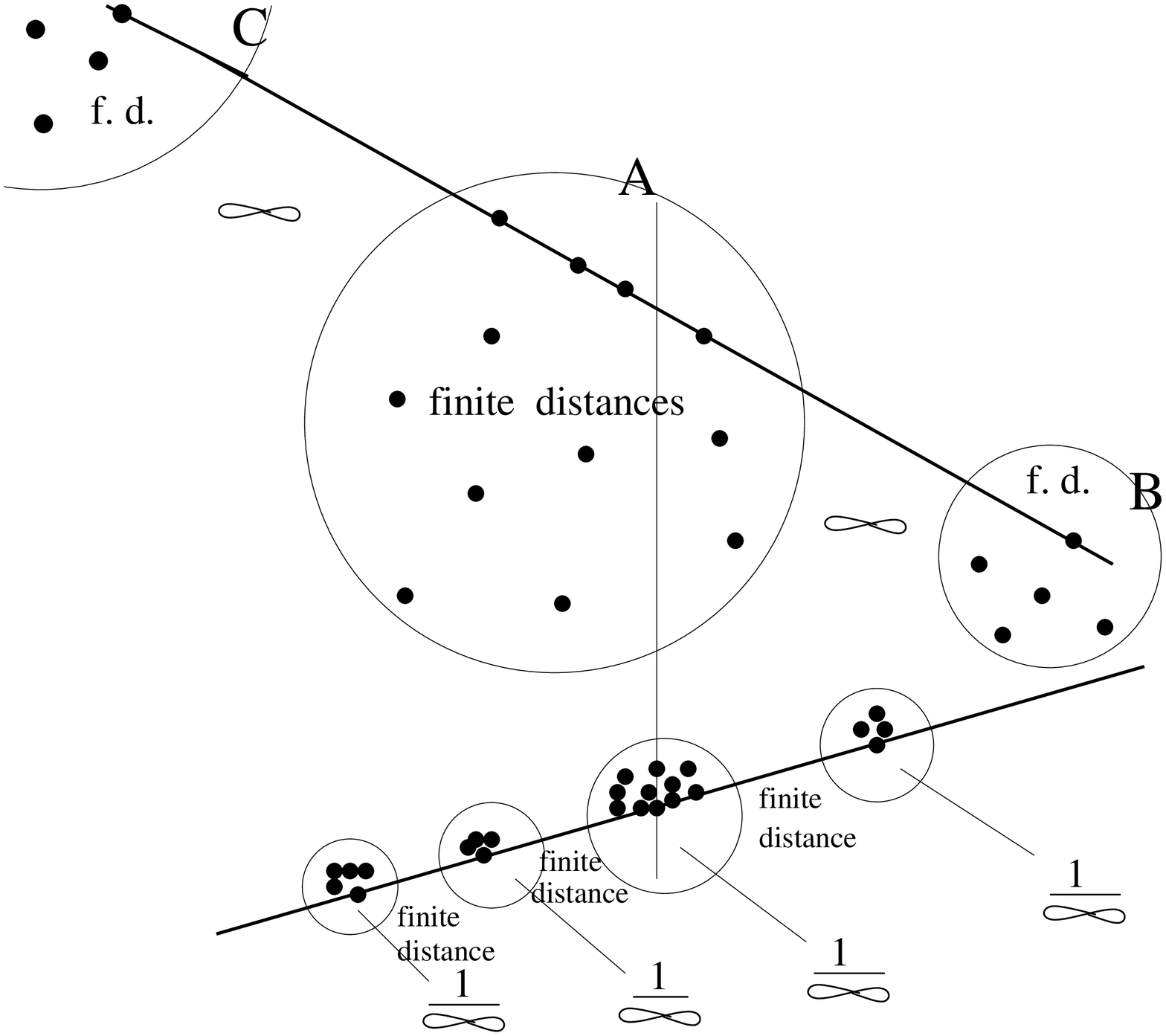}{120mm}{0}{A typical stratum S2 of codimension 1}
\subsection{\tt The direct image $\overline{p}_*(\omega)$ of a closed form
$\omega$ on $\overline{K(m,n;s)}$}
Consider the projection
$\overline{p}\colon\overline{K(m,n;s)}\to\overline{K(m,n)}$. Let $\omega$ be
a smooth differential form (not necessarily closed) on
$\overline{K(m,n;s)}$. The smoothness here means that it is a continuous form
such that the restriction of it to each open stratum is smooth. What is the
"right definition" of the pushforward $\overline{p}_*(\omega)$ which {\it
should be an element of the bicomplex} $B^{\mb\mb}(m,n)$? The usual
definition as the integration over the fiber is not good because it is a
trivial bundle over a single stratum $\sigma$, but the fiber changes from
stratum to stratum.
\begin{definition}
Let $\overline{p}\colon M_1\to M_2$ be a map of manifolds with corners. Let
$\omega$ be a homogeneous element of some degree $k-\dim M_1$ in $B^{\mb\mb}(M_1)$.
It means that $\omega=\omega_0^i+\omega^j_1+\omega_2^t+\dots$ where $\omega_0^i$ is
a form of degree $k$ on the $i$th stratum of codimension 0, $\omega_1^j$ is
a form of degree $k-1$ on the $j$th stratum of codimension 1, and so on.
Define the direct image $\overline{p}_*\colon B^{\mb\mb}(M_1)\to
B^{\mb\mb}(M_2)$ as the following map of degree 0:
The image $\overline{p}_*(\omega)$ is defined as the sum $\sigma_{m,n}
\overline{p}_*\omega_m^n$ over all $m,n$. Here by $\overline{p}_*$ we mean
the integration over the fiber. By the definition of a map of manifold with
corners, the image of each  stratum $\sigma$ in $M_1$ is a stratum
$\overline{p}_*(\sigma)$ in $M_2$, and $\dim F+\dim
\overline{p}_*(\sigma)=\dim\sigma$ where $F$ is the fiber over
$\overline{p}_*(\sigma)$. The degree of the form
$\overline{p}_*(\omega_m^n)$ is $\deg\omega_m^n-\dim F$, therefore $\deg
\omega_m^n-\dim\sigma=(\deg\omega_m^n-\dim
F)-\dim\overline{p}_*(\sigma)=\deg\overline{p}_*(\omega_m^n)-\dim\overline{p}_*(\sigma)$.
Therefore, we obtain a map $\overline{p}_*\colon B^{\mb\mb}(M_1)\to
B^{\mb\mb}(M_2)$ of degree 0.
\end{definition}
As we will see now, the following result is just an application of the
Stokes formula.
\begin{proposition}
The map $\overline{p}_*\colon B^{\mb\mb}(M_1)\to B^{\mb\mb}(M_2)$ is a map
of (total) complexes.
\begin{proof}
First of all, we need the following "relative" version of the Stokes
formula:
\begin{lemma}
Let $E$ and $B$ be smooth manifolds with boundary, and let $p\colon E\to B$
be the trivial bundle with fiber $F$ which is a compact manifold with boundary (that
is, $E=B\times F$). Let $\omega$ be a differential form on $E$, we denote by
$p_*\omega$ the direct image which is the integral of $\omega$ along the
fiber $F$. Then we have:
\begin{equation}\label{eq160_1}
p_*d\omega=(p|_{\partial F})_*\omega+dp_*\omega
\end{equation}
where $d$ is the de Rham differential.
\begin{proof}
As $p$ is a trivial bundle, we can decompose the de Rham differential on $E$
into the sum of its horizontal component $d_{hor}$, and its vertical
component $d_{vert}$, $d=d_{hor}+d_{vert}$. Therefore, $p_*d\omega=
\int_Fd_{hor}\omega+\int_Fd_{vert}\omega$. By the Stokes formula, the second
summand is $\int_{\partial F}\omega$. The first summand is
$d_{hor}\int_F\omega=dp_*\omega$.
\end{proof}
\end{lemma}

Now we prove the Proposition.
It is enough to consider the case when the element $\omega$ of the bicomplex
has only one component, which is a form $\Omega$ on a simplex $\sigma$.
We need to prove:
\begin{equation}\label{eq160_3}
\overline{p}_*(d\Omega)-\overline{p}_*\Omega|_{\partial\sigma}=
d\overline{p}_*\Omega-(\overline{p}_*\Omega)|_{\partial\overline{p}\sigma}
\end{equation}
which follows immediately from Lemma above and from the straightforward
relation
\begin{equation}\label{eq160_4}
\overline{p}_*\Omega|_{\partial\sigma}-(\overline{p}_*\Omega)|_{\overline{p}\partial\sigma}=
(\overline{p}|_{\partial F})_*\Omega
\end{equation}
where $F$ is the fiber over $\sigma$.
\end{proof}
\end{proposition}
\comment
\subsection{\tt The homotopical Gerstenhaber-Schack complex}
\subsubsection{\tt The Gerstenhaber-Schack space and a space homotopically
equivalent to it: an introduction.}
Let $A$ be a vector space (the reader can think about it as about the
(co)associative bialgebra with 0 product and 0 coproduct). The corresponding
Gerstenhaber-Schack space is, by definition, the following graded space:
\begin{equation}\label{eq333.1}
K^\mb=\oplus_{m,n\ge 1}\Hom(A^{\otimes m}, A^{\otimes n})[-m-n+2]
\end{equation}
Our goal is to equip $K^\mb$ with an $L_\infty$ structure associated with
the strata of codimension 1 in $\overline{K(m,n)}$. The main problem we meet
here is the following:

Consider $\Psi_1,\dots,\Psi_n\in\Hom(A^{\otimes m},A)$ and
$\Theta_1,\dots,\Theta_m\in\Hom(A,A^{\otimes n})$. We could attach to the
stratum of codimension 1 in the $\frac23$PROP $\overline{K(m,n)}$ $\sigma=
K_{1,1,\dots,1\ (m\ times)}^n\times K^{1,1,\dots,1\ (n\ times)}_m\hookrightarrow
\overline{K(m,n)}$ the "$L_\infty$-component"
$(\Psi_1\otimes\dots\otimes \Psi_n)\circledcirc
(\Theta_1\otimes\dots\otimes\Theta_m)\in\Hom(A^{\otimes m},A^{\otimes n})$
because according to the general principles the $L_\infty$ components should
correspond to the strata of codimension 1. It means, that we want to define
a map
\begin{equation}\label{eq333.2}
L_{m+n}\colon \wedge^{m+n}K^\mb\to K^\mb[2-m-n]
\end{equation}
as
\begin{equation}\label{eq333.3}
L_{m+n}(\Psi_1\wedge\dots\wedge\Psi_n\wedge\Theta_1\wedge\dots\wedge\Theta_m)=
\Alt ((\Psi_1\otimes\dots\otimes \Psi_n)\circledcirc
(\Theta_1\otimes\dots\otimes\Theta_m))
\end{equation}
Here we meet our first problem: the degree of
$L_{m+n}((\Psi_1\wedge\dots\wedge\Psi_n\wedge\Theta_1\wedge\dots\wedge\Theta_m))$
is by (\ref{eq333.2}) $n(m-1)+m(n-1)+1-m-n=2mn-2(m+n)+2$, but on the other
hand, it belongs to $\Hom(A^{\otimes m}, A^{\otimes n})$ and therefore has
the degree $m+n-2$. The construction below was invented to solve this
problem.

Replace the greaded vector space $K^\mb$ by the following complex:
\begin{equation}\label{eq333.4}
\begin{aligned}
\ &\overline{K}^\mb=\\
&\Hom(A,A)[0]\oplus \bigoplus_{m,n\ge 1, m+n\ge 3}(
\Hom(A^{\otimes m},A^{\otimes n})[-m-n+2]\otimes\\
&\otimes(\oplus_{all\  strata\  \sigma\  in\  K(m,n)}
(\Omega^\mb_{DR}(\sigma),d_{DR}),\partial))
\end{aligned}
\end{equation}
Let us explain what is written in the last equation: in the r.h.s. we take
the tensor product of $\Hom(A^{\otimes m},A^{\otimes n})$ with the sum of
the de Rham complexes on all (closed) strata in $\overline{K(m,n)}$ equipped
with the chain differential. For any imbedding
$\sigma_0\hookrightarrow\sigma$ of strata we have the corresponding map of
the de Rham complexes $\Omega^\mb_{DR}(\sigma)\to\Omega^\mb_{DR}(\sigma_0)$
which allows us to consider the sum of the de Rham complexes over all open
strata in $\overline{K(m,n)}$ as a bicomplex with the de Rham and the chain
differential. We consider the chain complex as negatively-graded, and the de
Rham complex as positively graded. It is clear that the bicomplex has only
1-dimensional cohomology in degree 0 because the space $\overline{K(m,n)}$
is contractible. Therefore, $\overline{K}^\mb$ is homotopically equivalent
to $K^\mb$.

Introduce the notation: $[\omega]_{\ell}$ means a form of some degree
$\deg\omega$ on a stratum of dimension $\ell$ of some space
$\overline{K(m,n)}$ considered as an element of the bicomplex. In
particular, the degree $\deg [\omega]_{\ell}$ in the bicomplex is
$\deg\omega -\ell$.

We consider $\overline{K}^\mb$ as a more right object than $K^\mb$. We
explain below that $\overline{K}^\mb$ has a natural $L_\infty$-structure.

Consider the previous construction with $L_{m+n}$ for $\overline{K}^\mb$.
For simplicity, consider $\wtilde{\Psi}_1=\Psi_1\otimes[\omega_1]_{m-2},\dots,
\wtilde{\Psi}_n=\Psi_n\otimes[\omega_n]_{m-2}$ and
$\wtilde{\Theta}_1=\Theta_1\otimes[\omega_1^\prime]_{n-2},\dots,\wtilde{\Theta}_m=
\Theta_m\otimes[\omega_m^\prime]_{n-2}$ where $Psi_i$s and $\Theta_j$s are as above,
and all $\omega$s and $\omega^\prime$s are some homogeneous differential
forms (of some degrees) on the {\it top dimension stratum} of
$\overline{K(m,1)}$ and $\overline{K(1,n)}$. The case of strata of not top
dimension is considered below.
Then the degree of
$L_{m+n}(\wtilde{\Psi}_1\wedge\dots\wedge\wtilde{\Psi}_n\wedge
\wtilde{\Theta}_1\wedge\dots\wedge\wtilde{\Theta}_m)$ is equal to
\begin{equation}\label{eq333.6}
n(m-1)+m(n-1)+\sum_{i=1}^n(\deg
\omega_i-m+2)+\sum_{j=1}^m(\deg\omega_j^\prime-n+2)+2-m-n
\end{equation}
(Here $\deg\omega_i-m+2$ is the degree of the element
$\omega\otimes\sigma_{max}$ in the bicomplex where $\sigma_{max}$ is the top
dimension stratum in $\overline{K(m,1)}$).
On the other hand, we can attach to this data, at first, the $\frac23$PROP
product $\Alt ((\Psi_1\otimes\dots\otimes \Psi_n)\circledcirc
(\Theta_1\otimes\dots\otimes\Theta_m))$, and, at second, a differential form
$\omega$ at the stratum $\sigma=
K_{1,1,\dots,1\ (m\ times)}^n\times K^{1,1,\dots,1\ (n\ times)}_m\hookrightarrow
\overline{K(m,n)}$ of {\it codimension 1} in $\overline{K(m,n)}$ of degree
$\deg\omega=\sum\deg\omega_i+\sum\deg\omega^\prime_j$. We explain the
construction of this $\omega$ in the Subsections below. Now let us compute the degree:
\begin{equation}\label{eq333.7}
\begin{aligned}
\ &\deg(\Alt ((\Psi_1\otimes\dots\otimes \Psi_n)\circledcirc
(\Theta_1\otimes\dots\otimes\Theta_m))\otimes[\omega]_{m+n-4}=\\
&m+n-2+(\deg\omega-(m+n-4))
\end{aligned}
\end{equation}
where $m+n-4$ in the formula above is the dimension of the stratum of
codimension 1 in $\overline{K(m,n)}$. We see that in the assumption
\begin{equation}\label{eq333.8}
\deg\omega=\sum\deg\omega_i+\sum\deg\omega_j^\prime
\end{equation}
the right hand sides of (\ref{eq333.6}) and of (\ref{eq333.7}) are equal,
and then the equality of the degrees for $L_{m+n}$ holds. If
$\sum\deg\omega_i+\sum\deg\omega_j^\prime>n+m-3=\dim K(m,n)$, we define the
corresponding operation as 0. This is the cause of our original problem with
definition of $L_{m+n}$ for $K^\mb$: we imbed $K^\mb\to\overline{K}^\mb$
(not canonically) mapping an element $\Psi\in\Hom(A^{\otimes m}, A^{\otimes
n})$ to some element $\Psi\otimes[\omega]_{m+n-3}$ where $\omega$ is a top
degree form on the open higher dimension stratum. Then the degrees of
$\omega_i$'s and of $\omega^\prime_j$'s are very big, and we obtain zero.

\subsubsection{\tt The construction of the form $\omega$ and the strata of higher
codimension}
First we should construct the form $\omega$ on the stratum
$K_{1,1,\dots,1\ (m\ times)}^n\times K^{1,1,\dots,1\ (n\ times)}_m$
of codimension 1
in $\overline{K(m,n)}$ starting from the forms $\omega_i$'s and
$\omega_j^\prime$'s promised in the previous Subsection. It is very easy:

Consider the projection $p_i\colon K^{1,1,\dots,1\ (n\ times)}_m\to K(m,1)$
which just forget about all upper lines (points) in $K^{1,1,\dots,1\ (n\
times)}$ except the $i$th ($1\le i\le n$). (Actually this projection is an
isomorphism). Take then the wedge product
$\wedge_{i=1}^np_j^*\omega_i=\omega^{(1)}$ which is a differential form on $K^{1,1,\dots,1\ (n\
times)}$. Analogously, consider the projection $p^\prime_j\colon K_{1,1,\dots,1\ (m\
times)}^n\to K(1,n)$ and define $\omega^{(2)}=\wedge_{j=1}^mp^{\prime
*}_j\omega_j^\prime$. Then define $\omega=\omega^{(1)}\boxtimes\omega^{(2)}$
as a differential form on $K_{1,1,\dots,1\ (m\ times)}^n\times K^{1,1,\dots,1\ (n\
times)}_m$ which is a stratum of codimension 1 in $\overline{K(m,n)}$.

Concerning the strata of higher codimension, the construction of the
previous Subsection can be generalized as follows: we can consider $\wtilde{\Psi}_1=
\Psi_1\otimes[\omega_1]_{a},\dots,
\wtilde{\Psi}_n=\Psi_n\otimes[\omega_n]_{a}$ where
$\omega_1,\dots,\omega_n$ are some differential forms on {\it the same for all}
$\omega$'s stratum of dimension $a$
in $\overline{K(m,1)}$. As well,
consider $\wtilde{\Theta}_1=\Theta_1\otimes[\omega_1^\prime]_{b},\dots,\wtilde{\Theta}_m=
\Theta_m\otimes[\omega_m^\prime]_{b}$ where
$\omega^\prime_1,\dots,\omega^\prime_m$ are some differential forms on {\it
the same for all} $\omega^\prime$'s
stratum of dimension $b$ in $\overline{K(1,n)}$.
Then the stratum $\sigma=K_{1,1,\dots,1\ (m\ times)}^n\times K^{1,1,\dots,1\ (n\
times)}_m\hookrightarrow\overline{K(m,n)}$ of codimension 1 defines the
correspondence which associates to a stratum of dimension $a$ in $K^{1,1,\dots,1\ (n\
times)}_m$ and a stratum of dimension $b$ in $K_{1,1,\dots,1\ (m\ times)}^n$
a stratum $\sigma_{tot}$ of dimension $a+b$ (of the codimension equal to the codimensions
of these strata plus 1) in $\overline{K(m,n)}$ (this stratum $\sigma_{tot}$
belongs to the closure of $\sigma$ in $\overline{K(m,n)}$). The degrees of
the forms are arbitrary. Then we define a form $\omega$ on $\sigma_{tot}$ of degree
$\deg\omega=\sum\deg\omega_i+\sum\deg\omega_j^\prime$, as above. Then we
define
\begin{equation}\label{eq333.10}
\begin{aligned}
\ &L_{m+n}(\wtilde{\Psi}_1\wedge\dots\wedge\wtilde{\Psi}_n\wedge
\wtilde{\Theta}_1\wedge\dots\wedge\wtilde{\Theta}_m)=\\
&(\Alt ((\Psi_1\otimes\dots\otimes \Psi_n)\circledcirc
(\Theta_1\otimes\dots\otimes\Theta_m))\otimes[\omega]_{\dim\sigma_{tot}}
\end{aligned}
\end{equation}
In the case when the forms $\omega_i$'s are defined on {\it different}
strata, or $\omega_j^\prime$'s are defined on different strata, we define $L_{m+n}(\wtilde{\Psi}_1\wedge\dots\wedge\wtilde{\Psi}_n\wedge
\wtilde{\Theta}_1\wedge\dots\wedge\wtilde{\Theta}_m)$ as 0.
\subsubsection{\tt The $L_\infty$ structure}
\subsubsection{\tt The Gerstenhaber-Schack differential}
\subsubsection{\tt The computation of the Gerstenhaber-Schack cohomology for
$A=S(V^*)$}
\endcomment
\section{\tt The formality equation and the quantization of Lie
bialgebras}
\subsection{\tt The first look at the formality equation}
The components $\U_k$ of the analog of Kontsevich $L_\infty$ morphism is a
sum over graphs. To each graph we attach two things: these are a map from
$\Hom(A^{\otimes m},A^{\otimes n})$, and a top degree differential form on a
configuration space. Moreover, this top degree differential form is the
product of the differential forms attached to the edges of the graph (we do
not assume that these are necessarily 1-forms). Let us try to outline the
situation, starting with these very general principles.

We want to construct a map
\begin{equation}\label{eq22.1}
\U_k\colon\wedge^k(\g_1)\to\g_2[1-k]
\end{equation}
where
as above $\g_1=\oplus_{m,n\ge 1}\Hom(\wedge^m(V),\wedge^n(V))[-m-n+2]$, and
$\g_2=\oplus_{m,n\ge 1}\Hom(A^{\otimes m},A^{\otimes n})[-m-n+2]$. We know
that $\g_1$ is the Poisson Lie algebra, and expect for an
$L_\infty$-structure on $\g_2$ for which $U$ is an $L_\infty$ map. This
explains the shift of degree on $1-k$ in (\ref{eq22.1}).

Let us compute the degrees. We compute $\U_k(\gamma_1,\dots,\gamma_k)$.
Suppose that $\gamma_i\in\Hom(\wedge^{a_i}V,\wedge^{b_i}V)$.
Suppose $\U_k(\gamma_1,\dots,\gamma_k)\in\Hom(A^{\otimes m},A^{\otimes n})$.
Then it follows from (\ref{eq22.1}) that
\begin{equation}\label{eq22.2}
\sum_{i=1}^k(a_i+b_i-2)+(1-k)=m+n-2
\end{equation}
The last equation is equivalent to
\begin{equation}\label{eq22.3}
\sum_{i=1}^ka_i+\sum_{i=1}^kb_i=m+n+3k-3
\end{equation}
Now we can express the left hand side trough the number of edges of the
graph. Namely, we distinguish the "inner" edges and the "external" edges.
We suppose that we have some boundary in the configuration space with $m+n$
points there, and the interior, with $k$ points. There are no edges between
$m+n$ boundary points. Then each edge contain {\it at least} one inner
vertex. We call the edge internal, if the both ends of it are inner, and
inner, if only one among the two ends is inner. Then it is clear that
\begin{equation}\label{eq22.4}
2\sharp E_{inner}+\sharp E_{external}=\sum_{i=1}^ka_i+\sum_{i=1}^kb_i
\end{equation}
(Indeed, the number of edges starting and ending at the $i$th inner point is
$a_i+b_i$ by the assumption, we count then each inner edge twice and each
external edge one time).
Finally, we have:
\begin{equation}\label{eq22.5}
2\sharp E_{inner}+\sharp E_{external}=m+n+3k-3
\end{equation}
We interpret the last equation as follows: the dimension of the
configuration space should be equal to $3k+m+n-3$, and we should attach a
2-form to each inner edge, and a 1-form to each external edge. The dimension
$3k+m+n-3$ of the configuration space means that the "boundary" points
should belong to some 1-dimensional space, inner points--to a 3-dimensional
space, and there is an action of a 3-dimensional group on the configuration
space.

We made this computation in the very beginning of this work together with
Maxim Kontsevich, and then he invented the spaces $K(m,n;s)$ as an
appropriate candidates.
\subsection{\tt The admissible graphs and the corresponding operators in $\Hom(A^{\otimes m},A^{\otimes n})$}
\subsubsection{\tt The admissible graphs}
We will integrate over the spaces $K(m,n;s)$ some differential forms
of the top degree, associated with {\it admissible graphs }. We associate with any
{\it inner } edge of this graph the Propagator 2-form constructed above,
and with any {\it external } edge the corresponding 1-form.

\begin{definition}
Admissible graph $\Gamma$ is an oriented graph with labels such that:
\begin{itemize}
\item[1)] the set of vertices $V_\Gamma$ is $\{1,\dots,s\} \sqcup
\{\underline{1},\dots,\underline{m}\}\sqcup \{\overline{1},\dots,\overline{n}\}$,
$3s+m+n\geq 3$, the vertices from the set $\{1,\dots,s\}$ are called
vertices of the first type, the vertices from the set $\{\underline{1},\dots,\underline{m}\}$
are called the lower vertices of the second type, and the vertices from
the set $\{\overline{1},\dots,\overline{n}\}$ are called the upper vertices of the second type,
\item[2)] every edge $(v_1,v_2)\in E_\Gamma$ starts at a vertex of the first type or at an
upper vertex of the second type, and ends at a vertex of the first type or at a lower vertex of second type,
$v_1\in\{1,\dots,s\}\sqcup \{\overline{1},\dots,\overline{n}\},\ v_2\in
\{1,\dots,s\}\sqcup \{\underline{1},\dots,\underline{m}\}$,
if both $v_1,v_2\in\{1,\dots,s\}$, the edge is called inner, other edges are called
external, there are no (external) edges between two vertices
of the second type,
\item[3)]  there are no simple loops, that is edges of the type $(v,v)$,
there are no multiple external edges, but there are can be multiple inner edges,
\item[4)] for every vertex $k$ of the first type the sets of edges
$$
{Star}(k)=\{(v_1,v_2)\in E_\Gamma |v_1=k\}
$$
and
$$
{End}(k)=\{(v_1,v_2)\in E_\Gamma |v_2=k\}
$$
are labeled by symbols $(s_k^1,\dots,s_k^{\sharp {Star}(k)})$ and
$(e_k^1,\dots,e_k^{\sharp {End}(k)})$.
\end{itemize}
\end{definition}
The simplest examples of admissible graphs are shown in Figure 8.
\begin{remark}
Notice that we do {\it not} fix the number of the edges of $\Gamma$ in this
definition.
\end{remark}
\sevafigc{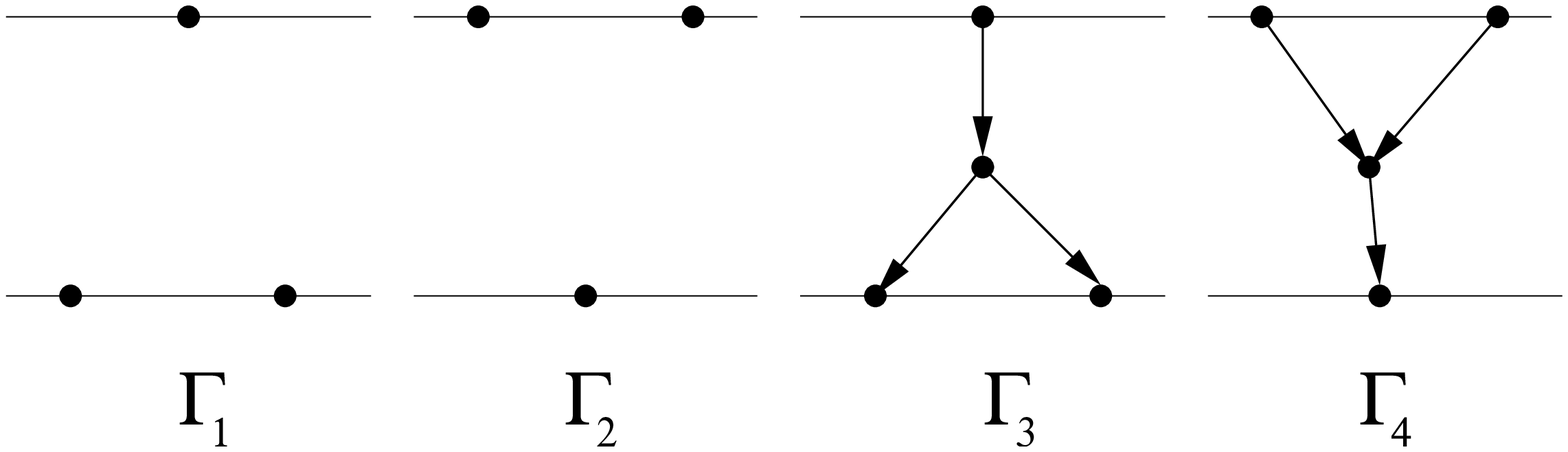}{100mm}{0}{The simplest admissible graphs}
\subsubsection{\tt The polydifferential operators associated with admissible
graphs}
For an admissible graph $\Gamma$ with $s$ vertices of first type, $(m,n)$
vertices of second type, we associate a map
$$
\Phi_\Gamma\colon{\otimes}^s \Lambda^\mb (V\oplus V^*)\to \Hom(A^{\otimes m}
\to A^{\otimes n})[1-s]
$$
Recall that $V$ is a Lie bialgebra, $A=S^\mb(V^*)$ is a free commutative
cocommutative associative bialgebra, and
$$
\deg(\Hom(\wedge^i(V),\wedge^j(V))=
\deg(\Hom(A^{\otimes i},A^{\otimes j})=i+j-2
$$
Let $\gamma_1,\dots,\gamma_s\in \Hom(\Lambda^\mb (V),\Lambda^\mb (V))$.
Then $\Phi_\Gamma(\gamma_1,\dots,\gamma_s)$ is non-zero only if
$\gamma_i\in\Hom (\wedge^{\sharp {Star}(i)}V,\wedge^{\sharp {End}(i)}V)$.

We are going to write a formula for
$$
\Phi_\Gamma(\gamma_1,\dots,\gamma_s)(f_1,\dots,f_{m})\in A^{\otimes n}
$$
The formula is the sum over all labelings of the edges of $\Gamma$ by indices
running from $1$ to $d$, $d=\dim V$:
\begin{equation}\label{eq1.55}
\Phi_\Gamma=\sum_{I\colon E_\Gamma\to\{1,\dots,d\}}\Phi_{\Gamma,I}
\end{equation}
where each $\Phi_{\Gamma,I}$ is

\begin{equation}\label{eq1.66}
\Phi_{\Gamma,I}=\Delta^{n}\left(\prod_{\begin{subarray}{c}
v\in V_\Gamma, \\ v\neq \overline{i}\  \text{for some}\ i\end{subarray}}\Psi_v\right)
\cdot\left(\otimes_{i=1}^{n}\Psi_{\overline{i}}\right)
\end{equation}
The product $\cdot$ here means the product
$$
(a_1\otimes\dots\otimes a_{n})\cdot
(b_1\otimes\dots\otimes b_{n})=((a_1\cdot b_1)\otimes\dots\otimes(a_{n}\cdot
b_{n})).
$$

At each vertex $v$ of the first type, the function $\Psi_v$ is a constant:
it is the matrix element for $\gamma_v$ ,
\begin{equation}\label{eq1.77}
\Psi_v=<\gamma_v(x^{I(s_v^1)}\wedge\dots\wedge
x^{I(s_v^{\sharp {Star}v})}),x^{I(e_v^1)}\wedge\dots\wedge
x^{I(e_v^{\sharp {End}v})}>
\end{equation}

Now at each {\it lower } vertex of the second type, the function $\Psi_v$ is a partial derivative
of $f_v$:
\begin{equation}\label{eq1.88}
\Psi_v=\left(\prod_{e\in E_\Gamma,e=(*,v)}\partial_{I(e)}\right)f_v
\end{equation}
and the function associated with each vertex $v$ of {\it upper } first type is
\begin{equation}\label{eq1.99}
\Psi_v=\left(\prod_{e\in E_\Gamma, e=(v,*)}x^{I(e)}\right)
\end{equation}

Now the formula for the summand $\Phi_{\Gamma,I}$ is given by formula (\ref{eq1.66}).

The formula (\ref{eq1.66}) is equivariant with respect to the linear group
${GL}(V)$ because it uses only invariant operations.

\begin{example}
For the first graph $\Gamma_1$ drawn in Figure~8, the map
$\Phi_{\Gamma_1}\colon A^{\otimes 2}\to A$ is the product:
$$
\Phi_{\Gamma_1}(f_1\otimes f_2)=f_1\cdot f_2
$$
For the second graph $\Gamma_2$ in Figure~8, the map
$\Phi_{\Gamma_2}\colon A\to A^{\otimes 2}$ is the coproduct:
$$
\Phi_{\Gamma_2}(f)=\Delta (f)
$$
Consider the third graph $\Gamma_3$ in Figure~8. It defines a map
$
\Phi_{\Gamma_3}\colon A^{\otimes 2}\to A
$
when in the unique inner vertex is placed the structure map $c_{ij}^k\colon\wedge ^2V\to V$.
By the definition above,
$$
\Phi_{\Gamma_3}(f_1\otimes f_2)=\sum_{i,j,k=1}^dc_{ij}^kx_k\partial_i(f_1)\partial_j(f_2)
=\{f_1,f_2\}
$$
where $\{,\}$ here stands for the Kostant-Kirillov Poisson bracket on $V^*$.

For the fourth graph $\Gamma_4$ drawn in Figure~8, we have a map
$
\Phi_{\Gamma_4}\colon A\to A^{\otimes 2}$  where in the unique inner vertex is placed
the cobracket $d_i^{jk}\colon V\to\wedge^2V$.
By our definition,
$$
\Phi_{\Gamma_4}(f)=\sum_{i,j,k=1}^d d_i^{jk}\Delta(\partial_i f)\cdot (x_j\otimes x_k)
$$
which is equal to the Poisson cobracket of $f$.
\end{example}
\subsection{\tt The differential form $\Omega_\Gamma$ associated with an
admissible graph $\Gamma$}
Consider an admissible graph $\Gamma$ with $m$ lower vertices of the first
type, $n$ upper vertices of the first type, and $s$ inner vertices (of the
second type). Suppose it has $\sharp E_{inner}$ edges connecting the inner
vertices and $\sharp E_{external}$ edges one the end-points of which is a
vertex of the first type (upper or inner).
We attach to the $\ell$th inner vertex of $\Gamma$ an element $\gamma_\ell$
of the Poisson Lie algebra $\g^\mb_1=\oplus_{m,n\ge
1}\Hom(\wedge^m(V),\wedge^n(V))[-m-n+2]$ where $V$ is a finite-dimensional
vector space. If $\gamma_\ell\in\Hom(\wedge^{m_\ell}(V),\wedge^{n_\ell}(V))$
than at the inner vertex $\ell$ there are $n_\ell$ edges starting at this
vertex, and $m_\ell$ edges ending at this vertex.
Then we have:
\begin{equation}\label{eqfin_1}
\sum_{\ell=1}^s\deg\gamma_\ell=2\sharp E_{inner}+\sharp E_{external}-2s
\end{equation}
Now we attach to the graph a form $\Omega_\Gamma$ on the space $\overline{K(m,n;s)}$ as
follows:

Each edge (inner or external) $\alpha$ of $\Gamma$ defines a map
$\jmath_\alpha\colon \overline{K(m,n;s)}\to \overline{K(0,0;2)}$ to the
3-dimensional Eye (see Section 2.2). We have constructed the Propagator 2-form
$\phi$
(for inner edges), and the Propagator 1-form $\phi_0$ (for external edges).
Now define
\begin{equation}\label{eqfin_2}
\Omega_\Gamma=\bigwedge_{\alpha\in
E_{inner}}\jmath_\alpha^*(\phi)\wedge\bigwedge_{\alpha\in
E_{external}}\jmath_\alpha^*(\phi_0)
\end{equation}
which is a form (with singularities) of degree $2\sharp E_{inner}+\sharp
E_{external}$ on the space $\overline{K(m,n;s)}$. (The labeling of $\Gamma$ allows to fix
the ordering in the wedge product). The singularities are only
$\delta$-singularities, and the form absolutely converges.
\begin{remark}
More precisely, we consider an object like the product of $\delta$-functions
which does not rigorously exist. For this we define all constructions above
with the Propagator and the configuration spaces defined from the {\it
truncated} 3-dimensional Eye (see Section 2.2.2) and then take the limit
when the rectangle $P_\lambda$ tends to the interval $(\lambda=0)$. This is
a kind of regularization we use here. In the sequel we always have in mind
this regularization and never consider these problems.
\end{remark}
In Section 3.2.2 we attached also an operator
$\Phi_\Gamma(\gamma_1,\dots,\gamma_s)\in\Hom(A^{\otimes m},A^{\otimes n})$
where $A=S(V^*)$.

Define now a map
\begin{equation}\label{eqfin_3}
\U_\Gamma\colon \wedge^s(\g^\mb_1)\to\wtilde{K^\mb}_{GS}[1-s]
\end{equation}
defined as
\begin{equation}\label{eqfin_4}
\U_\Gamma(\gamma_1\wedge\dots\wedge\gamma_s)=\frac1{s!}\Alt_{\gamma}\Phi_\Gamma(\gamma_1,\dots,\gamma_s)\otimes
\overline{p}_*{\Omega_\Gamma}
\end{equation}
Here $\overline{p}\colon \overline{K(m,n;s)}\to \overline{K(m,n)}$ is the natural projection, and $\overline{p}_*(\Omega_\Gamma)$ is considered as an element of
$B^{\mb\mb}(m,n)$ which has the only one nonzero component--the form
$\overline{p}_*(\Omega_\Gamma)$ on the top dimensional open stratum in
$\overline{K(m,n)}$.

Let us prove that the shift of grading is correct: indeed, by
(\ref{eqfin_1}), $\sum_\ell\deg\gamma_\ell=2\sharp E_{inner}+\sharp
E_{external}-2s$. One needs to prove that $2\sharp E_{inner}+\sharp
E_{external}-2s+(1-s)=\deg
\Phi_\Gamma+(\deg\overline{p}_*(\Omega_\Gamma)-(m+n-3))$. But
$\deg\Phi_\Gamma=m+n-2$ and $\deg\overline{p}_*(\Omega_\Gamma)=2\sharp E_{inner}+\sharp
E_{external}-3s$. We are done.

\subsection{\tt The "formality" equation: a Conjecture and a Theorem}
Fix $s\ge 0$. Define
\begin{equation}\label{eqfin_5}
\U_{s}(\gamma_1,\dots,\gamma_s)=\sum_{m,n,m+n\ge 3}\sum_{\Gamma\in\Gamma_{m,n;s}}\U_\Gamma(\gamma_1,\dots,\gamma_s)
\end{equation}
where $\Gamma_{m,n;s}$ are the graphs $\Gamma$ with $m$ lower vertices of
the first type, $n$ upper vertices of the first type, and $s$ inner
vertices.
When $\gamma_1,\dots,\gamma_s$ are fixed, the number $2\sharp
E_{inner}+\sharp E_{external}$ is fixed, and therefore the degrees $\deg
\Omega_\Gamma$ and $\deg\overline{p}_*(\Omega_\Gamma)$ are fixed. Consider
such $m,n$ that
$\deg\overline{p}_*\Omega_\Gamma=m+n-4$ (for 1 less than the
maximal).

The form $\overline{p}_*\Omega_\Gamma$ is not closed: by Lemma 2.4,
\begin{equation}\label{eqfin_6}
d\overline{p}_*\Omega_\Gamma=\overline{p}_*(d\Omega_\Gamma)-\overline{p}|_{\partial F}\Omega_\Gamma
\end{equation}
The first summand is 0 because the Propagator is closed. Finally, we have:
\begin{equation}\label{eqfin_7}
d\overline{p}_*\Omega_\Gamma+\overline{p}|_{\partial F}\Omega_\Gamma=0
\end{equation}
Both forms are top degree forms on $\overline{K(m,n)}$. Let us integrate the
equation (\ref{eqfin_7}) over the fundamental cycle in $\overline{K(m,n)}$.
We have:
\begin{equation}\label{eqfin_8}
\int_{\overline{K(m,n)}}d\overline{p}_*\Omega_\Gamma+\int_{\overline{K(m,n)}}\overline{p}|_{\partial F*}\Omega_\Gamma=0
\end{equation}
Consider the equation:
\begin{equation}\label{eqfin_9}
\sum_{\Gamma\in\Gamma_{m,n;s}^{m+n+3s-4}}\Phi_{\Gamma}(\gamma_1,\dots,\gamma_s)\cdot(\int_{\overline{K(m,n)}}d\overline{p}_*\Omega_\Gamma+\int_{\overline{K(m,n)}}\overline{p}|_{\partial F*}\Omega_\Gamma)=0
\end{equation}
where $\Gamma_{m,n;s}^{m+n+3s-4}$ denotes the graphs with $m$ lower vertices of
first type, $n$ upper vertices of first type, $s$ vertices of second type,
and for which $2\sharp E_{inner}+\sharp E_{external}=m+n+3s-4$, and $F$ is
the fiber.

By the usual Stokes formula, we can rewrite:
\begin{equation}\label{eqfin_10}
\int_{\overline{K(m,n)}}d\overline{p}_*\Omega_\Gamma=\int_{\partial\overline{K(m,n)}}\overline{p}_*\Omega_\Gamma
\end{equation}
It is clear that only boundary strata of codimension 1 will contribute to
these formulas. These strata have been described above.
Now we prove the following theorem:
\newpage
\begin{theorem*}
For fixed $s\ge 1$ we have the following relation for each $m,n$:
\begin{equation}\label{eqfinfin}
\begin{aligned}
\ &\sum_{\codim\sigma=1}\int_{\sigma}N_{\{w_i\}}\sum_{\{i_k\}\in\Sigma_s}\sum_{\{w_j\}}
\pm\left(\frac{\U_{w_1}(\gamma_{i_1},\dots,\gamma_{i_{w_1}})
\dots\U_{w_\ell-w_{\ell-1}}(\gamma_{i_{w_{\ell-1}}+1},\dots,\gamma_{i_{w_\ell}})}
{\U_{w_{\ell+1}-w_\ell}(\gamma_{i_{w_{\ell}+1}},\dots,\gamma_{i_{w_{\ell+1}}})
\dots
\U_{s-w_{\ell+\ell^\prime-1}}(\gamma_{i_{w_{\ell+\ell^\prime-1}+1}},\dots,\gamma_{i_s})}\right)_{\sigma}\\
&+\int_{\overline{K(m,n)}}\sum_{1\le i<j\le
s}\pm\U_{s-1}(\{\gamma_i,\gamma_j\},\gamma_1,\dots,\hat{\gamma}_i,\dots,\hat{\gamma}_j,\dots,\gamma_s)=0
\end{aligned}
\end{equation}
{\rm This is our "formality" theorem. Let us explain what is written here.
We take the sum over all strata $\sigma$ of codimension 1 in
$\overline{K(m,n)}$. We associate with each stratum of codimension 1 the corresponding operation
$\frac{\wtilde{\Psi}_1\wtilde{\Psi}_2\dots\wtilde{\Psi}_{n_1+1}}{\wtilde{\Theta}_1\wtilde{\Theta}_2\dots
\wtilde{\Theta}_{m_1+1}}$ where all $\wtilde{\Psi}$'s and $\wtilde{\Theta}$'s are our components $\U_\alpha$. We
write this operation now as $\left(\frac{\U_{w_1}(\gamma_{i_1},\dots,\gamma_{i_{w_1}})
\dots\U_{w_\ell-w_{\ell-1}}(\gamma_{i_{w_{\ell-1}}+1},\dots,\gamma_{i_{w_\ell}})}
{\U_{w_{\ell+1}-w_\ell}(\gamma_{i_{w_{\ell}+1}},\dots,\gamma_{i_{w_{\ell+1}}})
\dots
\U_{s-w_{\ell+\ell^\prime-1}}(\gamma_{i_{w_{\ell+\ell^\prime-1}+1}},\dots,\gamma_{i_s})}\right)_{\sigma}$.
The numbers $\ell$ and $\ell^\prime$ are uniquely defined from the combinatorics of $\sigma$.
The numbers $\{w_i\}$ are not defined uniquely, we take the summation over all possibilities. Finally,
we alternate over $\gamma_i$'s, taking the sum over all permutations from the permutation group $\Sigma_s$.
{\it The idea is that in the components $\U_i$ the corresponding factor in $B^{\mb\mb}$ could be not the top
degree form, we take the operation on them, and in the result (when we integrate) only the top degree
component will contribute} (see also the definition of the Integral map in Section 3.4.4).{\it Therefore, our condition on $s,m,n$ and the graphs
$\Gamma$ is that $\Gamma\in\Gamma_{m,n;s}^{m+n+3s-4}$. It means that in (\ref{eqfinfin}) we take the sum over all
such graphs.} The combinatorial factor $N_{\{w_i\}}$ is $\frac1{\prod_{i=1}^{\ell+\ell^\prime}w_i!}$.}
\end{theorem*}
We prove the Theorem in Sections 3.4.2-3.4.3 below. First of all, discuss that the
Gerstenhaber-Schack differential [GS] is hidden in the first summand of the
equation (\ref{eqfinfin}).

\subsubsection{\tt The "formality" relation and the Gerstenhaber-Schack
differential}
First of all, recall what the Gerstenhaber-Schack differential of a
(co)associative bialgebra $A$ is. Recall, that it is a differential on the
Gerstenhaber-Schack complex
\begin{equation}\label{eqeqeq_1}
K_{GS}(A)=\oplus_{m,n\ge 1}\Hom(A^{\otimes m}, A^{\otimes n})[-m-n+2]
\end{equation}
Now let $\Psi\colon A^{\otimes m}\to A^{\otimes n}\in K_{GS}^{m+n-2}(A)$.
We are going to define the Gerstenhaber-Schack differential
$d_\GS (\Psi)\in \Hom (A^{\otimes (m+1)},A^{\otimes n})\oplus\Hom (A^{\otimes m},A^{\otimes (n+1)})$. Denote the projections of $d_\GS$ to the first summand by
$(d_\GS)_1$, and the projection to the second summand by $(d_\GS)_2$.
The formulas for $(d_\GS)_1$ and $(d_\GS)_2$ are:
\begin{equation}\label{eq0.6fin}
\begin{aligned}
\ &(d_\GS)_1(\Psi)(a_0\otimes\dots\otimes a_m)=\ \ \ \ \ \ \ \ \ \ \ \ \ \ \ \ \ \ \ \ \ \ \hspace{7cm}\\
&\Delta^n(a_0)*\Psi(a_1\otimes\dots\otimes a_m)\\
&+\sum_{i=0}^{m-1}(-1)^{i+1}\Psi(a_0\otimes\dots\otimes (a_i* a_{i+1})\otimes
\dots\otimes a_m)\\
&(-1)^{m-1}\Psi(a_0\otimes\dots\otimes a_{m-1})*\Delta^n(a_m)
\end{aligned}
\end{equation}
and
\begin{equation}\label{eq0.7fin}
\begin{aligned}
\ &(d_\GS)_2(\Psi)(a_1\otimes\dots\otimes a_m)=\\
&(\Delta^{(1)}(a_1)*\Delta^{(1)}(a_2)*\dots *\Delta^{(1)}(a_m))\otimes
\Psi(\Delta^{(2)}(a_1)\otimes\dots\otimes \Delta^{(2)}(a_m))\\
&+\sum_{i=1}^n(-1)^i\Delta_{i}\Psi(a_1\otimes\dots\otimes a_m)\\
&+(-1)^{n+1}\Psi(\Delta^{(1)}(a_1)\otimes\Delta^{(1)}(a_2)\otimes\dots \otimes\Delta^{(1)}(a_m))\otimes(\Delta^{(2)}(a_1)*\Delta^{(2)}(a_2)*\dots *\Delta^{(2)}(a_m))
\end{aligned}
\end{equation}
Here we symbolically write $\Delta(a)=\Delta^{(1)}(a)\otimes\Delta^{(2)}(a)$ having in mind the sum of several such
terms, $\Delta(a)=\sum_i\Delta^{(1)}_i(a)\otimes\Delta^{(2)}_i(a)$, where
$\Delta$ is the coproduct in $A$, $*$ is the product in $A$, and
$\Delta_i=\Id\otimes\dots\otimes\Id\otimes\Delta\otimes\Id\otimes\dots\otimes\Id$
where $\Delta$ is applied to the $i$-th factor.

Now we show that the all terms in formulas (\ref{eq0.6fin}) and
(\ref{eq0.7fin}) are corresponded to some boundary strata of codimension 1
in $\overline{K(m+1,n)}$ or $\overline{K(m,n+1)}$.

\sevafigc{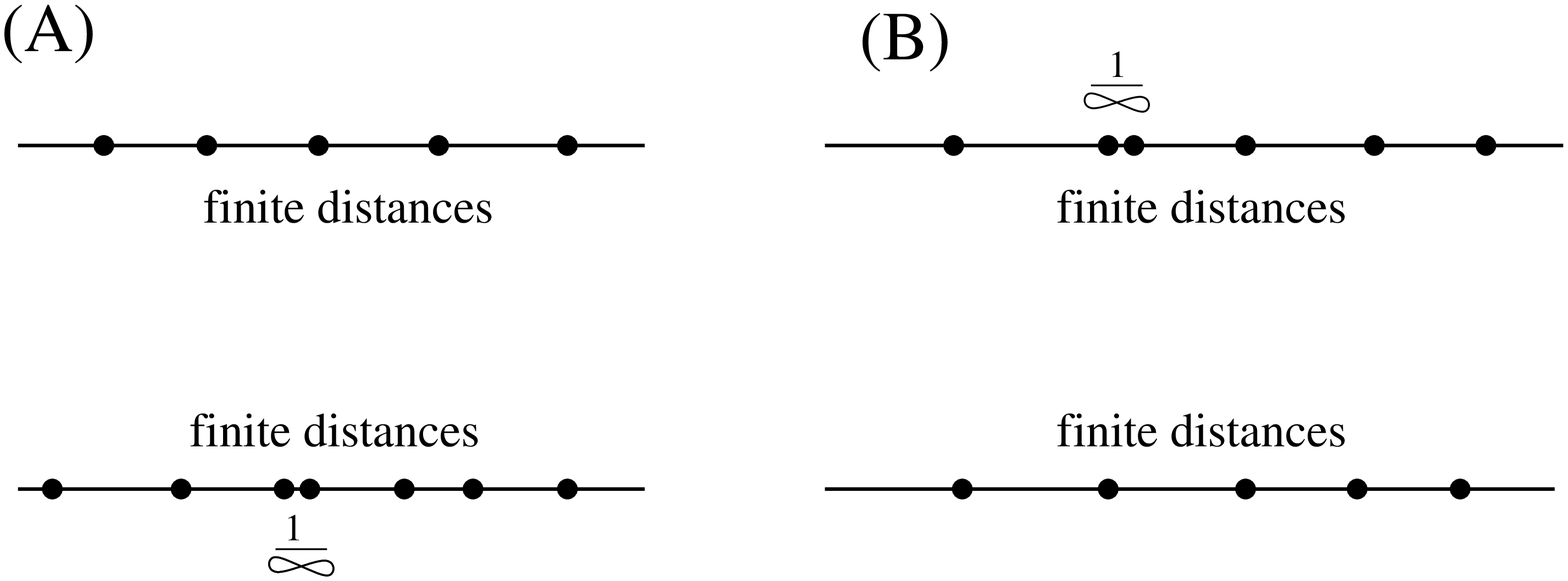}{140mm}{0}{Non-boundary terms in the
Gerstenhaber-Schack differential}

The non-boundary terms (the second lines in the right hand-sides of the
formulas for the Gerstenhaber-Schack differential) are corresponded to the
strata drawn in Figure~9(A) for $(d_{GS})_1$ and in Figure~9(B) for
$(d_{GS})_2$. (See Figure 2 above for the picture of a general stratum of
codimension 1). In the Figure~9(A) the stratum is $K(m,n)\times K(2,1)$.
The space $K(2,1)$ is a point, and we take the function (zero-form) 1 as the
corresponding element in $B^{\mb\mb}(2,1)$. Analogously for Figure~9(B).

The boundary terms (the first and the third lines in the r.h.s. of the
formulas) of the Gerstenhaber-Schack differential are drawn in Figure 10
(for $(d_{GS})_1$) and Figure 11 (for $(d_{GS})_2)$ below.
\sevafigc{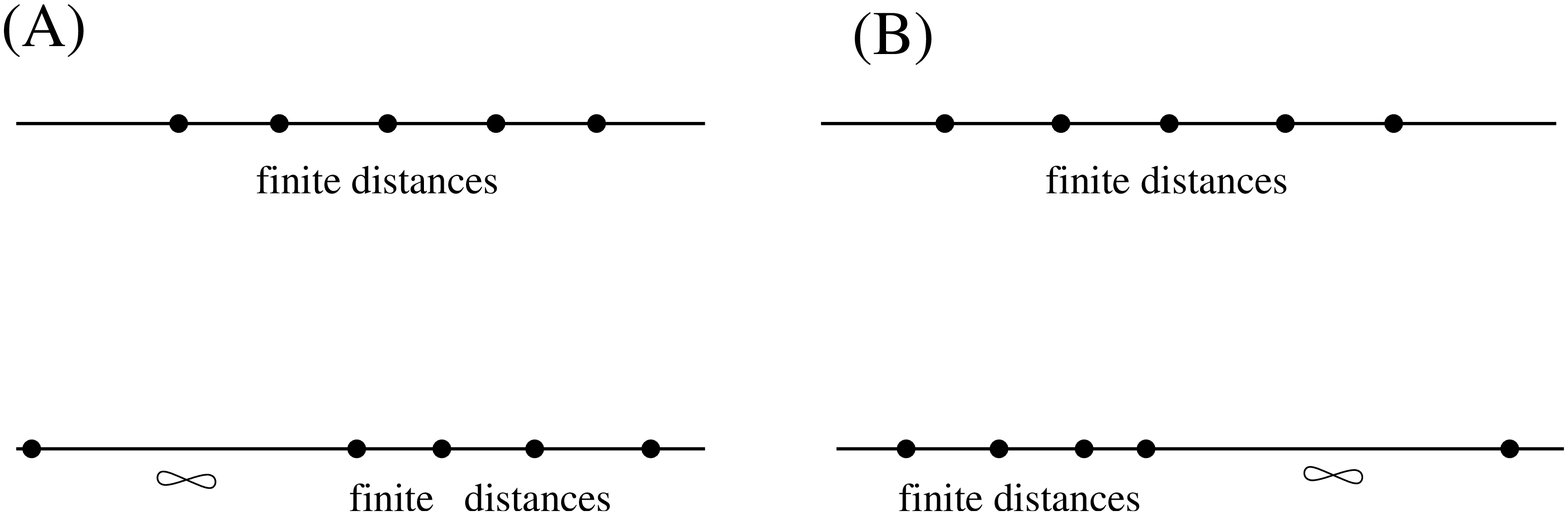}{140mm}{0}{The two boundary terms in $(d_{GS})_1$}
\sevafigc{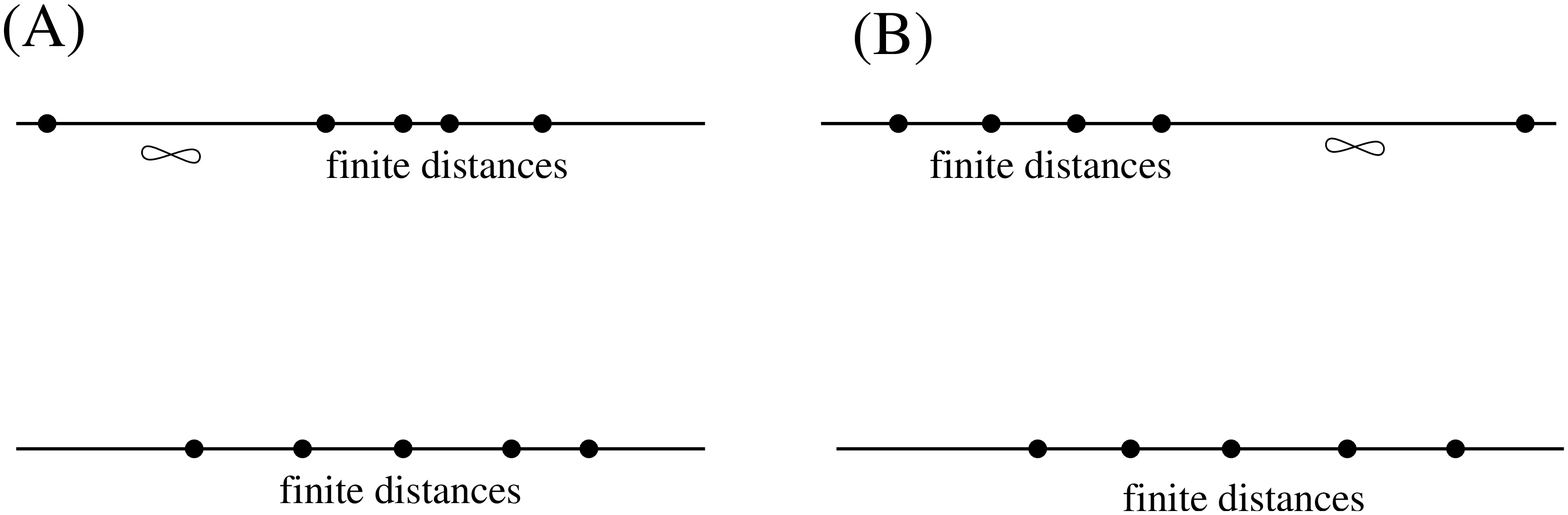}{140mm}{0}{The two boundary terms in $(d_{GS})_2$}
In the Figure~10 (A),(B) it is shown the two boundary terms (the first and
the third lines, correspondingly) of the formula (\ref{eq0.6fin}). The
boundary stratum for Figure~10(A) is $K_{1,m}^n\times K_2^{1,1,\dots,1}$.
The second factor is 0-dimensional. We need to construct a form on
$K_{1,m}^n$ to define the composition with values in $B^{\mb\mb}$. The
canonical map $K_{1,m}^n\to K(m,n)$ is an isomorphism, and we just take the
pull-back of the form $\omega$ on $K(m,n)$ where
$\wtilde{\Psi}=\Psi\otimes\omega$ and we compute
$(d_{GS})_1(\wtilde{\Psi})$. Analogously for the Figure~10(B), and for the
Figure~11(A),(B).

We see that we can rewrite some of the summands in the first line of
(\ref{eqfinfin}) to distinguish the Gerstenhaber-Schack differential among
these terms.
Roughly speaking, these are exactly the strata of codimension 1 such that
they are the products of two factors one of which is a point.

Now we are passing to the proof of the Formality Theorem.
\subsubsection{\tt We begin to prove the Formality Theorem}
It follows from formulas (\ref{eqfin_9}) and(\ref{eqfin_10})
\begin{equation}\label{eqvesna1}
\sum_{\Gamma\in\Gamma_{m,n;s}^{m+n+3s-4}}\Phi_{\Gamma}(\gamma_1,\dots,\gamma_s)
\cdot(\int_{\partial\overline{K(m,n)}}\overline{p}_*\Omega_\Gamma+
\int_{\overline{K(m,n)}}\overline{p}|_{\partial F*}\Omega_\Gamma)=0
\end{equation}
for fixed $m,n,s,\gamma_1,\dots,\gamma_s$.

{\it We claim that the first summand of (\ref{eqvesna1}) is equal to
the first summand of (\ref{eqfinfin}), and the second summand of
(\ref{eqvesna1}) is equal to the second summand of (\ref{eqfinfin})}.
At first, by the dimensional reasons only the strata of codimension 1 do
contribute to (\ref{eqvesna1}). First concern on the first summand.

Consider a boundary stratum $\sigma$ in $\overline{K(m,n;s)}$ of the type S2
(see Figure~7). There are $\gamma_i$'s in the inner points. We should prove
that the corresponding graphs are of a special form which give exactly
\begin{equation}\label{eqvesna2}
\int_{\sigma}\frac1{s!}\sum_{\{i_k\}\in\Sigma_s}\sum_{\{w_j\}}
\pm\left(\frac{\U_{w_1}(\gamma_{i_1},\dots,\gamma_{i_{w_1}})
\dots\U_{w_\ell-w_{\ell-1}}(\gamma_{i_{w_{\ell-1}}+1},\dots,\gamma_{i_{w_\ell}})}
{\U_{w_{\ell+1}-w_\ell}(\gamma_{i_{w_{\ell}+1}},\dots,\gamma_{i_{w_{\ell+1}}})
\dots
\U_{s-w_{\ell+\ell^\prime-1}}(\gamma_{i_{w_{\ell+\ell^\prime-1}+1}},\dots,\gamma_{i_s})}\right)_{\sigma}
\end{equation}
We mean that the differential forms corresponded to other graphs give 0. Let
us first draw the graphs which give (\ref{eqvesna2}).

Consider for simplicity the case of $\overline{K(2,2;s)}$. Schematically,
the picture is shown in Figure~12 below. Here in the picture the distances
between the points on the boundary lines are (a bit informally speaking)
infinite, and the inner points are close to the 4 external points. So, the
inner points form the 4 groups (relatively to which point on the boundary it
is closed). Between points of each group there is a "inner life", it means
that there are some edges between the points inside each group. We show by
the thin lines the only edges (oriented) between the points in different
groups. They are all oriented from the top to the bottom.
This picture is very good to show what kind of edges we have, but it is
little bit informal. When we are interesting which distances are $\infty$
and which are finite, Figure~13 below is better.
\sevafigc{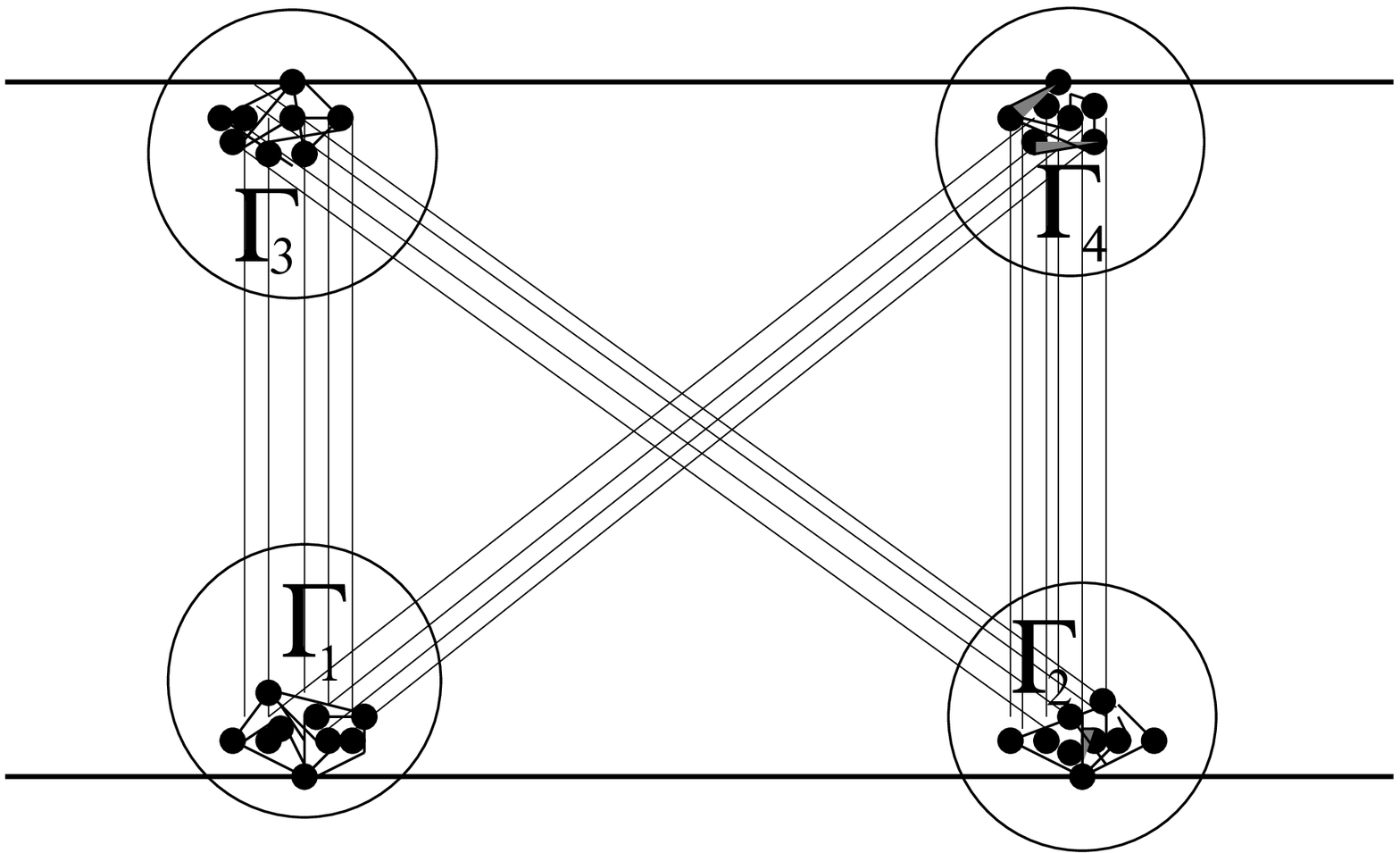}{120mm}{0}{A typical graph which contribute to the boundary
strata of the second type}
\sevafigc{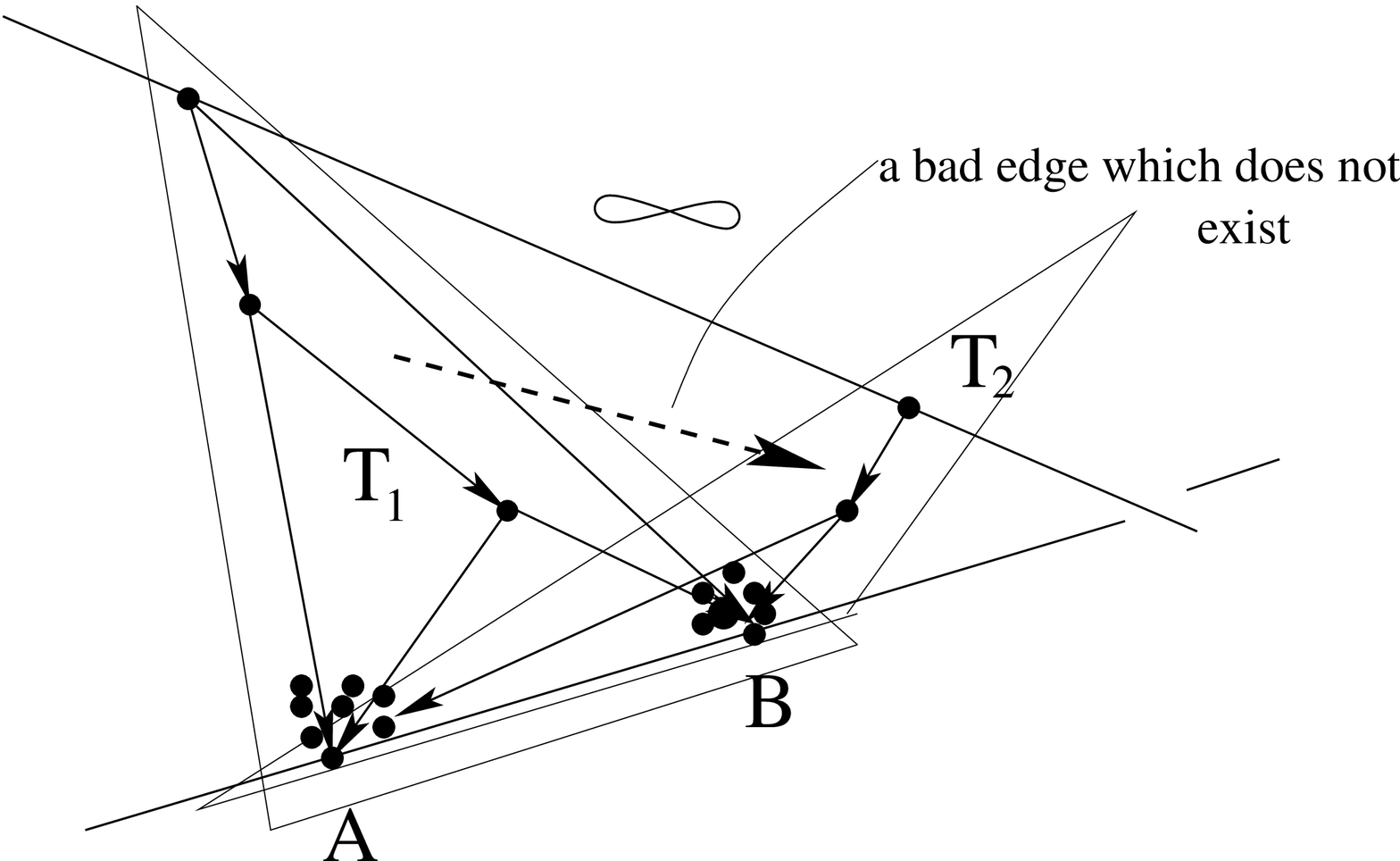}{120mm}{0}{A better picture for a stratum of the type
S2}
Here the two points ${\bf A}$ and ${\bf B}$ are in a finite distance from
each other, and the distance between the two upper boundary points is
$\infty$. The inner points at a finite domain are divided into 2 groups (the
triangles ${\bf T_1}$ and ${\bf T_2}$ in Figure~13). {\bf There are no edges
from a point of a one triangle to a point of another}. There are also inner
points infinitely close to the boundary points ${\bf A}$ and ${\bf B}$ (the
distances between them are of order $\frac1\infty$). There are some edges
between these points. When we apply the infinite transform $(0,0,\infty)\in
G^{(3)}$ to Figure~13, the configuration will be "symmetric": the upper line
will look as the lower line before the transform, and wise versa. The
picture in Figure~12 shows better what kind of edges we have.

Compute the weight $\int_{\sigma}\overline{p}_*\Omega_\Gamma$ for $\sigma$
and $\Gamma$ shown in Figures~12,13. {\it It is clear that the integral
factorizes into the product of 4 integrals}. The graph $\Gamma$ has
$m+n+3s-4$ edges (here in the example $m=n=2$), and each factor in the product is the weight of a graph
$\Gamma_i$ with $m_i+n_i+3s_i-3$ edges ($i=1\dots 4$).
These graphs $\Gamma_i$ are admissible graphs, which contribute to $\U_i$ in
formula (\ref{eqvesna2}). Now we should prove that the corresponding
operator $\Phi_\Gamma$ is equal to $\frac{\Psi_1\Psi_2}{\Theta_1\Theta_2}$
where $\Psi_1$ is $\Phi_{\Gamma_1}$ where $\Gamma_1$ is the graph inside the
triangle ${\bf T_1}$ in Figure~13, $\Psi_2=\Phi_{\Gamma_2}$ where $\Gamma_2$
is the graph inside ${\bf T_2}$, and $\Theta_1$ and $\Theta_2$ are
$\Phi_{\Gamma_3}$ and $\Phi_{\Gamma_4}$ where to see $\Gamma_3$ and
$\Gamma_4$ we should first apply to Figure~13 the infinite transform
$(0,0,\infty)\in G^{(3)}$. It follows from the definition of $\Phi_\Gamma$.

Now we have the following lemma:
\begin{lemma}
Each graph $\Gamma$ which contributes to the first summand of
(\ref{eqvesna1}) factorizes as is shown in Figures~12,13 into "disjoint
union" of 4 graphs (that is, there are no "bad edges" shown in Figure~13).
\begin{proof}
Consider an edge connecting a point from the triangle ${\bf T_1}$ with a
point from the triangle ${\bf T_2}$ (see the "bad edge" on Figure~13). Consider the
corresponding point of the 3-dimensional Eye (see Figure~5) corresponding to
this oriented pair of points. It is clear that this point belongs to the
boundary of the tetrahedron, but not to the closed interval $(\lambda=0)$.
Then the value of the Propagator 2-form on this edge is 0 by the
construction of the Propagator in Section 2.2.2. All other graphs belong to
the "factorizable" pictures.
\end{proof}
\end{lemma}

We considered here the case when $m=n=2$. The general case is absolutely
analogous.

It remains to consider the second summand of (\ref{eqvesna1}).
In the second summand of (\ref{eqvesna1}) we have the boundary strata of
$\overline{K(mn;s)}$ of the first type (see Section 2.3). It is clear that
the strata S1.2 do not contribute to the integral: as in the lemma above,
the corresponding "limit" edges belong to the boundary of the Tetrahedron
but not to the interval $(\lambda=0)$, and the Propagator vanishes on them.
If there is no edge connecting the limit point to the finite configuration,
the integral vanishes by the dimensional reasons.

There remain the strata S1.1 when several inner points move close to each
other. We claim that only the case when two inner points move close to each
other and are connected by a single edge does contribute. It follows from
the "Kontsevich lemma in dimensions $\ge3$ which we need in dimension 3. It is true also in dimension 2 but the
proof (given in [K1]) is more complicated. In dimensions $\ge 3$ this lemma
was found by Kontsevich in his study of the Chern-Simons theory (see, e.g. [K4])
but we do not know any place where the proof is written. Below we reproduce
the original Kontsevich's proof.

\subsubsection{\tt Kontsevich Lemma in dimension $\ge 3$}
Consider the configuration space $\mathrm{Conf}_n(\mathbb{R^d})$ defined as follows:
\begin{equation}\label{eqvesna3}
\mathrm{Conf}_n(\mathbb{R}^d)=\{(p_1,\dots,p_n\in\mathbb{R}^d), p_i\ne p_j\
{for}\ i\ne j\}
\end{equation}
There is a $(d+1)$-dimensional group $G(d)$ acting on the space
$\mathrm{Conf}_n(\mathbb{R}^d)$. This is the group of all $d$ linearly
independent shifts and 1-parametric family of dilatations. Consider the
quotient
\begin{equation}\label{eqvesna4}
\mathrm{C}_n(\mathbb{R}^d)=\mathrm{Conf}_n(\mathbb{R}^d)/G(d)
\end{equation}
If $n\ge 2$ it is a smooth manifold of dimension $nd-(d+1)$. For any two
points $\{p_i,p_j\}$ there is the restriction map $t_{ij}\colon
\mathrm{C}_n(\mathbb{R}^d)\to\mathrm{C}_2(\mathbb{R}^d)\simeq S^{d-1}$.
Consider the volume form $\Omega$ on the sphere $S^{d-1}$. Denote the
pull-back $t_{ij}^*\Omega$ by $\Omega_{ij}$.

\begin{lemma}
Let $d\ge 3$.
Consider an (oriented) graph $\Gamma$ with $n$ vertices and $e$ edges such that
$e(d-1)=nd-d-1$. Associate with $\Gamma$ a top degree form $\Omega_\Gamma$
on $\mathrm{C}_n(\mathbb{R}^d)$:
\begin{equation}\label{eqvesna5}
\Omega_\Gamma=\bigwedge_{(ij)\in E_\Gamma}\Omega_{ij}
\end{equation}
(the order in the wedge product is irrelevant) where $E_\Gamma$ is the set
of the edges of $\Gamma$.
Then the integral $\int_{\mathrm{C}_n(\mathbb{R}^d)}\Omega_\Gamma$ always
converges and is nonzero only in the case $n=2, e=1$.
\begin{proof}
We prove that when $d\ge 3$ there exists at least 1 vertex of $\Gamma$ of
the valence $\le 2$ (the valence is defined for $\Gamma$ as for an
non-oriented graph). Indeed, suppose that the valences of the all vertexes
are $\ge 3$. Then the number of edges $\sharp E_\Gamma\ge \frac{3n}2$. By
our assumption, $e(d-1)=nd-d-1$, therefore, $nd-d-1\ge (\frac{3n}2)(d-1)$.
The last inequality has no solutions for $d\ge 3$.

Denote by $v$ a vertex of valence $\le 2$. Consider independently the 3
cases when it is 0,1,2 and prove in all cases that
$\int_{\mathrm{C}_n(\mathbb{R}^d)}\Omega_\Gamma=0$.

If the valence of $v$ is 0 then the integral is 0 by dimensional reasons
because it is an integral of a form of degree $nd-d-1$ over a space of
dimension $nd-2d-1$. Analogously, when the valence of $v$ is equal to 1, put
another endpoint of this edge to a fixed point (we can do it using the group
$G(d)$). Then the integral by $v$ is an integral of $(d-1)$-form by
$\mathbb{R}^d$ which gives 0.

Consider the case when the valence of $v$ is 2. We want to prove that
$\int_{\mathrm{C}_n(\mathbb{R}^d)}\Omega_\Gamma=0$.
\sevafigc{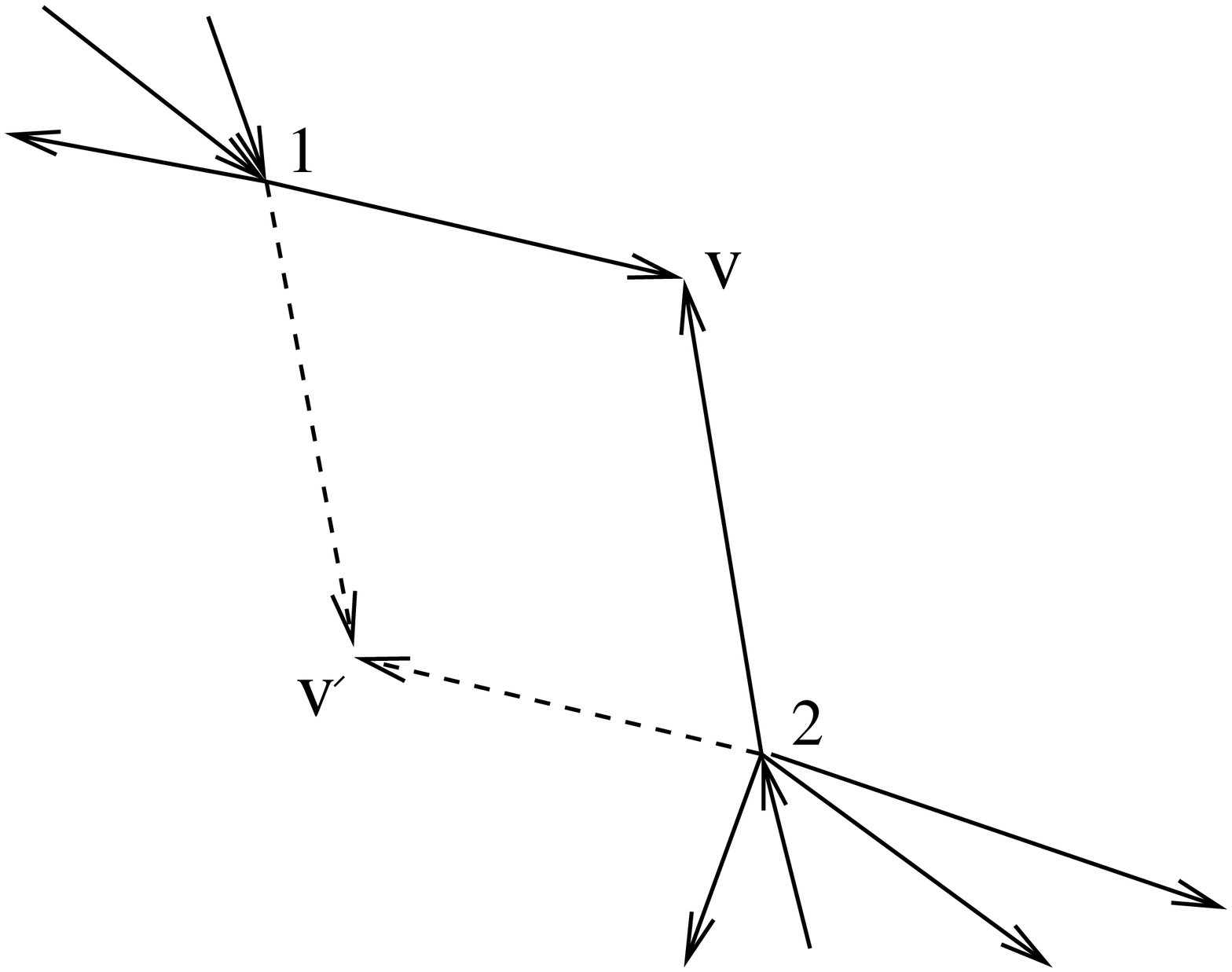}{100mm}{0}{The proof of the Kontsevich lemma}
Consider the symmetric integral over $v^\prime$ (see Figure~14).
The sign of the transform $\vartheta\colon v\mapsto v^\prime$ is $(-1)^d$.
Suppose the points {\bf 1} and {\bf 2} are fixed and we integrate only over
$v$. We have:
\begin{equation}\label{eqvesna6}
\int_{v\in\mathbb{R}^d}\Omega_{1v}\wedge\Omega_{2v}=(-1)^d\int_{v^\prime\in\mathbb{R}^d}
\Omega_{1v}\wedge\Omega_{2v}
\end{equation}
We have: $\Omega_{1v^\prime}=(-1)^{d}\Omega_{2v}$ and
$\Omega_{2v^\prime}=(-1)^{d}\Omega_{1v}$. Then we have:
\begin{equation}\label{eqvesna7}
\int_{v\in\mathbb{R}^d}\Omega_{1v}\wedge\Omega_{2v}=
(-1)^d\int_{v^\prime\in\mathbb{R}^d}\Omega_{2v^\prime}\wedge\Omega_{1v^\prime}
\end{equation}
But
$\Omega_{2v^\prime}\wedge\Omega_{1v^\prime}=(-1)^{(d-1)^2}\Omega_{1v^\prime}\wedge\Omega_{2v^\prime}$,
and the integral equal to itself multiplied on $(-1)^{d+(d-1)^2}$ which is
equal to $-1$ for all $d$. Therefore, the integral is 0.
\end{proof}
\end{lemma}

We finish to prove Theorem 3.4. In the second summand of (\ref{eqvesna1})
only the boundary strata of type S1.1 contribute. The integral factorizes to
the product of two integrals one of them is an  integral of the type of
Kontsevich lemma (in dimension 3). This lemma claims that the only case
which contributes to the integral is the case when two inner points
connected by a 1 edge move close to each other. It gives exactly the second
summands in (\ref{eqfinfin}).

Theorem 3.4 is proven.
\qed
\subsubsection{\tt The "Formality Conjecture"}
Consider the bicomplex $B^{\mb\mb}(m,n)$. As it was noticed before, it has
the only cohomology in degree 0 which is 1-dimensional. There is a very
natural map
\begin{equation}\label{eqfin10_1}
\int\colon B^{\mb\mb}(m,n)\to\mathbb{C}[0]
\end{equation}
which is a quasi-isomorphism. Namely, let
$\Omega=\omega_{m+n-3}+\omega_{m+n-4}+\dots+\omega_1+\omega_0$ be a general element
of $B^{\mb\mb}(m,n)[0]$ where $\omega_i$ is a linear combination of {\it top
degree} forms on strata of codimension $i$.
Define now the map $\int$ as
\begin{equation}\label{eqfin10_2}
\int(\Omega)=\int_{K(m,n)}\omega_{m+n-3}+\int_{\codim
\sigma=1}\omega_{m+n-4}+\dots+\int_{\codim\sigma=m+n-4}\omega_1+\int_{\codim\sigma=m+n-3}\omega_0
\end{equation}
and define the map $\int$ as 0 on $B^{\mb\mb}(m,n)[k]$ for $k\ne 0$. It
follows immediately from the Stokes formula that the map $\int$ vanishes on
the boundaries in $B^{\mb\mb}(m,n)$, and it defines a quasi-isomorphism.

It motivates the following "Formality Conjecture":
\begin{conjecture}
There exist the components
$\wtilde{\U}_\Gamma:=\Alt_{\gamma}\Phi_\Gamma(\gamma_1,\dots,\gamma_s)\otimes(
\overline{p}_*{\Omega_\Gamma}+D\wtilde{\Omega}_\Gamma)$ where $D$ is the
differential in $B^{\mb\mb}$ (that is, $\wtilde{\U}_\Gamma$ differs from
$\U_\Gamma$ on a boundary in $B^{\mb\mb}$) and the corresponding
$\wtilde{\U}_s=\sum_{m,n}\sum_{\Gamma\in\Gamma_{m,n;s}}\wtilde{\U}_\Gamma$
such that the following "formality on the level of complexes" holds (here we
do not fix $m,n$, only fix $s$):
\begin{equation}\label{eqfinfin1}
\begin{aligned}
\ &\sum_{\codim\sigma=1}N_{\{w_i\}}\sum_{\{i_k\}\in\Sigma_s}\sum_{\{w_j\}}
\pm\left(\frac{\wtilde{\U}_{w_1}(\gamma_{i_1},\dots,\gamma_{i_{w_1}})
\dots\wtilde{\U}_{w_\ell-w_{\ell-1}}(\gamma_{i_{w_{\ell-1}}+1},\dots,\gamma_{i_{w_\ell}})}
{\wtilde{\U}_{w_{\ell+1}-w_\ell}(\gamma_{i_{w_{\ell}+1}},\dots,\gamma_{i_{w_{\ell+1}}})
\dots
\wtilde{\U}_{s-w_{\ell+\ell^\prime-1}}(\gamma_{i_{w_{\ell+\ell^\prime-1}+1}},\dots,\gamma_{i_s})}\right)_{\sigma}\\
&+\sum_{1\le i<j\le
s}\pm\wtilde{\U}_{s-1}(\{\gamma_i,\gamma_j\},\gamma_1,\dots,\hat{\gamma}_i,\dots,\hat{\gamma}_j,\dots,\gamma_s)=0
\end{aligned}
\end{equation}
\end{conjecture}

The meaning of this Conjecture is that if we would "add" some other terms to
our structure depending on the strata of codimension 1 to $\wtilde{K}_{GS}$
to have a "pure" $L_\infty$ structure, the Conjecture above would express
that this structure is formal, that is, is equivalent to its cohomology
$\g_1^\mb=\oplus_{m,n\ge 1}\Hom(\wedge^m(V),\wedge^n(V))[-m-n+2]$. In the
Subsection below we show how the Lie bracket on $\g^\mb_1$ (which is a
Poisson Lie bracket) is hidden in the operations corresponding to the strata
of codimension 1 in $\overline{K(m,n)}$.
\subsubsection{\tt The strata of codimension 1 and the Lie bracket on
$\g_1^\mb$}
First of all, formulate the following result:
\begin{lemma}
\begin{itemize}
\item[(i)] For any (co)associative bialgebra, the Gerstenhaber-Schack
differential is indeed a differential, that is,
$((d_{GS})_1+(d_{GS})_2)^2=0$,
\item[(ii)] in the case of the free commutative cocommutative (co)associative
bialgebra $A=S(V^*)$ where $V$ is a finite-dimensional vector space, the
Gerstenhaber-Schack cohomology of $\g_2^\mb=K^\mb_{GS}$ is isomorphic to
$\g_1^\mb=\oplus_{m,n\ge 1}\Hom(\wedge^m(V),\wedge^n(V))[-m-n+2]$,
\item[(iii)] the following analog of the Hochschild-Kostant-Rosenberg map
$\varphi\colon\g_1^\mb\to\g^\mb_2$ (in the case $A=S(V^*)$) defines an
isomorphism on cohomology where the map $\varphi$ is the operation
$\Phi_\Gamma$ corresponding to the following graph $\Gamma$:
\sevafigc{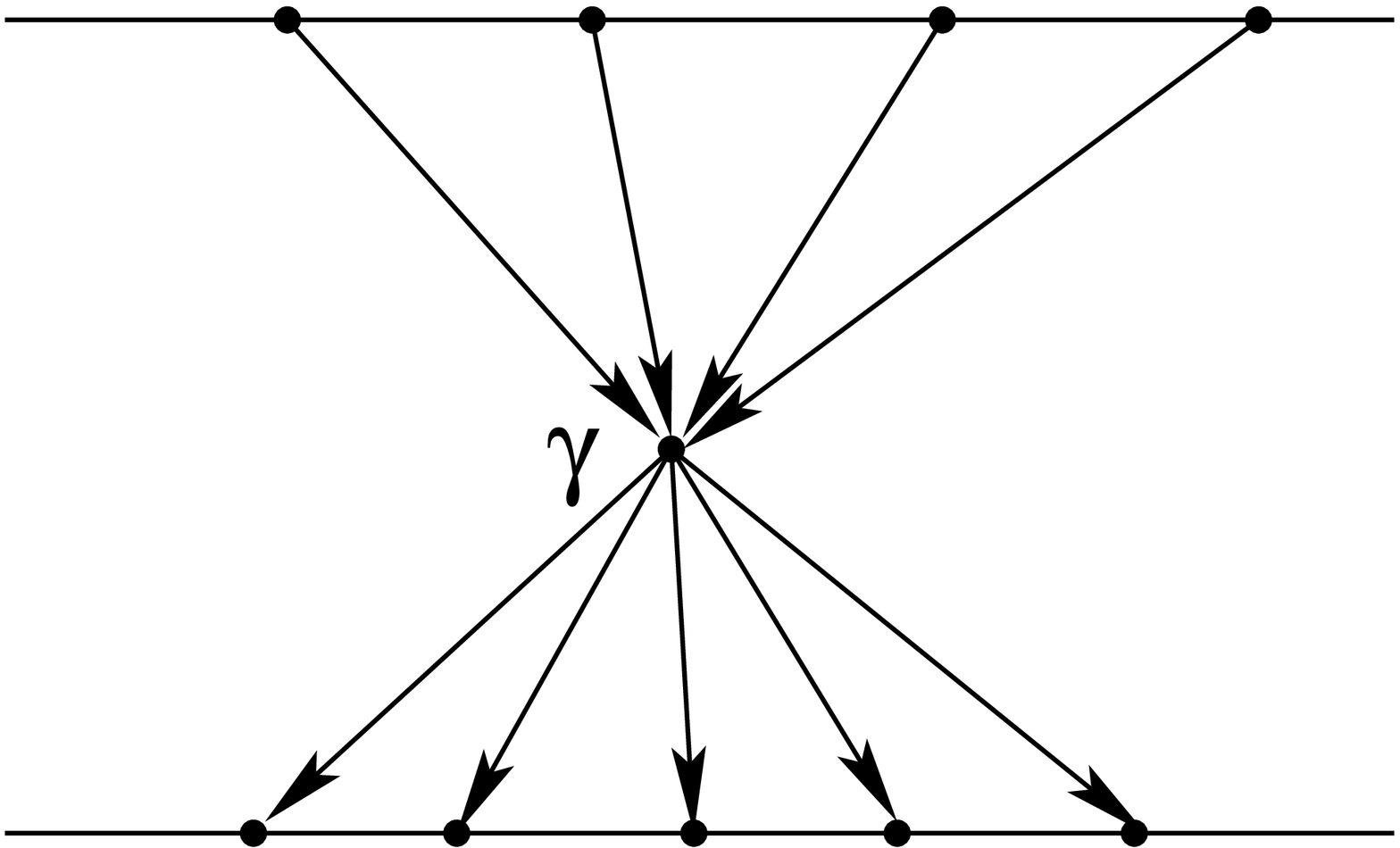}{80mm}{0}{The analog of the
Hochschild-Kostant-Rosenberg map}
\end{itemize}
\begin{proof}
We just sketch the proof here. We need to compute the Gerstenhaber-Schack
cohomology in the case of bialgebra $A=S(V^*)$.
The two components of the Gerstenhaber-Schack differential form a bicomplex.
We can use any of the two canonical spectral sequences to compute the
cohomology. Compute first the cohomology of the differential $(d_{GS})_1$.
We have the Hochschild cohomological complexes here $\Hoch^\mb(A,M_k)$ where
$M_k=A^{\otimes k}$, and the bimodule structure is given by the left and the
right multiplications on $\Delta^{k-1}(a)$. First we need to compute this
Hochschild cohomology (for all $k$). For this we need a more direct
description of the bimodule $M_k$.
For any Lie algebra $\g$ consider the left action of $\g$ on the universal
enveloping algebra $U(\g)$. Consider the tensor product of $k$ of such
$\g$-modules. Then we obtain the corresponding $U(\g)$-module structure on
$\U(\g)^{\otimes k}$. In the case of Abelian Lie algebra $\g=V^*$ this
structure on $S^(V^*)^{\otimes k}$ coincides with the left action on our
bimodule. The analogous construction works as well for the right action.
Therefore, we need to decompose the $\g$-(bi)module $U(\g)^{\otimes k}$, and
then we get automatically a decomposition of the corresponding
$U(\g)$-(bi)module.
It is easy to do in the case of an Abelian Lie algebra $\g$ (for $\g=V^*$).
Namely, $S(V^*)^{\otimes k}$ is the algebra of polynomials $\{f(v_1,\dots,
v_k)\}$ where $v_i\in V$. The action of $V^*$ is $f\mapsto
(\xi(v_1)+\dots+\xi(v_k))\cdot f$ for $\xi\in V^*$. Now polynomials
$f(v_1+\dots+v_k)$ form a submodule (with respect to the both left and right
actions)  which is isomorphic to the bimodule $A=S(V^*)$. It is clear then
that as an $A$-bimodule, $A^{\otimes k}=A\otimes A^{\otimes k-1}$ where in
the bimodule structure $A$ acts (tautologically) only on the first factor
(that is, it is a direct sum of infinitely many copies of $A$). We can
identify $A^{\otimes k-1}$ with $\{f(v_1-v_2,v_1-v_3,\dots,v_1-v_k)|\}$.

Therefore, by the Hochschild-Kostant-Rosenberg theorem, the cohomology of
the differential $(d_{GS})_1$ in the $k$-th row is $T^\mb_{poly}(V)\otimes
A^{\otimes k-1}$. Now we can compute the second differential on the image of
this cohomology in $K^\mb_{GS}$ with respect to the
Hochschild-Kostant-Rosenberg map, and then we get that the second term $E_2$ of the spectral sequence
is isomorphic to $\g_1^\mb$. We need to prove independently that under the
map of the item (iii) of this Lemma $\g_1^\mb$ is ${\it imbedded}$ to the
cohomology.
\end{proof}
\end{lemma}
Now the following question arises: which operations  on the (homotopical)
Gerstenhaber-Schack complex $\wtilde{K}_{GS}^\mb$ among the operations corresponded to the
strata of codimension 1 in $\overline{K(m,n)}$ define on the cohomology the Poisson Lie bracket on
$\g_1^\mb$? Here we answer this question.

Let $\Psi\in\Hom(A^{\otimes m_1+1},A^{\otimes n_0})$ and $\Theta\in \Hom(A^{\otimes m_0},
A^{\otimes n_1+1})$. Define $\wtilde{\Psi}=\Psi\otimes\omega_1$ and
$\wtilde{\Theta}=\Theta\otimes\omega_2$ where the form $\omega_1$ is a top
degree form on the top dimensional open stratum in
$\overline{K(m_1+1,n_0)}$, and the form $\omega_2$ is a top degree form on
the top dimensional open stratum in $\overline{K(m_0,n_1+1)}$. Define all
other data in the operation corresponded to the stratum of codimension 1 in
$\overline{K(m,n)}$ drawn in Figure~2 as $*^{m_1}\otimes 1$ and
$\Delta^{n_1}\otimes 1$ when the 1's above are considered as zero degree
differential forms on $\overline{K(m_1+1,1)}$ and on $\overline{K(1,n_1+1)}$
correspondingly. Then the composition
\begin{equation}\label{eqfin0022}
\frac{\wtilde{\Psi}_1\wtilde{\Psi}_2\dots\wtilde{\Psi}_{n_1+1}}{\wtilde{\Theta}_1\wtilde{\Theta}_2\dots
\wtilde{\Theta}_{m_1+1}}:=\frac{\Psi_1\Psi_2\dots\Psi_{n_1+1}}{\Theta_1\Theta_2\dots\Theta_{m_1+1}}\otimes[\omega]_\sigma
\end{equation}
and the form $\omega$ in the r.h.s is a {\it top degree form} on the stratum
$\sigma$ in $\overline{K(m,n)}$.

Even if we suppose that this operation when $\omega_1$ and $\omega_2$ are
not top degree forms is zero, it is well-defined (the latter means it is
compatible with the differential in $B^{\mb\mb}$).

We claim that on the cohomology (of both differentials in $B^{\mb\mb}$ and
of the Gerstenhaber-Schack differential) this operation defines exactly the
Poisson Lie bracket in $\g_1^\mb$ (when we take the sum over all $\sigma$
with fixed $m_0,m_1,n_0,n_1$). Moreover, our Hochschild-Kostant-Rosenberg
map is a "Lie algebra map" up to a boundary. This is an initial point for
the construction of the "formality" morphism.
\subsection{\tt Deformation quantization of Lie bialgebras}
Let $V$ be a Lie bialgebra. Recall that it means that we have
$\alpha\in\Hom(\wedge^2V,V)$ and $\beta\in\Hom(V,\wedge^2V)$ such that
$\{\alpha,\alpha\}=0$, $\{\beta,\beta\}=0$, and $\{\alpha,\beta\}=0$ (where
$\{,\}$ denotes the Poisson Lie bracket in $\g_1^\mb$).

Define a product and a coproduct on $S(V^*)$ as follows:
\begin{equation}\label{vesnavesna1}
f*g=\sum_{\ell_1,\ell_2\ge
0}\frac1{\ell_1!}\frac1{\ell_2!}\sum_{\Gamma\in\Gamma_{2,1;\ell_1+\ell_2}^{3(\ell_1+\ell_2)}}\hbar^{\ell_1}_1\hbar^{\ell_2}_2
\U_\Gamma(\alpha,\dots,\alpha,\beta,\dots,\beta)(f\otimes g)
\end{equation}
and
\begin{equation}\label{vesnavesna2}
\Delta_*(f)=\sum_{\ell_1,\ell_2\ge
0}\frac1{\ell_1!}\frac1{\ell_2!}\sum_{\Gamma\in\Gamma_{1,2;\ell_1+\ell_2}^{3(\ell_1+\ell_2)}}
\hbar_1^{\ell_1}\hbar_2^{\ell_2}\U_{\Gamma}(\alpha,\dots,\alpha,\beta,\dots,\beta)(f)
\end{equation}
Here in the right-hand sides of the formulas we have $\ell_1$ of $\alpha$'s
end $\ell_2$ of $\beta$'s. The values of $\U_\Gamma$ are top degree forms on
$\overline{K(2,1)}$ and $\overline{K(1,2)}$ which are just points, and we
identify the function on them with the numbers.
\begin{theorem*}
The product $f*g$ and the coproduct $\Delta_*(f)$ defined above satisfy the
axioms (i)-(iii) of (co)associative bialgebras (see the first page of the
Introduction).
\begin{proof}
To prove the associativity, apply the formality theorem to the space
$\overline{K(3,1;\ell_1+\ell_2)}$ and to
$\gamma_1,\dots,\gamma_{\ell_1}=\alpha$,
$\gamma_{\ell_1+1},\dots,\gamma_{\ell_1+\ell_2}=\beta$.

To prove the coassociativity, apply the formality theorem to the space
$\overline{K(1,3;\ell_1+\ell_2)}$ and to the $\gamma_i$'s as above.

To prove the compatibility, apply the formality theorem to the space
$\overline{K(2,2;\ell_1+\ell_2)}$ and to the $\gamma_i$'s as above.
\end{proof}
\end{theorem*}

Notice that we obtain a 2-parametric deformation quantization of Lie
bialgebras.
\subsection{\tt An informal discussion about why we have not an $L_\infty$
structure here}
Probably the most strangeness among the strange things in this paper is that
the right-hand sides of our "formality" do not obey the $L_\infty$ Jacobi
identity. Here we want to explain a geometrical cause for that.

Suppose that we have a collection of open manifolds $M_k$ and of their
compactifications $\overline{M}_k$ which all are manifold with corners such
that any face of codimension $\ell$ is a product of $\ell+1$ different
$M_i$'s (maybe these are not the most right conditions for the description
in general). Then when we compute the boundary operator $\partial M_k$ we
get a quadratic expression in $M_i$'s which is the sum over the all strata
of codimension 1. Then suppose we have a stratum of codimension 1 of $M_k$
which is a product $M_i\times M_j$. How to compute the boundary
$\partial(M_i\times M_j)$? The most natural is to suppose that the Leibniz rule
for $\partial$ holds. (Of course, we have in mind that the strata are
labeled by some combinatorial objects like trees, and the Leibniz rule
should be understood in this sense).
For example, in the case of Stasheff polyhedra the Leibniz rule holds.

Suppose we have something like an "algebraic representation" of the family
$\{M_i\}$. It means that we associate to each space $M_k$ a graded vector
space $V_k$ and to each stratum of codimension 1 an operation on
$\oplus_k V_k$. Suppose that these operations are compatible with the
boundary operations. The latter means that if a stratum $\sigma_1$ can be
obtained by a "degeneration of codimension 1" from a stratum $\sigma$, and a
stratum $\sigma_2$ can be obtained by a degeneration of codimension 1 from a
stratum $\sigma_1$, then it implies that the corresponding compositions can
be also obtained functorially one from another. {\it Then the Jacobi
identity is an immediate consequence of the Leibniz rule for the boundary
operator} $\partial$. This is more or less in the spirit of the Markl's
paper [M1]. It follows from the identity $\partial^2=0$ for the boundary
operator.

But we can imagine the situation when the Leibniz rule is not satisfied. For
instance, it could be an operator of the second order. It means that to
compute the boundary of any stratum we need to know the boundaries of strata
of codimension 0 and of strata of codimension 1. Then the Jacobi identity
will be replaced by a more complicated structure. Probably, it will be an
odd vector field $Q$ on $L[1]$ such that $Q^3=0$ (not $Q^2=0$ as in the
$L_\infty$-case).

In our situation the problem is that the Leibniz rule is not satisfied. Let
us explain how it works:

Consider the space $\overline{K(m,n)}$ constructed in the paper. We
canonically identify a boundary stratum $\sigma$ of codimension 1 with a product
of two spaces of different type, $K_{m_1+1}^{1,\dots,1, n_0,1,\dots,1}\times
K_{1,\dots,1,m_0,1,\dots,1}^{n_1+1}$. Then we compute the boundary of
$\sigma$ of codimension 1. Probably, for the spaces $K_{m_1+1}^{1,\dots,1,
n_0,1,\dots,1}$ and $K_{1,\dots,1,m_0,1,\dots,1}^{n_1+1}$
the Leibniz rule is satisfied. The map $i\colon K_{m_1+1}^{1,\dots,1, n_0,1,\dots,1}\times
K_{1,\dots,1,m_0,1,\dots,1}^{n_1+1}\to \overline{K(m,n)}$ is an imbedding. But the
canonical extension $\overline{K_{m_1+1}^{1,\dots,1, n_0,1,\dots,1}\times
K_{1,\dots,1,m_0,1,\dots,1}^{n_1+1}}\to\overline{K(m,n)}$ is {\it not} an
imbedding. In particular, it maps the boundary of codimension 1 to the
boundary of $\sigma$ of codimension 1 surjectively, but it contracts some components
to points. Therefore, the Leibniz rule is satisfied modulo these
components.

It is very interesting to understand the whole structure arising from this
construction.
\subsection*{Acknowledgements}
I am grateful to Maxim Kontsevich for sharing with me with his construction
of the spaces $K(m,n)$ and $K(m,n;s)$ and for many discussions. When I constructed the CROC compactification of
$K(m,n)$ in [Sh1], I discussed it many times with Borya Feigin, and these discussions
finally explained to me the necessity of a Stasheff-type compactification for
the usual (commutative) deformation theory of (co)associative bialgebras.
As well I discussed with Borya Feigin almost all important points of the paper,
and his remarks many times clarified the situation for me.
Discussions on the Stasheff-type compactification and on the corresponding Lie algebra
with Giovanni Felder were very useful
for me. I am also grateful to Giovanni Felder for very helpful discussions
on the Propagator.
I am gateful to Alberto Cattaneo for his explanation to me the Kontsevich
proof of Lemma 3.4.3, and to Martin Markl for the reading of a part of this
paper and for his remarks and corrections.
I would like to thank the ETH (Zurich) for the financial
support and for the very stimulating atmosphere.
\bigskip

P.S. The author afraids that some fractions with factorials in the formulas
concerning the "formality" are computed wrong :)

\bigskip
\bigskip
Dept. of Math., ETH-Zentrum, 8092 Zurich, SWITZERLAND\\
e-mail: {\tt borya@mccme.ru, borya@math.ethz.ch}

\end{document}